\documentclass[12pt]{amsart}

\usepackage{amsmath}
\usepackage{amsbsy}
\usepackage{amssymb}
\usepackage{euscript}
\usepackage{amscd}
\usepackage{eq}

\newcommand{\bba}{{\mathbb A}}

\newcommand{\bbc}{{\mathbb C}}

\newcommand{\bbq}{{\mathbb Q}}
\newcommand{\bbr}{{\mathbb R}}

\newcommand{\bbz}{{\mathbb Z}}

\newcommand{\al}{{\alpha}}

\newcommand{\gam}{{\gamma}}

\newcommand{\Gam}{{\Gamma}}
\newcommand{\del}{{\delta}}

\newcommand{\Del}{{\Delta}}

\newcommand{\vep}{{\varepsilon}}

\newcommand{\ka}{{\kappa}}

\newcommand{\lam}{{\lambda}}

\newcommand{\lamb}{\underline{\lambda}}

\newcommand{\sig}{{\sigma}}
\newcommand{\Sig}{\Sigma}

\newcommand{\om}{{\omega}}
\newcommand{\Om}{{\Omega}}
\newcommand{\hp}{\widehat{\Phi}}

\newcommand{\gc}{{\mathfrak c}}
\newcommand{\gC}{{\mathfrak C}}

\newcommand{\g }{{\mathfrak g}}

\newcommand{\gi}{{\mathfrak i}}

\newcommand{\gM}{{\mathfrak M}}

\newcommand{\gp}{{\mathfrak p}}

\newcommand{\gV}{{\mathfrak V}}

\newcommand{\gEx}{{\mathfrak {Ex}}}

\newcommand{\esN}{{\EuScript{N}}}
\newcommand{\esR}{{\EuScript{R}}}
\newcommand{\esL}{{\EuScript{L}}}

\newcommand{\cS}{{\mathcal S}}
\newcommand{\co}{{\mathcal O}}
\newcommand{\cD}{{\mathcal D}}

\newcommand{\xspan}{{\operatorname{Span}}}
\newcommand{\aff}{{\operatorname {Aff}}}

\newcommand{\tr}{{\operatorname {Tr}}}

\newcommand{\lev}{{\operatorname{lev}}}
\newcommand{\n}{{\operatorname {N}}}
\newcommand{\m}{{\operatorname {M}}}

\newcommand{\im}{{\operatorname {Im}}}
\newcommand{\re}{{\operatorname {Re}}}
\newcommand{\gal}{{\operatorname{Gal}}}
\newcommand{\gl}{{\operatorname{GL}}}

\newcommand{\orth}{{\operatorname{O}}}

\newcommand{\sst}{{\text{ss}}}

\newcommand{\vol}{{\operatorname{vol}}}
\newcommand{\ord}{{\operatorname{ord}}}

\newcommand{\rep}{representation}

\newcommand{\bs}{Schwartz--Bruhat function}

\newcommand{\pv}{prehomogeneous vector space}

\newcommand{\res}{\operatorname{Res}}

\newcommand{\supp}{\operatorname{supp}}

\newcommand{\A}{\bba}
\newcommand{\Z}{\bbz}
\newcommand{\Q}{\bbq}
\newcommand{\R}{\bbr}
\newcommand{\C}{\bbc}

\newcommand{\ma}{\bba^{\times}}
\newcommand{\mr}{\R_+}

\newcommand{\mk}{k^{\times}}

\newcommand{\rg}{G_{k}}

\newcommand{\rv}{V_{k}}

\newcommand{\md}{d^{\times}}
\newcommand{\gMf}{\gM_{\text f}}
\newcommand{\gMdy}{\gM_{\text {dy}}}

\newcommand{\gMsp}{\gM_{\text {sp}}}

\newcommand{\akiadd}
{Department of Mathematics\\ Oklahoma State University \\
Stillwater OK 74078}
\newcommand{\kableadd}%
{Department of Mathematics\\ Cornell University\\
Ithaca NY 14853}
\newcommand{\akiemail}{yukie{@}math.okstate.edu}

\newcommand{\kableemail}{kable{@}math.cornell.edu}

\newcommand{\sub}{\subseteq}
\newcommand{\bk}{\backslash}
\newcommand{\ti}{\widetilde}

\newcommand{\ccd}{,\cdots,}

\newcommand{\beq}{\begin{equation}
\begin{aligned}}
\newcommand{\eeq}{\end{aligned}
\end{equation}}

\newtheorem{thm}{Theorem}[section]
\newtheorem{lem}[thm]{Lemma}
\newtheorem{cor}[thm]{Corollary}
\newtheorem{prop}[thm]{Proposition}

\theoremstyle{definition}

\newtheorem{defn}[thm]{Definition}

\newtheorem{cond}[thm]{Condition}

\theoremstyle{remark}

\newcommand{\kt}{{\widetilde{k}}}
\newcommand{\GL}[1]{{\operatorname{GL}(#1)}}
\newcommand{\intl}{\int\limits}
\newcommand{\calo}{{\mathcal{O}}}

\renewcommand{\k}{\kappa}

\newcommand{\trv}[1]{{\operatorname{Tr}_{\kt_v/k_v}(#1)}}
\newcommand{\normv}[1]{{\operatorname{N}_{\kt_v/k_v}(#1)}}

\newcommand{\diag}{\operatorname{diag}}
\newcommand{\cleq}{\preccurlyeq}
\newcommand{\tl}[1]{{\tilde{#1}}}

\begin{document}

\title[Mean value theorem]
{The mean value of the product of class numbers of paired quadratic fields I}
\author{Anthony C. Kable}
\address{\kableadd}
\email{\kableemail}
\author{Akihiko Yukie}
\address{\akiadd}
\email{\akiemail}
\date{\today}
\begin{abstract} 
Let $k$ be a number field and $\ti k$ a fixed 
quadratic extension of $k$.  
In this paper and its companion, we find the mean value of the 
product of class numbers and regulators of two quadratic extensions 
$F,F^*\not=\ti k$ contained in the biquadratic extensions of $k$
containing $\ti k$.   
\end{abstract}
\maketitle

\section{Introduction}\label{intro}

If $k$ is a number field then let $\Del_k$, $h_k$ 
and $R_k$ be the absolute discriminant 
(which is an integer), the class number and the regulator,
respectively.  
We fix a number field $k$ and a 
quadratic extension $\ti k$ of $k$.
If $F\not=\ti k$ 
is another quadratic extension of $k$, let $\ti F$ be the
compositum of $F$ and $\ti k$.  Then $\ti F$ is a biquadratic 
extension of $k$ and so contains precisely three quadratic 
extensions, $\ti k$, $F$ and, say, $F^*$ of $k$.  We say that 
$F$ and $F^*$ are \emph{paired}.  In this paper and its companion 
\cite{kable-yukie-pbh-II}, we shall 
find  the mean value of $h_FR_Fh_{F^*}R_{F^*}$ or, equivalently, 
the mean value of $h_{\ti F}R_{\ti F}$ with respect to $|\Del_F|$.  

Our main results are Theorem \ref{theorem-main} 
and Corollaries \ref{cor-main1} and \ref{cor-main2}
in which $k$ is an arbitrary number field and $F$ runs 
through quadratic extensions with 
given local behaviors at a fixed finite number of places.
However, for the sake of simplicity, we state our results here
assuming that $k=\Q$ and that $F$ runs through either real or 
imaginary quadratic extensions of $\Q$ without any further 
local conditions.  

Let $\ti k = \Q(\sqrt{d_0})$ where $d_0\not=1$ is a square
free integer.  Suppose $|\Del_{\Q(\sqrt{d_0})}| 
= \prod_p p^{\ti\del_p(d_0)}$ is the prime decomposition.    
Note that $\ti\del_p(d_0)>0$ 
if and only if $p$ is ramified in $\Q(\sqrt{d_0})$.  
Moreover, if $p\not=2$ is ramified in $\Q(\sqrt{d_0})$ 
then $\ti\del_p(d_0)=1$, and if $p=2$ then 
$\ti\del_p(d_0)=2$ when $d_0\equiv 3 \;(4)$
and  $\ti\del_p(d_0)=3$ when 
$d_0$ is an even number.   Note that if  $d_0\equiv 1,5 \;(8)$ 
then the prime $2$ is split or inert in $\Q(\sqrt{d_0})$, respectively. 

For any prime number $p$, we put 
\begin{equation*} 
E'_p(d_0) = \begin{cases} 
1-3p^{-3}+2p^{-4}+p^{-5}-2p^{-6} 
& \text{if $p$ is split in}\; \ti k, \\
(1+p^{-2})(1-p^{-2}-p^{-3}+p^{-4}) 
& \text{if $p$ is inert in}\; \ti k, \\
(1-p^{-1})(1+p^{-2}-p^{-3}+p^{-2\ti\del_p(d_0)-2\lfloor\ti\del_p(d_0)
/2\rfloor-1})
& \text{if $p$ is ramified in }\; \ti k
,\end{cases}
\end{equation*} 
where $\lfloor\ti\del_p(d_0)/2\rfloor$ is the largest integer less
than or equal to $\ti\del_p(d_0)/2$.

We define 
\begin{equation*} 
\begin{aligned}
c_+ (d_0) & = \begin{cases} 16 & d_0>0, \\ 
               8\pi   & d_0<0, \end{cases} \quad 
c_- (d_0) = \begin{cases} 4\pi^2 & d_0>0, \\ 
               8\pi   & d_0<0, \end{cases} \\
M(d_0) & = |\Del_{\Q(\sqrt{d_0})}|^{\frac 12} 
\zeta_{\Q(\sqrt{d_0})}(2)\prod_p E'_p(d_0)
.\end{aligned}
\end{equation*} 

The following theorems are special cases of 
Corollaries \ref{cor-main1}, \ref{cor-main2}.  
\begin{thm}\label{simple-mainthm1}
With either choice of sign we have
\begin{equation*} 
\lim_{X\to\infty} X^{-2} 
\sum_{\substack{[F:\Q]=2, \\ 0< \pm \Del_F < X}}
h_FR_Fh_{F^*}R_{F^*} 
= c_{\pm}(d_0)^{-1} M(d_0)
.\end{equation*} 
\end{thm} 
\begin{thm}\label{simple-mainthm2}
With either choice of sign we have
\begin{equation*} 
\lim_{X\to\infty} X^{-2} 
\sum_{\substack{[F:\Q]=2, \\ 0< \pm \Del_F < X}} 
h_{F(\sqrt{d_0})}R_{F(\sqrt{d_0})}
=c_{\pm}(d_0)^{-1} h_{\Q(\sqrt{d_0})} R_{\Q(\sqrt{d_0})} M(d_0)
.\end{equation*} 
\end{thm} 

Note that in Theorem \ref{simple-mainthm1} if $d_0>0$ and 
$\Del_F<0$ then both $F$ and $F^*$ are imaginary quadratic 
fields and so Theorem \ref{simple-mainthm1} states that 
\begin{equation*}
\lim_{X\to\infty}X^{-2}
\sum_{\substack{[F:\Q]=2, \\ 0 < -\Del_F < X}} h_Fh_{F^*}
=\tfrac{1}{4\pi^2}M(d_0)
,\end{equation*}
which reflects the titles of this paper and \cite{kable-yukie-pbh-II}.  

Theorems of this kind are called \emph{density theorems}.  
Many density theorems are known in number theory including, for
example, the prime number theorem, the theorem of Davenport-Heilbronn
\cite{dahea}, \cite{daheb} on the density of 
the number of cubic fields and the 
theorem of Goldfeld-Hoffstein \cite{goho}
on the density of class number
times regulator of quadratic fields.  

Among the three density theorems we quoted above, 
the prime number theorem, which is probably the 
best known density theorem, is of a more multiplicative 
nature than the other two theorems and our result
has more similarities to these. We would like to point out that the
Euler factor $1-p^{-2}-p^{-3}+p^{-4}$, which appears in $E_p'(d_0)$ in
our result when $p$ is inert, also occurred in the Goldfeld-Hoffstein
theorem at every odd prime. We do not as yet understand the
significance of this coincidence.

The original proof of the Davenport--Heilbronn theorem  used  
the  ``fundamental domain method" and the original 
proof of the Goldfeld--Hoffstein 
theorem used  Eisenstein series 
of half-integral weight.  However, these two theorems
can also be proved using the zeta function theory of \pv s.
The Davenport--Heilbronn theorem corresponds to the 
space of binary cubic forms and the Goldfeld--Hoffstein theorem
corresponds to the space of binary quadratic forms.  
The global theory of these two cases was investigated 
extensively by Shintani in \cite{shintania}, \cite{shintanib}.
The local theory and the proof of the density theorem, 
which use the global theory carried out by Shintani,  
were done by Datskovsky and Wright
\cite{dawra}, \cite{dawrb} in the first case 
and by Datskovsky \cite{dats} in the 
second (also correcting a minor
error in the constant appearing in the Goldfeld-Hoffstein theorem).    
This zeta function theory of \pv s is the approach we 
take to prove Theorems \ref{simple-mainthm1} and \ref{simple-mainthm2}.

We now recall the definition of \pv s.  
Let $G$ be a reductive group and 
$V$ a \rep{} of $G$ both of which are defined over 
an arbitrary field $k$ of characteristic zero.   
For simplicity, we assume that $V$ is an irreducible
\rep{} of $G$.  
\begin{defn} \label{pvdefn}
The pair $(G,V)$ is called a \emph{\pv}{} if 
\begin{itemize}
\item[(1)] there exists a Zariski open $G$-orbit in $V$ and 
\item[(2)] there exists a non-constant polynomial $P(x)\in k[V]$
and a rational character $\chi(g)$ of $G$ such that 
$P(gx) = \chi(g)P(x)$ for all $g\in G$ and $x\in V$.  
\end{itemize}
\end{defn} 
Any polynomial $P(x)$ in the above definition is called a relative 
invariant polynomial.  It is known that if $P(x)$ is the 
relative invariant polynomial of the lowest degree then
any other relative invariant polynomial is a constant multiple of 
a power of $P(x)$.  So if we put $V^{\sst}=\{x\in V\mid P(x)\not=0\}$
then this definition does not depend on the choice of $P(x)$.

The notion of \pv s was introduced by M. Sato in the early 1960's.  
The principal parts of global zeta functions 
for some \pv s have been determined 
by Shintani \cite{shintania}, \cite{shintanib}, 
and the second author \cite{yukiec}, \cite{yukiei}.  
Roughly speaking, the global zeta function is a 
counting function for the unnormalized Tamagawa numbers of the
stabilizers of points in $\rv^{\sst}$.  This interpretation of  
expected density theorems for \pv s  
is discussed in the introduction to
\cite{wryu} and in section 5 of \cite{kayu}, p. 342,
in some cases including those we will consider in this paper.   
Unfortunately, the global zeta function is not exactly 
this counting function and  Datskovsky and Wright formulated 
in \cite{dawrb} what we call the filtering process 
to deal with this difficulty.    

To explain the need for the filtering process we consider the  space of 
binary quadratic forms.  Gauss made a conjecture in \cite{gauss}
on the density of class number times regulator of orders 
in quadratic fields. This conjecture was 
proved by Lipschutz \cite{lipschutz}
in the case of imaginary quadratic fields and by Siegel \cite{siegele}
in the case of real quadratic fields,
and much work has been done on the error term estimate also 
(see Shintani \cite{shintanib}, pp. 44, 45 and  
Chamizo--Iwaniec \cite{chamizo-iwaniec-errorterm}, for example).  
However, each quadratic field has infinitely many orders and so 
we must filter out this repetition in order to obtain the density of 
class number times regulator for quadratic fields. 

In order to apply the filtering process it is necessary to 
carry out at least the following steps:
\begin{itemize} 
\item[(1)] Find the principal part  of the global zeta function 
at its rightmost pole.  
\item[(2)] Find a uniform estimate of the standard local zeta functions.
\item[(3)] Find the local densities.  
\end{itemize}
Note that, despite Tauberian theory, (1) is necessary even to show
the existence of the density.  The standard local zeta functions
will be defined in section \ref{plan}.  
If we apply the filtering process
the constant in the density theorem will have an 
Euler product and we call
the Euler factor the \emph{local density}.  Also, we must point out 
that 
the present formulation of the filtering process does not allow
us to use the poles of the global zeta function
other than the rightmost pole, as can be done in the case of integral
equivalence classes. It is an important problem 
in the future to improve the filtering process so that we 
can get error term estimates. However, although it does not, in its
current form, yield an error term, our approach does appear to be
the only one presently available which allows the field $k$ to be a
general number field rather than just $\Q$.

Let $\aff^n$ be $n$-dimensional 
affine space regarded as a variety 
over  the ground field $k$.   
Let $\ti k$ be a fixed quadratic extension of 
$k$, $W$  the space of binary $\ti k$-Hermitian forms
and $\m(2,2)$ the 
space of $2\times 2$ matrices.  We regard $\gl(2)_{\ti k}$ 
as a group over $k$.  In this paper and \cite{kable-yukie-pbh-II}, 
we consider the following two \pv s:
\begin{itemize} 
\item[(1)] $G = \gl(2)\times \gl(2)\times \gl(2),
\quad V = \m(2,2)\otimes \aff^2$, 
\item[(2)] $G = \gl(2)_{\ti k}\times \gl(2),
\quad V = W\otimes \aff^2$.
\end{itemize}

Case (2) is a $k$-form of case (1).  
We gave an interpretation of the 
expected density theorem for case (2) in 
section 5 of \cite{kayu}.  Let $k$ be a number field
and $\ti G$ be the image of $G$ in $\gl(V)$.    
For $x\in  \rv^{\sst}$,
let $\ti G^{\circ}_x$ be the identity component of the 
stabilizer.  In case (2), the orbit space $\rg\bk \rv^{\sst}$
corresponds bijectively with quadratic extensions of $k$ and, 
if $x$ corresponds to fields other than $k$ and $\ti k$, 
the weighting factor in the density theorem is 
the unnormalized Tamagawa number of $\ti G^{\circ}_x$, which is more or 
less $h_FR_Fh_{F^*}R_{F^*}$ or $h_{\ti F}R_{\ti F}$.  
The  principal part at the rightmost pole
of the global 
zeta function for this case was  obtained  in \cite{yukieh}, 
Corollary 8.16. Therefore it remains to carry out steps (2) and (3) 
of the filtering process.    

Let $v$ be a finite place of a number field $k$ and $k_v$
its completion at $v$.  The local zeta functions we 
consider are certain integrals over $G_{k_v}$-orbits in  
$V^{\sst}_{k_v}$.  The analogous 
integral over the set $V^{\sst}_{k_v}$ is called the 
\emph{Igusa zeta function}.  Igusa has made significant contributions
to the computation of this type of integral
(see \cite{igusaa}, \cite{igusai}, \cite{igusaf}
\cite{igusae},  \cite{igusaj}, \cite{igusak}), 
and the explicit form of the Igusa zeta
function is known for many cases.  However, we need information
on integrals over orbits and we cannot deduce 
a uniform estimate from the present knowledge of 
Igusa zeta functions.  Datskovsky and Wright \cite{dawra} and 
Datskovsky \cite{dats} accomplished the uniform estimate 
of the standard local zeta functions by explicitly computing them
at all finite places.  However, as the rank of the 
group grows, it becomes increasingly difficult to compute the
explicit forms of the standard local zeta functions, especially 
at special places such as dyadic places, and we have to be 
abstemious with our labor. So we shall only 
prove a uniform estimate of the standard local zeta functions 
at all but finitely many places, without
finding their  explicit forms.    

We follow Datskovsky's approach in \cite{dats}
(which can also be seen implicitly in \cite{daheb})
to find the local densities. We must consider biquadratic extensions
and consequently the dyadic places of $k$ are difficult and technical
to handle, given the possible appearance of wild ramification.
For this reason, we have a separate paper
\cite{kable-yukie-pbh-II} to deal with 
biquadratic extensions generated by two ramified quadratic
extensions over a dyadic field.  However, the reader 
should be able to find all the main ingredients for proving 
Theorems \ref{simple-mainthm1} and \ref{simple-mainthm2}
in this part.

For the rest of the introduction
we discuss the organization of this paper.
Throughout, except for section \ref{space},
$k$ is a fixed number field and $\ti k$ is a fixed 
quadratic extension of $k$.  In section \ref{space},
$k$ is an arbitrary field of characteristic zero and 
$\ti k$ is a quadratic extension of it.
In section \ref{notation}
we describe notation we use throughout the paper.
In section \ref{space} we review  from \cite{kayu} 
the interpretation of the 
orbit space $\rg\bk\rv^{\sst}$ for the \pv s (1) and (2) above
and fix 
parametrizations of the stabilizers of certain points in  
$\rv^{\sst}$.     
In section \ref{measure} we fix various normalizations regarding 
the invariant measure on $\gl(2)$ both locally and globally.  
In section \ref{formulation} we define a measure on the stabilizer 
of each point in $V^{\sst}$, both locally and globally, which is 
in some sense canonical and prove that the volume of 
$\ti G^{\circ}_{x\,\A}/\ti G^{\circ}_{x\, k}$ 
is the unnormalized Tamagawa number of 
$\ti G^{\circ}_x$.   As we mentioned above, this volume is the weighting  
factor in the density theorem.  
We also introduce the local zeta functions.  
In section \ref{plan} we first define and review the
analytic properties of the global zeta function.  Then we define 
the standard local zeta functions 
and express the global zeta function in terms of them, 
thus making it more or less a counting function 
for $h_FR_Fh_{F^*}R_{F^*}$.  The final and 
most important purpose of this section is to review
the filtering process and to identify the conditions under which 
it works. Assuming these conditions, we then deduce 
a preliminary density theorem involving certain as yet unevaluated
constants.
In section \ref{final} we list the values of those constants 
from later sections
and state the final form of the density theorem.  Therefore, sections
\ref{plan} and \ref{final} are the heart of this paper and 
\cite{kable-yukie-pbh-II}.  After finishing these sections, the reader 
should understand the outline of the proof of  our result.  
Later sections are devoted to verifying the conditions 
mentioned above and to evaluating the constants involved.  
In section \ref{sec-omega} we introduce the notion of an omega set
following Datskovsky \cite{dats}.  This set is used to describe 
the local zeta function, which is an 
integral over a subset of the vector 
space, as an integral over the group.   
In section \ref{sec-estimate} we prove a uniform estimate for the
standard local zeta functions.  
In sections \ref{sec-iv}--\ref{sec-infinite}
we find the constants needed to describe
the final form of the density theorem, except for a few
cases at dyadic places.  Those dyadic cases are handled in part II 
and a separate introduction describing the particular technical
difficulties associated with these cases will be given there.

\section {Notation} \label{notation} 

This section is confined to establishing our basic notational
conventions. Additional notation required throughout the paper will
be introduced and explained in the next three sections. 
More specialized notation will be introduced in the
section where it is required.

If $X$ is a finite set then $\# X$ will denote its cardinality. The
standard symbols $\Q$, $\R$, $\C$ and $\Z$ will denote respectively
the rational, real and complex numbers and the rational integers.
If $a\in\R$ then the largest integer $z$ such that $z\leq a$ is
denoted $\lfloor a\rfloor$ and the smallest integer $z$ such that
$z\geq a$ by $\lceil a\rceil$. The set of positive real numbers is
denoted $\R_{+}$. If $R$ is any ring then $R^{\times}$ is the set
of invertible elements of $R$ and if $V$ is a variety defined over
$R$ then $V_R$ denotes its $R$-points. If $G$ is an algebraic group
then $G^{\circ}$ denotes its identity component.

Both $k$ and $\kt$ are number fields and so each number theoretic
object we introduce for $k$ has its counterpart for $\kt$.
Generally the notation for the $\kt$ object will be derived from
that of the $k$ object by adding a tilde. Let $\gM$,
$\gM_{\infty}$, $\gM_{\text{f}}$, $\gM_{\text{dy}}$, 
$\gM_{\R}$ and $\gM_{\C}$
denote respectively the set of all places of $k$, all infinite
places, all finite places, all dyadic places (those dividing the
place of $\Q$ at 2), all real places and all complex places.
(Correspondingly we have $\ti \gM$ and so on.) Let $\gM_{\text{rm}}$,
$\gM_{\text{in}}$ and $\gM_{\text{sp}}$ be the sets of places of $k$ which are
respectively ramified, inert and split on extension to $\kt$.
Recall that a real place of $k$ which lies under a complex place of
$\kt$ is regarded as ramified.

Let $\calo$ be the ring of integers of $k$. If $v\in\gM$ then $k_v$
denotes the completion of $k$ at $v$ and $|\;|_v$ denotes the normalized
absolute value on $k_v$. If $v\in\gMf$ then $\calo_v$ denotes
the ring of integers of $k_v$, $\pi_v$ a uniformizer in $\calo_v$,
$\gp_v$ the maximal ideal of $\calo_v$ and $q_v$ the cardinality of
$\calo_v/\gp_v$. If $a\in k_v$ and $(a)=\gp_v^i$ then we write 
$\ord_{k_v}(a)=i$. If $\gi$ is a fractional ideal in $k_v$ and
$a-b\in\gi$ then we write $a\equiv b\;(\gi)$ or $a\equiv b\;(c)$ if $c$
generates $\gi$.

If $k_1/k_2$ is a finite extension either of local fields or of
number fields then we shall write $\Del_{k_1/k_2}$ for the relative
discriminant of the extension; it is an ideal in the ring of
integers of $k_2$. The symbol $\Del_{k_1}$ will stand for
$\Del_{k_1/\Q_p}$ or $\Del_{k_1/\Q}$ according as the situation is
local or global. To ease the notational burden we shall use the
same symbol, $\Del_{k_1}$, for the classical absolute discriminant
of $k_1$ over $\Q$. Since this number generates the ideal
$\Del_{k_1}$, the resulting notational identification is harmless.
If $\gi$ is a fractional ideal in the number field $k_1$ and $v$ is
a finite place of $k_1$ then we write $\gi_v$ for the closure of
$\gi$ in $k_{1,v}$. It is a fractional ideal in $k_{1,v}$. If $\gi$
is integral then we put $\esN(\gi)=\#(\calo_{k_1}/\gi)$. Note that
$\esN(\gi)=\prod_v\esN_v(\gi_v)$, where the product is over all
finite places of $k_1$ and $\esN_v(\gp_v^a)=q_v^a$ for $a\in\Z$. This
formula serves to extend the domain of $\esN$ to all fractional 
ideals in $k_1$.
We shall use the notation $\tr_{k_1/k_2}$ and $\n_{k_1/k_2}$ for
the trace and the norm in the extension $k_1/k_2$.

Returning to $k$, we let $r_1$, $r_2$, $h_k$, $R_k$ and $e_k$ be
respectively the number of real places, the number of complex
places, the class number, the regulator and the number of roots of
unity contained in $k$. It will be convenient to set
\begin{equation}
\gC_k=2^{r_1}(2\pi)^{r_2}h_kR_ke_k^{-1}
.\end{equation}

We assume that the reader is familiar with the basic definitions
and facts concerning adeles. These may be found in \cite{weilc}.
The ring of adeles, the group of ideles 
and the adelic absolute value of $k$ 
are denoted by  $\A$, $\ma$ and $|\;|$ respectively.  
When we have to show the number field or the local field 
on which we consider the absolute value, we may use notation 
such as $|\;|_F$.   There is a natural
inclusion $\A\to \ti \A$, 
under which an adele $(a_v)_v$ corresponds to the adele
$(b_w)_w$ such that $b_w=a_v$ if $w|v$. 
Let $\A^1=\{t\in \ma\mid |t|=1\}.$
Using the identification $\ti k\otimes_{k} \A \cong \ti \A$ ,
the norm map 
$\n_{\ti k/k}$ can be extended to a map from $\ti \A$ to $\A$.  
It is known (see p. 139 in \cite{weilc}) that 
$|\n_{\ti k/k}(t) | = |t|_{\ti \A}$ for $t\in \ti \A$. 
Suppose $[k:\Q]=n$.  Then $[\ti k:\Q]=2n$.  
For $\lam\in\mr$, $\lamb\in \ma$ 
is the idele whose component at any infinite place is $\lam^{\frac 1n}$
and whose component at any finite place is $1$.  
Also $\ti \lamb\in \ti\A^{\times}$ 
is the idele whose component at 
any infinite place is $\lam^{\frac 1{2n}}$
and whose component at any finite place is $1$.  
Clearly $\lamb = \ti \lamb{}^2$.  Since 
$|\lamb|=\lam$ and $|\ti \lamb|_{\ti \A}=\lam$ we conclude
that $|\lamb|_{\ti \A}= \lam^2$.  When we have to show the 
number field on which we consider $\lamb$, we use the notation 
such as $\lamb_F$.

If $V$ is a vector space over $k$ we let $V_{\A}$ be its adelization
and $V_{\infty}$ and $V_{\text{f}}$ its infinite and finite parts.  
Let $\cS(V_{\A})$, $\cS(V_{\infty})$, $\cS(V_{\text{f}})$ and $\cS(V_{k_v})$
be the spaces of 
\bs s on each of the indicated domains.

We choose a Haar measure  $dx$ on $\A$ so 
that $\int_{\A/k} dx = 1$.
For any  $v\in \gM_{\text{f}}$, 
we choose a Haar measure $dx_v$ on $k_v$ 
so that $\int_{\co_v}dx_v=1$.  
We use the ordinary Lebesgue measure $dx_v$ for $v$ real, 
and  $dx_v\wedge d\bar{x}_v$ for $v$ imaginary.  
Then $dx=|\Delta_k|^{-\frac 1 2}\prod_v
dx_v$  (see \cite{weilc}, p. 91).  
\par

We define a Haar measure 
$\md t^1$ on $\A^1$ so that 
$\int_{\A^1/\mk} \md t^1 =  1$.  Using 
this measure, we choose a Haar
measure $\md t$ on $\ma$  so that 
$$
\int_{\ma} f(t) \md t  = \int_0^{\infty}   
\int_{\bba^1}  f(\lamb t^1) \md \lam \md t^1
,$$ 
where $\md \lam = \lam^{-1}d\lam$.
For any $v\in \gM_{\text{f}}$, we choose a Haar measure $\md t_v$
on $\mk_v$ so that $\int_{\co_v^{\times}} \md t_v = 1$.  
Let $\md t_v(x)=|x|_v^{-1}d x_v$ if $v$ is real, and
$\md t_v(x)=|x|_v^{-1}d x_v\wedge d\bar{x}_v$ 
if $v$ is imaginary.  
Then $\md t=\gC_k^{-1}\prod_v \md t_v$  
(see \cite{weilc}, p. 95).  We later have to 
compare the global measure and the product of 
local measures and for that purpose it is 
convenient to notate the product of local measures
on $\A,\A^{\times}$ as follows 
\begin{equation} \label{productmeasures} 
d_{\text{pr}} x = \prod_v dx_v\,,\; 
d_{\text{pr}}^{\times} t = \prod_v d^{\times} t_v\,
.\end{equation}

\par
Let $\zeta_k(s)$ be the Dedekind zeta function of $k$.
We define
\begin{equation}\label{Zkdefn}
Z_k(s)=|\Delta_k|^{\frac s 2}
\left(\pi^{-\frac s 2}\Gam\left(\frac s 2\right)\right)^{r_1}
\left((2\pi)^{1-s}\Gam(s)\right)^{r_2}\zeta_k(s)\,
.\end{equation}
This definition differs from that in \cite{weilc}, p.129 by the
inclusion of the $|\Del_{k}|^{s/2}$ factor and from that in 
\cite{yukiec} by a factor of $(2\pi)^{r_2}$. It is adopted here as
the most convenient for our purposes. We note that it was the
quotient $Z_{k}(s)/Z_{k}(s+1)$ rather than $Z_{k}(s)$ itself which
played a significant role in \cite{yukiec} and this quotient is
unchanged here.
It is known 
(\cite{weilc}, p.129) that 
\begin{equation} \label{zetaresidues} 
\operatorname{Res}_{s=1} \zeta_k(s) = 
|\Delta_k|^{-{\frac 12}}\gC_k\text{\quad and so\quad}
\operatorname{Res}_{s=1} Z_k(s) = \gC_k\,
.\end{equation}

Finally, we introduce the following notation
\begin{equation} 
a(t_1,t_2)= \pmatrix t_1 & 0\\ 0 & t_2\endpmatrix,\;
n(u) = \pmatrix 1 & 0\\ u & 1\endpmatrix\,.
\end{equation}

\section{A review of the orbit space} \label{space}

This section is devoted to defining the prehomogeneous vector
spaces which are at the heart of this work and reviewing their
fundamental properties. Arithmetic plays no role here, so in this
section $k$ may be any field of characteristic zero and $\kt$ any
quadratic extension of $k$. We denote the non-identity element of
$\gal(\kt/k)$ by $\sigma$.

A matrix $x\in\m(2,2)_{\kt}$ is called \emph{Hermitian} if
${}^{t}x=x^{\sigma}$. The set of all Hermitian matrices in
$\m(2,2)_{\kt}$ forms a $k$-vector space which we shall denote by
$W$. The elements of $W$ are also referred to as binary Hermitian
forms.

We define and discuss the two spaces we require in parallel as far
as possible; they will be distinguished as cases $(1)$ and $(2)$.
Let
\begin{equation}
V=\begin{cases} \m(2,2)\otimes\aff^2 & \mbox{\quad in case (1),}
\\
W\otimes\aff^2 & \mbox{\quad in case (2),}
\end{cases}
\end{equation}
where $\aff^n$ is $n$-dimensional affine space regarded as a
variety over $k$. Let 
\begin{equation}
G=\begin{cases} \GL{2}\times\GL{2}\times\GL{2} &
\mbox{\quad in case (1),} \\
\GL{2}_{\kt}\times\GL{2} &
\mbox{\quad in case (2),}
\end{cases}
\end{equation}
where $\GL{2}_{\kt}$ is regarded as an algebraic group over $k$ by
restriction of scalars. If $g\in G$ then we shall write
$g=(g_1,g_2,g_3)$ in case (1) and $g=(g_1,g_2)$ in case (2).
It will be convenient to identify
$x=(x_1,x_2)\in V$ with the $2\times 2$-matrix $M_x(v)=v_1x_1+v_2x_2$
of linear forms in the variables $v_1$ and $v_2$, which we collect
into the row vector $v=(v_1,v_2)$. With this identification, we
define a rational action of $G$ on $V$ via
\begin{equation}
M_{gx}(v)=\begin{cases} g_1M_x(vg_3)\,{}^{t}g_2 &
\mbox{\quad in case (1),} \\
g_1M_x(vg_2)\,{}^{t}g_1^{\sigma} &
\mbox{\quad in case (2).}
\end{cases}
\end{equation}
In both cases we define $F_x(v) = -\det M_x(v)$.  
Then
\begin{equation}\label{emap}
F_{gx}(v) = \begin{cases}
\det g_1\det g_2 F_x(v g_3) & \mbox{\quad in case (1),} \\
\n_{\ti k/k}(\det g_1) F_x(v g_2) & \mbox{\quad in case (2).}
\end{cases}
\end{equation} 
We let $P(x)$ be the discriminant of the binary quadratic form
$F_x(v)$. 
Then $P(x)\in k[V]$ and $P(gx)=\chi(g)P(x)$ where 
\begin{equation} 
\chi(g) = \begin{cases} 
(\det g_1\det g_2\det g_3)^2 & \mbox{\quad in case (1),} \\
(\n_{\ti k/k}(\det g_1)\det g_3)^2 & \mbox{\quad in case (2).}
\end{cases}
\end{equation}
A calculation shows that 
$P(x)$ is not identically zero and so it is a  relative invariant
polynomial for $(G,V)$ in each case.  
We let $V^{\sst}$ denote the complement of the hypersurface defined
by $P(x)=0$ in $V$.

We define $\ti T=\ker(G\to \gl(V))$; in case (1) 
\begin{equation} 
\ti T =
\{ (t_1I_2,t_2I_2,t_3I_2)\mid 
t_1,t_2,t_3\in \gl(1),\; t_1t_2t_3=1\}
\end{equation}
and in case (2)
\begin{equation}
\ti T=\{ (t_1I_2,t_2I_2)\mid 
t_1\in \gl(1)_{\ti k},t_2\in \gl(1),\; 
\n_{\ti k/k}(t_1)t_2=1\} \,.
\end{equation}
It will be convenient to introduce standard coordinates on $G$ and
$V$. Elements of $G$ have the form $g=(g_1,g_2,g_3)$ or
$g=(g_1,g_2)$. In either case we shall write
\begin{equation}\label{groupcoords}
g_i=\begin{pmatrix}g_{i11}&g_{i12}\\g_{i21}&g_{i22}\end{pmatrix}
\end{equation}
for each $i$. Elements of $V$ are vectors $x=(x_1,x_2)$. We shall
put
\begin{equation}\label{case(1)coords}
x_i=\begin{pmatrix}x_{i11}&x_{i12}\\x_{i21}&x_{i22}\end{pmatrix}
\end{equation}
in case (1) and
\begin{equation}\label{case(2)coords}
x_i=\begin{pmatrix}x_{i0}&x_{i1}\\x_{i1}^{\sigma}&x_{i2}
\end{pmatrix}
\end{equation}
in case (2).

In the language of Galois descent, case (2) is a $k$-form of case
(1); they become isomorphic on extension of scalars from $k$ to
$\kt$. Indeed it is well known that, as $\kt$-varieties,
\begin{equation}
G\times\kt \cong \GL{2}\times\GL{2}\times\GL{2}
\end{equation}
and
\begin{equation}
W\times\kt \cong \m(2,2)
\end{equation}
so that
\begin{equation}
V\times\kt\cong\m(2,2)\otimes\aff^2,
\end{equation}
and a calculation shows that the induced action of $G\times\kt$ on
$V\times\kt$ is that of case (1). The Galois automorphism $\sigma$
induces a $k$-automorphism of the $k$-varieties $G$ and $V$ which we
denote by $i(\sigma)$. If $(g_1,g_2,g_3)\in G_{\kt}$ then
$i(\sigma)(g_1,g_2,g_3)=(g_2^{\sigma},g_1^{\sigma},g_3^{\sigma})$
and if $x\in W_{\kt}$ then $i(\sigma)x={}^{t}\!x^{\sigma}$, where
$\sigma$ as a superscript denotes the entry-by-entry action of
$\sigma$. In particular, $G_k$ is embedded in $G_{\kt}\cong(G\times
\kt)_{\kt}$ via the map $(g_1,g_2)\mapsto(g_1,g_1^{\sigma},g_2)$.

We are now ready to recall the description of the space of
non-singular orbits in $V_k$.

\begin{defn} Let $\gEx_2$ be the set of isomorphism classes of 
extensions of $k$ of degree at most two.  
\end{defn}
It is proved in \cite{wryu}, pp. 305--310 and \cite{kayu}, p. 324 that 
$\rg\bk \rv^{\sst}$ corresponds bijectively 
with $\gEx_2$.  Moreover if $x\in V$ then the corresponding field 
is generated by the roots of $F_x(v)=0$.  We denote this field by 
$k(x)$.  

Suppose that $p(z)=z^2+a_1z+a_2\in k[z]$ has distinct roots $\al_1$
and $\al_2$. We collect these into a set $\al=\{\al_1,\al_2\}$
since the numbering is arbitrary. Define $w_p\in V_k$ by
\begin{equation} 
w_p = \left(\pmatrix 0 & 1 \\ 1 & a_1\endpmatrix,
\pmatrix 1 & a_1\\ a_1 & a_1^2-a_2\endpmatrix\right)\,;
\end{equation} 
a computation shows that $F_{w_p}(z,1)=p(z)$ and so $w_p\in
V_k^{\sst}$ and $k(w_p)=k(\al)$ is the splitting field of $p$.
Let
\begin{equation}\label{wdefn}
w = \left(\pmatrix 1 & 0\\ 0 & 0\endpmatrix,
\pmatrix 0 & 0\\ 0 & 1\endpmatrix\right)\,,
\end{equation}
\begin{equation}
h_{\al}=\begin{pmatrix}1&-1\\-\al_1&\al_2\end{pmatrix}
\end{equation}
and then define $g_p\in G_{k(w_p)}$ by
\begin{equation}\label{g_pdefn}
g_p=\begin{cases}
(h_{\al},h_{\al},(\al_2-\al_1)^{-1}h_{\al}) &
\mbox{\quad in case (1) or when $k(w_p)=\kt$,} \\
(h_{\al},(\al_2-\al_1)^{-1}h_{\al}) &
\mbox{\quad otherwise.}
\end{cases}
\end{equation}
With these definitions it is easy to check that $w_p=g_pw$.

We close this section with a detailed description of the
$k$-rational points of the stabilizer $G_{w_p}$. Similar
descriptions were derived in \cite{kayu} and \cite{wryu} and,
although we are using different orbital representatives here, the
arguments are so similar that they will only be sketched. The
method is as follows: We begin with a description of $G_w$ as a
$k$-variety; this is given in section 3 of \cite{wryu} for case
(1) and in section 2 of \cite{kayu} for case (2). Then we find,
by direct calculation, the $k$-rational points in $g_p
G_{w\,k(w_p)} g_p^{-1}$ and this gives us $G_{w_p\,k}$.

If we let
\begin{equation}\label{standardtoruselement}
t=\begin{cases}
(a(t_{11},t_{12}),a(t_{21},t_{22}),a(t_{31},t_{32})) &
\mbox{\quad in case (1),} \\
(a(t_{11},t_{12}),a(t_{21},t_{22})) &
\mbox{\quad in case (2),}
\end{cases}
\end{equation}
then
\begin{equation}\label{stab(0)}
G_{w\,k}^{\circ}=\begin{cases}
\{t\mid t_{ij}\in k^{\times}, t_{1j}t_{2j}t_{3j}=1\ \forall\ 
i,j\} &
\mbox{\quad in case (1),} \\
\{t\mid t_{1j}\in\kt^{\times}, t_{2j}\in k^{\times},
\n_{\kt/k}(t_{1j})t_{2j}=1\ \forall\ j\} &
\mbox{\quad in case (2),}
\end{cases}
\end{equation}
and so $G_{w\,k}^{\circ}\cong\GL{1}_k^4$ in case (1) and
$G_{w\,k}^{\circ}\cong\GL{1}_{\kt}^2$ in case (2). Putting
\begin{equation}
\tau=\begin{pmatrix}0&1\\1&0\end{pmatrix}
\end{equation}
the class of $(\tau,\tau,\tau)$ in case (1) or of
$(\tau,\tau)$ in case (2) generates $G_{w\,k}/G_{w\,k}^{\circ}$.

Now let
\begin{equation}\label{telem}
t=\begin{cases}
(a(t_{11},t_{12}),a(t_{21},t_{22})) &
\mbox{\quad in case (2) when $k(w_p)\neq\kt$,} \\
(a(t_{11},t_{12}),a(t_{21},t_{22}),a(t_{31},t_{32})) &
\mbox{\quad otherwise.}
\end{cases}
\end{equation}
We assume that $k(w_p)/k$ is quadratic, since if $k(w_p)=k$ then
$G_{w_p\,k}$ is conjugate to $G_{w\,k}$ over $k$. Let $\nu$ be the
non-trivial element of $\gal(k(w_p)/k)$, which may also be thought
of as an element of $\gal(\kt(w_p)/\kt)$ when $k(w_p)\neq\kt$. Here
$\kt(w_p)$ denotes the compositum of $\kt$ and $k(w_p)$.

In case (1), $G_{w_p\,k}^{\circ}$ is
\begin{equation}\label{stab(1)}
\{g_ptg_p^{-1}\mid t_{ij}\in k(w_p)^{\times},
t_{i1}=t_{i2}^{\nu}, t_{1j}t_{2j}t_{3j}=1\ \forall\  i,j\}
\end{equation}
and so $G_{w_p\,k}^{\circ}\cong\GL{1}_{k(w_p)}\times\GL{1}_{k(w_p)}$.
In case (2) when $k(w_p)=\kt$, $G_{w_p\,k}^{\circ}$ is
\begin{equation}\label{stab(2)*}
\{g_ptg_p^{-1}\mid
t_{ij}\in\kt^{\times},
t_{12}^{\sigma}=t_{21},
t_{11}^{\sigma}=t_{22},
t_{1j}t_{2j}t_{3j}=1\ \forall\ i,j\}
\end{equation}
and so $G_{w_p\,k}^{\circ}\cong\GL{1}_{\kt}\times\GL{1}_{\kt}$. In
case (2) when $k(w_p)\neq\kt$, $G_{w_p\,k}^{\circ}$ is
\begin{equation}\label{stab(2)}
\{
g_ptg_p^{-1}\mid
t_{1j}\in\kt(w_p)^{\times},t_{2j}\in k(w_p)^{\times},
t_{11}^{\nu}=t_{12},
\n_{\kt(w_p)/k(w_p)}(t_{1j})t_{2j}=1
\ \forall\ j\}
\end{equation}
and so $G_{w_p\,k}^{\circ}\cong\GL{1}_{\kt(w_p)}$.
In every instance, $G_{w_p\,k}/G_{w_p\,k}^{\circ}$ is generated by
$g_p(\tau,\tau,\tau)g_p^{-1}$ or $g_p(\tau,\tau)g_p^{-1}$ as the
case may be.

It will be convenient to have an explicit description of how
$G_{w_p\,k}$ is embedded in $G_k$ in each case. To this end, define
\begin{equation}  \label{apdefn} 
A_p(c,d) = \pmatrix c & -d\\ a_2 d & c-a_1d\endpmatrix
\mbox{\quad and\quad}
\tau_p = \pmatrix -1 & 0\\ -a_1 & 1\endpmatrix\,.
\end{equation} 
It is easy to check that any matrix which has both
$\left(\begin{smallmatrix}1\\-\al_1\end{smallmatrix}\right)$
and
$\left(\begin{smallmatrix}1\\-\al_2\end{smallmatrix}\right)$
as eigenvectors must equal $A_p(c,d)$ for some $c$ and $d$. 
 Consequently, the set of all such matrices is closed
under multiplication, any two such matrices commute and if such a
matrix is invertible then its inverse lies in the same set.
\begin{lem}\label{stabexpl(1)}
In case (1), $G_{w_p\,k}^{\circ}$ consists of elements of $G_{k}$ of
the form 
\begin{equation}\label{stabform(1)}
(A_p(c_1,d_1),A_p(c_2,d_2),A_p(c_3,d_3))
\end{equation}
where
$c_i,d_i\in k$, $\det(A_p(c_i,d_i))\neq0$ for $i=1,2$ and $(c_3,d_3)$ is
related to $(c_1,d_1,c_2,d_2)$ by the equation
\begin{equation}\label{explicit(1)}
A_p(c_3,d_3)=A_p(c_1,d_1)^{-1}A_p(c_2,d_2)^{-1}\,.
\end{equation}
Moreover $[G_{w_p\,k}:G_{w_p\,k}^{\circ}]=2$ and
$G_{w_p\,k}/G_{w_p\,k}^{\circ}$ is generated by the class of
$(\tau_p,\tau_p,\tau_p)$.
\end{lem}
\begin{proof}
Suppose first that $k(w_p)=k$. Then
$G_{w_p\,k}^{\circ}=g_pG_{w\,k}^{\circ}g_p^{-1}$ and, by 
(\ref{stab(0)}), the elements of $G_{w\,k}^{\circ}$ may be
characterized as those $(g_1,g_2,g_3)\in G_{k}$ such that
$\left(\begin{smallmatrix}1\\0\end{smallmatrix}\right)$
and
$\left(\begin{smallmatrix}0\\-1\end{smallmatrix}\right)$
are both eigenvectors for each $g_i$ and $g_1g_2g_3=I_2$. Since
$h_{\al}
\left(\begin{smallmatrix}1\\0\end{smallmatrix}\right)=
\left(\begin{smallmatrix}1\\-\al_1\end{smallmatrix}\right)$ and
$h_{\al}
\left(\begin{smallmatrix}0\\-1\end{smallmatrix}\right)=
\left(\begin{smallmatrix}1\\-\al_2\end{smallmatrix}\right)$,
the first claim follows.
If $k(w_p)\neq k$ then calculation gives
$h_{\al}a(t,t^{\nu})h_{\al}^{-1}=A_p(c,d)$ where $t=c+d\al_1\in
k(w_p)$. With this observation, the first claim follows in this
case from (\ref{stab(1)}). Finally, $h_{\al}\tau
h_{\al}^{-1}=\tau_p$ and the second claim is established.
\end{proof}
\begin{lem}\label{stabexpl(2)}
In case (2), $G_{w_p\,k}^{\circ}$ consists of elements of $G_k$ of
the form
\begin{equation}\label{stabform(2)}
(A_p(c_1,d_1),A_p(c_2,d_2))
\end{equation}
where $c_1,d_1\in\kt$, $c_2,d_2\in k$, $\det(A_p(c_1,d_1))\neq 0$
and $(c_2,d_2)$ is related to $(c_1,d_1)$ by the equation
\begin{equation}\label{explicit(2)}
A_p(c_2,d_2)=A_p(c_1,d_1)^{-1}A_p(c_1^{\sigma},
d_1^{\sigma})^{-1}\,.
\end{equation}
Moreover, $[G_{w_p\,k}:G_{w_p\,k}^{\circ}]=2$ and
$G_{w_p\,k}/G_{w_p\,k}^{\circ}$ is generated by the class of
$(\tau_p,\tau_p)$.
\end{lem}
\begin{proof}
If $k(w_p)=k$ then, by (\ref{stab(0)}), $G_{w\,k}^{\circ}$ may be
characterized as the set of $(g_1,g_2)$ in $G_k$ such that
$\left(\begin{smallmatrix}1\\0\end{smallmatrix}\right)$
and
$\left(\begin{smallmatrix}0\\-1\end{smallmatrix}\right)$
are eigenvectors of $g_1$ and $g_1g_1^{\sigma}g_2=I_2$. Since
$G_{w_p\,k}^{\circ}=g_p G_{w\,k}^{\circ} g_p^{-1}$ and
$h_{\al}^{\sigma}=h_{\al}$, the claim follows. If $k(w_p)\neq
k,\kt$ then $h_{\al}^{\sigma}=h_{\al}$ and a similar argument works
on setting $A_p(c_1,d_1)=h_{\al}a(t_{11},t_{11}^{\nu})h_{\al}^{-1}$
in the notation of (\ref{stab(2)}).

This leaves the case where $k(w_p)=\kt$. We use the notation of 
(\ref{stab(2)*}). If we set
$g_1=h_{\al}a(t_{11},t_{12})h_{\al}^{-1}=A_p(c_1,d_1)$ for some
$c_1,d_1\in\kt$ then, using the equation
$h_{\al}^{\sigma}=-h_{\al}\tau$, we have
$g_1^{\sigma}=h_{\al}a(t_{12}^{\sigma},t_{11}^{\sigma})h_{\al}^{-1}$
and so $g_2=g_1^{-1}g_1^{-\sigma}$ is
$h_{\al}a((t_{11}t_{12}^{\sigma})^{-1},(
t_{11}^{\sigma}t_{12})^{-1})h_{\al}^{-1}$. Thus $(g_1,g_2)\in
G_{w_p\,k}^{\circ}$. Finally, we have $\tau_p=h_{\al}\tau
h_{\al}^{-1}$ and the last claim follows from this.
\end{proof}

\section{An invariant measure on $\gl(2)$} \label{measure}

Assume now that $k$ is a number field. In this section we choose an
invariant measure on $\GL{2}$ in both the local and adelic
situations.

Let $T\subseteq\GL{2}$ be the set of diagonal matrices and
$N\subseteq\GL{2}$ be the set of lower-triangular matrices whose
diagonal entries are $1$. Then $B=TN$ is a Borel subgroup of
$\GL{2}$. Let $T_+=\{\lambda=a(\underline{\lambda}_1,
\underline{\lambda}_2)\mid\lambda_1,\lambda_2\in\mathbb{R}_+\}$ and
$K=\prod_{v\in\mathfrak{M}}K_v$ where $K_v=\operatorname{O}(2)$ if
$v\in\mathfrak{M}_{\mathbb{R}}$, $K_v=\operatorname{U}(2)$ if
$v\in\mathfrak{M}_{\mathbb{C}}$ and $K_v=\GL{2}_{\calo_v}$ if
$v\in\mathfrak{M}_{\text{f}}$. The group $\GL{2}_{\mathbb{A}}$ has the
Iwasawa decomposition $\GL{2}_{\mathbb{A}}=K T_{\mathbb{A}}
N_{\mathbb{A}}$ and so any element $g\in\GL{2}_{\mathbb{A}}$ can be
expressed as $g=\k(g)t(g)n(u(g))$ where $\k(g)\in K$,
$t(g)=a(t_1(g),t_2(g))$ and $u(g)\in\mathbb{A}$.

The measure $du$ on $\mathbb{A}$ defined in section \ref{notation}
induces an invariant measure on $N_{\mathbb{A}}$. Since $K$ is
compact we can choose an invariant measure $d\k$ on it so that the
total volume of $K$ is $1$. On $T_{\mathbb{A}}$ we put $\md t=\md
t_1\,\md t_2$ for $t=a(t_1,t_2)$, where $\md t_j$ is the measure on
$\mathbb{A}^{\times}$ defined in section \ref{notation}. Then
$db=|t_1t_2^{-1}|^{-1}\,\md t\, du$ defines an invariant measure on
$B_{\mathbb{A}}$ and $dg=d\k \,db$ defines an invariant measure on
$\GL{2}_{\mathbb{A}}$.

We make parallel definitions of invariant measures on
$\GL{2}_{k_v}$, $K_v$, $B_{k_v}$, $N_{k_v}$ and $T_{k_v}$, which we
denote by $dg_v$, $d\k_v$, $db_v$, $du_v$ and $\md t_v$,
respectively.  As in section \ref{notation}, we denote
the product of local measures on $G_{\A}$ as follows
\begin{equation}
d_{\text{pr}} g = \prod_v dg_v\,
.\end{equation}
Then (see section \ref{notation}) we have
\begin{equation} \label{glrelation}
du=|\Delta_k|^{-1/2}\prod_{v}du_v,\quad
\md t=\mathfrak{C}_{k}^{-2}\prod_{v}\md t_v\text{\quad and
so\quad}
dg = |\Delta_k|^{-1/2}\gC_k^{-2} 
d_{\text{pr}} g
.\end{equation}

Let $\GL{2}_{\mathbb{A}}^{0}=\{g\in\GL{2}_{\mathbb{A}}\mid
|\det(g)|=1\}$. If, for $\lambda\in\mathbb{R}_+$, we define
$c(\lambda)=a(\underline{\lambda},\underline{\lambda})$ then any
element of $\GL{2}_{\mathbb{A}}$ may be written uniquely as
$g=c(\lambda)g^0$ with $g^0\in\GL{2}_{\mathbb{A}}^0$. We choose a
Haar measure on $\GL{2}_{\mathbb{A}}^0$ so that $dg=2\md \lambda\,
dg^0$. It is well-known that the volume of
$\GL{2}_{\A}^{0}/\GL{2}_{k}$ with respect to $dg^0$ is
\begin{equation}
\gV_{k}=1/\res_{s=1}(Z_{k}(s)/Z_{k}(s+1))=\gC_{k}^{-1}Z_{k}(2)
.\end{equation}

As in section \ref{notation}, we note that all these definitions
apply equally well to the number field $\kt$ and yield a measure
on $\GL{2}_{\widetilde{\mathbb{A}}}$ and so on.
Having chosen an invariant measure on $\GL{2}$ both locally and
adelically, we also get local and adelic invariant measures on $G$
by taking the relevant product measures in each case.

\section{The canonical measure on the stabilizer}\label{formulation}

In this section we shall define a measure on $G_{x\,\bba}^{\circ}$
for $x\in V_k^{\sst}$ which is canonical (in a sense made precise
by Proposition \ref{canonicity}) and compute the volume of
$G_{x\,\bba}^{\circ}/{\ti T}_{\bba} G_{x\,k}^{\circ}$ under this
measure. We also make a canonical choice of measure on the
stabilizer quotient $G_{\bba}/G_{x\,\bba}^{\circ}$ and define
constants $b_{x,v}$ which will play an essential role in what
follows.

Before we begin this task it will be convenient for bookkeeping
purposes to attach to each orbit in $V_{k_v}^{\sst}$ 
where $v\in\gM$, an index which
records the arithmetic properties of $v$ and of the extension of
$k_v$ corresponding to the orbit. The orbit corresponding to $k_v$
itself will have index (sp), (in) or (rm) according as $v$ is in
$\gM_{\text{sp}}$, $\gM_{\text{in}}$ or $\gM_{\text{rm}}$. The
orbit corresponding to the unique unramified quadratic extension of
$k_v$ will have index (sp ur), (in ur) and (rm ur) for
$v\in\gM_{\text{sp}}$, $v\in\gM_{\text{in}}$ and
$v\in\gM_{\text{rm}}$ respectively. An orbit corresponding to a
ramified quadratic extension of $k_v$ will have index (sp rm) if
$v\in\gM_{\text{sp}}$ and (in rm) if $v\in\gM_{\text{in}}$. If
$v\in\gM_{\text{rm}}$ then the orbits corresponding to ramified
quadratic extensions of $k_v$ are subdivided into three types; the
one corresponding to $\kt$ has index (rm rm)*, those corresponding
to quadratic extensions $k(x)/k$ such that $k(x)\neq\kt$ and
$\kt\cdot k(x)/\kt$ is unramified have index (rm rm ur) and those
corresponding to quadratic extensions $k(x)/k$ such that
$k(x)\neq\kt$ and $\kt\cdot k(x)/\kt$ is ramified have index (rm rm
rm). This last type can occur only if $v\in\gM_{\text{dy}}$. 

Let $v\in\gM$ and $x\in V_{k_v}^{\sst}$. If
$v\notin\gM_{\text{sp}}$ then $v$ extends uniquely to a place of
$\kt$ which we also denote by $v$. In this case $\kt_v\cong
k_v\otimes_k\kt$. We denote by $\kt_v(x)$ the compositum of $\kt_v$
and $k_v(x)$.

From section \ref{space} we know that the group $G_{x\,k_v}^{\circ}$
may be determined up to isomorphism solely from the index of the
orbit of $x$. In fact, if we define
\begin{equation}\label{Hxdefn}
H_{x \, k_v}=\begin{cases}
(k_v^{\times})^4 &\text{\qquad (sp),} \\
(k_v(x)^{\times})^2 &\text{\qquad (sp ur), (sp rm),} \\
(\kt_v^{\times})^2 &\text{\qquad (in), (rm), (in ur), (rm rm)*,} \\
\kt_v(x)^{\times} &\text{\qquad otherwise,}
\end{cases}
\end{equation}
for each of the various indices then $G_{x\,k_v}^{\circ}\cong
H_{x\, k_v}$ in all cases. We may regard $H_{x\, k_v}$ as the
$k_v$-points of an algebraic group $H_x$ defined over $\calo_v$ and we
shall do so below. 

As in section \ref{space}, if $k_v(x)/k_v$ is quadratic then we
shall write $\nu$ for the generator of $\gal(k_v(x)/k_v)$. If
$\kt_v(x)\neq\kt_v$ then $\nu$ may also be regarded as the
generator of $\gal(\kt_v(x)/\kt_v)$. Also the \emph{type} of $x\in
V_{k_v}^{\sst}$ will be the index attached to the orbit $G_{k_v}x$.

We wish to introduce parameterizations for the elements of the
stabilizer in the various cases.
If $x$ is a point of type (sp), we write
\begin{equation}\label{stabelm1} 
s_x(t_x) = (a(t_{11},t_{12}),a(t_{21},t_{22}),
a((t_{11}t_{21})^{-1},(t_{12}t_{22})^{-1}))\,,
\end{equation}
where $t_x = (t_{11},\dots,t_{22})\in (\mk_v)^4$.
Let $s_{x1}(t_x),s_{x2}(t_x),s_{x3}(t_x)$ be the 
three components of $s_x(t_x)$.   
If $x$ is a point of type (sp ur) or (sp rm), we write 
\begin{equation}\label{stabelm2} 
s_x(t_x) = (a(t_{11},t_{11}^{\nu}),a(t_{21},t_{21}^{\nu}),
a((t_{11}t_{21})^{-1},(t_{11}^{\nu}t_{21}^{\nu})^{-1}))\,,
\end{equation} 
where $t_x = (t_{11},t_{21})\in (k(x)_v^{\times})^2$.  
We use the notation $s_{x1}(t_x)$ et cetera for this case also.
If $x$ is a point of type (in) or (rm) then we write
\begin{equation}\label{stabelm5} 
s_x(t_x) = (a(t_{11},t_{12}),
a(\n_{\ti k_v/k_v}(t_{11}^{-1}),\n_{\ti k_v/k_v}(t_{12}^{-1}))\,,
\end{equation} 
where $t_x = (t_{11},t_{12})\in (\ti k_v^{\times})^2$.  
We use the notation $s_{x1}(t_x)$ et cetera for this case also.
If $x$ is a point of type (in ur) or (rm rm)* then we write
\begin{equation}\label{stabelm3} 
s_x(t_x) = (a(t_{11},t_{12}),a(t_{12}^{\sig},t_{11}^{\sig}),
a((t_{11}t_{12}^{\sig})^{-1},(t_{11}^{\sig}t_{12})^{-1}))\,,
\end{equation} 
where $t_x = (t_{11},t_{12})\in (\ti k_v^{\times})^2$.  
We use the notation $s_{x1}(t_x)$ et cetera for this case also.
Finally if $x$ is a point of type (in rm), (rm ur), 
(rm rm ur) or (rm rm rm) then we write 
\begin{equation}\label{stabelm4} 
s_x(t_x) = (a(t_{11},t_{11}^{\nu}),
a(\n_{\ti k_v(x)/k_v(x)}(t_{11}^{-1}),
\n_{\ti k_v(x)/k_v(x)}(t_{11}^{-1})^{\nu}))\,,
\end{equation} 
where $t_x = t_{11}\in \ti k_v(x)^{\times}$.  
We use the notation $s_{x1}(t_x)$ et cetera for this case also.
On $H_{x\, k_v}$ we define an invariant measure $dt_{x,v}$ as follows:
\begin{equation}
dt_{x,v}=\begin{cases}
\md t_{11v}\,\md t_{12v}\, \md t_{21v}\,\md t_{22v}
&\text{\qquad (sp),} \\
\md t_{11v}\,\md t_{21v}
&\text{\qquad (sp ur), (sp rm),} \\
\md t_{11v}\,\md t_{12v}
&\text{\qquad (in), (rm), (in ur), (rm rm)*,} \\
\md t_{11v}
&\text{\qquad otherwise.}
\end{cases}
\end{equation}
We note that if $v\in\gM_{\text{f}}$ then the volume of
$H_{x\, \calo_v}$ under this measure is $1$ in every case.

Suppose that $x\in V_{k_v}^{\sst}$ corresponds to a quadratic
extension of $k_v$. Then it is possible to choose an element
$g_x\in G_{k(x)}$ such that $x=g_x w$. Consider the following
condition on such an element.
\begin{cond}\label{ghomcond}
$g_x^{-1}g_x^{\nu}=(-\tau,-\tau,\tau)$ or $(-\tau,\tau)$.
\end{cond}
It is possible to find $g_x$ satisfying this condition for any $x$.
Indeed, $x=g_{xw_p}w_{p}$ for some $g_{xw_p}\in G_{k_v}$ and some
choice of $p$. Then $x=g_{xw_p}g_pw$ and $g_x=g_{xw_p}g_p\in
G_{k_v(x)}$ satisfies the condition.
\begin{prop}
If $g_x$ satisfies Condition \ref{ghomcond} then
\begin{equation}
G_{x\,k_v}^{\circ}=g_x\{s_x(t_x)\mid t_x\in H_{x\, k_v}\}g_x^{-1}\,.
\end{equation}
\end{prop}
\begin{proof}
We have $k_v(x)=k_v(w_p)$ for some $p$. Since $g_x$ and $g_p$ both
satisfy Condition \ref{ghomcond}, $g_xg_p^{-1}\in G_{k_v}$ and if
we put $h=g_xg_p^{-1}$ then $hw_p=x$ and so
$G_{x\,k_v}^{\circ}=hG_{w_p\,k_v}^{\circ}h^{-1}$. From section \ref{space},
\begin{equation}
G_{w_p\,k_v}^{\circ}=g_p\{s_x(t_x)\mid t_x\in H_{x\, k_v}\}g_p^{-1}
\end{equation}
and the conclusion follows.
\end{proof}
If $x=g_x w$ with $g_x\in G_{k_v}$ then we need not impose any
condition on $g_x$.

Suppose now that $g_x\in G_{k_v(x)}$, $x=g_xw$ and $g_x$ satisfies
Condition \ref{ghomcond} if $k_v(x)\neq k_v$. Then we can define an
isomorphism $\theta_{g_x}:G_{x\,k_v}^{\circ}\to H_{x\, k_v}$ by setting
$\theta_{g_x}(g_xs_x(t_x)g_x^{-1})=t_x$. If $g_{x1}$ and $g_{x2}$
are two such elements then let $h=g_{x2}g_{x1}^{-1}$. Using the
condition, we see that $h\in G_{x\,k_v}$. Also
\begin{equation}
\begin{aligned}
\theta_{g_{x1}}(g) &=
s_x^{-1}(g_{x1}^{-1}gg_{x1}) \\
&=s_x^{-1}(g_{x2}^{-1}hgh^{-1}g_{x2}) \\
&=\theta_{g_{x2}}(hgh^{-1})
\end{aligned}
\end{equation}
and so $\theta_{g_{x1}}$ and $\theta_{g_{x2}}$ differ by the
automorphism $g\mapsto hgh^{-1}$ of $G_{x\,k_v}^{\circ}$.
Since $G_{x\,k_v}^{\circ}$ is abelian and $G_{x\,k_v}/G_{x\,k_v}^{\circ}$
has order two, the automorphism $g\mapsto hgh^{-1}$ depends only on
the class of $h$ in $G_{x\,k_v}/G_{x\,k_v}^{\circ}$ and either is the
identity (if $h\in G_{x\,k_v}^{\circ}$) or squares to the identity
(if $h\notin G_{x\,k_v}^{\circ}$). In either case, this automorphism
is measure preserving and hence we may make the following
definition without ambiguity.
\begin{defn}\label{measure''}
Let $dg_{x,v}''=\theta_{g_x}^*(dt_{x,v})$ for any choice of $g_x\in
G_{k_v(x)}$ such that $g_xw=x$ and $g_x$ satisfies Condition
\ref{ghomcond} if $k_v(x)\neq k_v$.
\end{defn}
This establishes a choice of invariant measure on
$G_{x\,k_v}^{\circ}$ for each $x\in V_{k_v}^{\sst}$.

We have
\begin{equation}
\ti T_{k_v}=\begin{cases}
\{(t_1I_2,t_2I_2,(t_1t_2)^{-1}I_2)\}
&\text{\qquad in case (1),} \\
\{(t_1I_2,\n_{\kt_v/k_v}(t_1)^{-1}I_2)\}
&\text{\qquad in case (2),}
\end{cases}
\end{equation}
and so 
$\ti T_{k_v}\cong(k_v^{\times})^2$ in case (1) and $\ti
T_{k_v}\cong \kt_v^{\times}$ in case (2). We use the measure
\begin{equation}\label{tmeasure}
\md\ti t_v=\begin{cases}
\md t_{1v}\,\md t_{2v}
&\text{\qquad in case (1),} \\
\md t_{1v}
&\text{\qquad in case (2),}
\end{cases}
\end{equation}
on this group. We let $d\ti g_{x,v}''$ be the measure on 
$G_{x\,k_v}^{\circ}/\ti T_{k_v}$ such that $dg_{x,v}''=d\ti g_{x,v}''\,
\md \ti t_v$.

It is to achieve the following result that we have taken such pains
with the definition of the measures.
\begin{prop}\label{canonicity}
Suppose that $x,y\in V_{k_v}^{\text{\upshape ss}}$ and that $y=g_{xy}x$ for some
$g_{xy}\in G_{k_v}$. Let $i_{g_{xy}}:G_{y\,k_v}^{\circ}\to
G_{x\,k_v}^{\circ}$ be the isomorphism
$i_{g_{xy}}(g)=g_{xy}^{-1}gg_{xy}$. Then
\begin{equation}
dg_{y,v}''=i_{g_{xy}}^*(dg_{x,v}'')
\qquad \text{and}\qquad
d\ti g_{y,v}''=i_{g_{xy}}^*(d\ti g_{x,v}'')\,.
\end{equation}
\end{prop}
\begin{proof}
Let $g_x$ be chosen as above and put $g_y=g_{xy}g_x$. Then $g_y\in
G_{k_v(y)}=G_{k_v(x)}$, $g_yw=y$ and if $k_v(y)\neq k$ then
$g_y^{-1}g_y^{\nu}=g_x^{-1}g_x^{\nu}$, so that $g_y$ satisfies
Condition \ref{ghomcond} in this case. It follows that
\begin{equation}
\begin{aligned}
i_{g_{xy}}^*(dg_{x,v}'') &=
i_{g_{xy}}^*\theta_{g_x}^*(dt_{x,v}) \\
&=(\theta_{g_x}i_{g_{xy}})^*(dt_{x,v}) \\
&=\theta_{g_y}^*(dt_{y,v}) \\
&=dg_{y,v}''
\end{aligned}
\end{equation}
because $H_{x\, k_v}=H_{y\, k_v}$ and $dt_{x,v}=dt_{y,v}$. This
establishes the first claim and the second then follows from the
observation that $i_{g_{xy}}|_{\ti T_{k_v}}$ is the identity map.
\end{proof}
We choose a left invariant measure $dg _{x,v}'$ on 
$G_{k_v}/G^{\circ}_{x\,k_v}$ so that 
if $\Phi\in \cS(V_{k_v})$ then
\begin{equation}\label{lint}
\int_{G_{k_v}/G^{\circ}_{x\,k_v}} 
|P(g_{x,v}'x)|_v^s \Phi(g_{x,v}'x) dg_{x,v}' 
= \int_{G_{k_v}x} |P(y)|_v^{s-2} \Phi(y)dy \,,
\end{equation} 
where $dy$ is the Haar measure such that 
the volume of $V_{\co_v}$ is one if $v\in\gM_{\text{f}}$, Lebesgue
measure if $v\in\gM_{\R}$ and $2^8$ times Lebesgue measure if
$v\in\gM_{\C}$. This is possible because $|P(y)|_v^{-2}\,dy$ is a
$G_{k_v}$-invariant measure on $V_{k_v}^{\sst}$ and each of the
orbits $G_{k_v}x$ is an open set in $V_{k_v}^{\sst}$.
Note that 
\begin{equation} \label{oint}
\begin{aligned}
{} & \int_{G_{k_v}/G^{\circ}_{x\,k_v}} 
|\chi(g_{x,v}')|_v^s \Phi(g_{x,v}'x) dg_{x,v}'  \\
& = |P(x)|_v^{-s} \int_{G_{k_v}/G^{\circ}_{k_v}} 
|P(g_{x,v}'x)|_v^s \Phi(g_{x,v}'x) dg_{x,v}' 
\end{aligned}
\end{equation}
and
so, from (\ref{lint}),
this integral converges absolutely at least when $\re(s)>2$. 
If $g_{xy}\in G_{k_v}$ satisfies
$y=g_{xy}x$ and $i_{g_{xy}}$ is the inner automorphism $g\mapsto
g_{xy}^{-1}gg_{xy}$ of $G_{k_v}$ then
$i_{g_{xy}}(G_{y\,k_v}^{\circ})=G_{x\,k_v}^{\circ}$ and so $i_{g_{xy}}$
induces a map $i_{g_{xy}}: G_{k_v}/ G_{y\,k_v}^{\circ}\to
G_{k_v}/G_{x\,k_v}^{\circ}$. Since the integral on the right hand
side of (\ref{lint}) depends only on the orbit of $x$, it follows
that $i_{g_{xy}}^*(dg_{x,v}')=dg_{y,v}'$.
\begin{defn}\label{bxdefinition}
For $v\in\gM$ and $x\in V_{k_v}^{\sst}$ we let $b_{x,v}>0$ be the
constant verifying $dg_v=b_{x,v}\,dg_{x,v}'\,dg_{x,v}''$, where
$dg_v$ is the measure on $G_{k_v}$ chosen at the end of section
\ref{measure}.
\end{defn}
\begin{defn}\label{lzeta}  
For $\Phi\in \cS(V_{k_v})$ and $s\in\C$
we define
\begin{equation*}
\begin{aligned}
Z_{x,v}(\Phi,s) 
& = b_{x,v}  \int_{G_{k_v}/G^{\circ}_{x\,k_v}} 
|\chi(g_{x,v}')|_v^s \Phi(g_{x,v}'x) dg_{x,v}' \\
& = b_{x,v}|P(x)|_v^{-s} \int_{G_{k_v}x} |P(y)|_v^{s-2} \Phi(y)dy\,.
\end{aligned}
\end{equation*}
\end{defn}

\begin{prop}\label{bxindep} 
If {\rm $x,y\in V^{\text{ss}}_{k_v}$} and 
$G_{k_v}x=G_{k_v}y$ then $b_{x,v}=b_{y,v}$.  
\end{prop} 
\begin{proof}  Since the group $G_{k_v}$ is unimodular
$i_{g_{x,y}}^* dg_v = dg_v$.  
So 
\begin{equation*}
\begin{aligned}
dg_v & = b_{y,v} dg_{y,v}' dg''_{y,v} \\
& = b_{y,v} i_{g_{x,y}}^* dg_{x,v}'i_{g_{x,y}}^* dg_{x,v}''
& = b_{y,v}b_{x,v}^{-1} i_{g_{x,y}}^* dg_v \\
& = b_{y,v}b_{x,v}^{-1} dg_v\,. 
\end{aligned}
\end{equation*}
Therefore $b_{x,v}=b_{y,v}$.
\end{proof}
Let $d_{\text{pr}}g_x'' = 
\prod_v dg_{x,v}'',\; d_{\text{pr}}\ti g_x'' 
= \prod_v d\ti g_{x,v}''$
and $\md_{\text{pr}} \ti t = \prod_v \md \ti t_v$, 
where $\md \ti t_v$ is defined
in (\ref{tmeasure}).  
\begin{prop}\label{object}
Suppose {\rm $x\in V^{\sst}_k$} and $k(x)\not=k,\ti k$.  
Then, with respect to the measure {\rm $d_{\text{pr}}\ti g_x''$}, 
the volume of $G^{\circ}_{x\,\A}/\ti T_{\A}G^{\circ}_{x\,k}$ 
is ${2\gC_{\ti k(x)}}/{\gC_{\ti k}}$.
\end{prop}
\begin{proof}  Identifying $\ti T$ with $\gl(1)_{\ti k}$ and 
$G^{\circ}_x$ with $\gl(1)_{\ti k(x)}$, we define
$\ti T^1_{\A}$ (resp. $G^{\circ1}_{x\,\A}$) 
to be the set of ideles of $\ti k$ (resp. $\ti k(x)$)
with absolute value one.  Let $\md_{\text{pr}} \ti t^1$ and
$d_{\text{pr}}g_x''{}^1$ be the measures on 
$\ti T^1_{\A}$ and $G^{\circ1}_{x\,\A}$, such that 
$d_{\text{pr}}g_x'' = \md \lam d_{\text{pr}}g_x''{}^1$, 
$\md_{\text{pr}} \ti t = \md \lam \md_{\text{pr}} \ti t^1$ for  
\begin{equation*}
g_x'' = \lamb_{\ti k(x)} g_x''{}^1,\; 
\ti t = \lamb_{\ti k} \ti t^1
.\end{equation*}
Note that if $\lam\in \mr$ then
the absolute value of $\lamb_{\ti k}$ as an idele of 
$\ti k(x)$ is $\lam^2$.  
Therefore, $d_{\text{pr}}g_x'' = 2\md \lam d_{\text{pr}}g_x''{}^1$
for $g_x'' = \lamb_{\ti k}g_x''{}^1$.  
Since $d_{\text{pr}}g_x'' = 
\md _{\text{pr}}\ti td_{\text{pr}}\ti g_x''$ this implies that
$2d_{\text{pr}}g_x''{}^1 = 
\md_{\text{pr}} \ti t^1d_{\text{pr}}\ti g_x''$. 
So 
\begin{equation*}
\begin{aligned}
2 \int_{G^{\circ1}_{x\,\A}/G^{\circ}_{x\,k}} d_{\text{pr}}g_x''{}^1
& = \int_{G^{\circ1}_{x\,\A}/G^{\circ}_{x\,k}\ti T_{\A}^1} 
d_{\text{pr}}\ti g_x'' 
\int_{\ti T^1_{\A}/\ti T_k} \md_{\text{pr}} \ti t^1 \\
& = \vol(G^{\circ1}_{x\,\A}/\ti T_{\A}^1G^{\circ}_{x\,k}) 
\int_{\ti T^1_{\A}/\ti T_k} \md_{\text{pr}} \ti t^1 \,.
\end{aligned}
\end{equation*}
Since 
\begin{equation*}
\int_{G^{\circ1}_{x\,\A}/G^{\circ}_{x\,k}} 
d_{\text{pr}}\ti g_x''{}^1
= \gC_{\ti k(x)}\text{\qquad and \qquad}
\int_{\ti T^1_{\A}/\ti T_k} \md_{\text{pr}} \ti t^1 
= \gC_{\ti k}
\end{equation*}
this proves the proposition.  
\end{proof}

\section{A preliminary mean value theorem and \\ the formulation 
of its proof}\label{plan}

In this section we introduce the global zeta function of the
prehomogeneous vector space $(G,V)$ and recall from \cite{yukieh}
its most basic analytic properties. The zeta function is
approximately the Dirichlet generating series for the sequence
$\vol(G_{x\,\A}^{\circ}/\ti T_{\A}G_{x\,k}^{\circ})$. If it were
exactly this generating series then our work would be almost
complete, since Tauberian theory would allow us to extract the mean
value of the coefficients from the analytic behavior of the series.
Unfortunately, the actual zeta function contains an additional
factor in each term and we proceed to explain the filtering process
by which this difficulty may be surmounted. This leads us, on the
basis of a number of assumptions, to a preliminary form of the mean
value theorem that is our goal. The validity of these assumptions
is demonstrated in later sections. The final form of the theorem,
which differs from the preliminary form mostly in being more
explicit, is given in the next section.

Let $G_{\A}=G_{1\A}\times G_{2\A}$, let $dg_1$ and $dg_2$ be the
measures on $G_{1\A}$ and $G_{2\A}$ which were defined in section
\ref{measure} and put $dg=dg_1\,dg_2$ for $g=(g_1,g_2)$; this is a
Haar measure on $G_{\A}$. Write $\ti G=G/\ti T$, 
so that $V$ is a faithful representation of $\ti G$. 
Since $\ti T\cong \GL{1}_{\kt}$ as groups
over $k$, the first Galois cohomology group of $\ti T$ is trivial
and it follows that $\ti G_F\cong G_F/\ti T_F$ for any field
$F\supseteq k$. Thus $\ti G_{\A}\cong G_{\A}/\ti T_{\A}$ and $\ti
G_{\A}/\ti G_{k}\cong G_{\A}/\ti T_{\A}G_k$. Let
$\md_{\text{pr}}\ti t$ be the measure on $\ti T_{\A}$ defined
immediately before Proposition \ref{object}. Then $\md\ti
t=\gC_{\kt}^{-1}\md_{\text{pr}}\ti t$ is the measure on $\ti
T_{\A}$ compatible under the isomorphism $\ti
T_{\A}\cong\ti\A^{\times}$ with the measure defined on
$\ti\A^{\times}$ in section \ref{notation}. We choose the measure
$d\ti g$ on $\ti G_{\A}$ which satisfies $dg=d\ti g\,\md\ti t$.
Similarly, we choose the measure $d\ti g_v$ on $\ti G_{k_v}$ which
satisfies $dg_v=d\ti g_v\,\md\ti t_v$. Let $d_{\text{pr}}\ti
g=\prod_v d\ti g_v$. Using (\ref{glrelation}) we obtain
\begin{equation}\label{tigmeasure}
d\ti g=|\Del_{k}\Del_{\kt}|^{-1/2}\gC_{k}^{-2}\gC_{\kt}^{-1}
d_{\text{pr}}\ti g\,.
\end{equation}

\begin{defn}\label{gzeta}
Let $L_0=\{x\in V_k^{\sst}\mid k(x)\neq k,\kt\}$.
For $\Phi\in\cS(V_{\A})$ and $s\in\C$ we define 
\begin{equation*}
Z(\Phi,s) = \int_{G_{\A}/\ti T_{\A}G_k}
|\chi(\ti g)|^s \sum_{x\in L_0} \Phi(\ti gx) d\ti g\,.
\end{equation*}
\end{defn}
The integral $Z(\Phi,s)$ is called the \emph{global zeta function}
of $(G,V)$. It was proved in \cite{yukieh} that the integral
converges (absolutely and uniformly on compacta) if 
$\re(s)$ is sufficiently large.  
However, a slightly different formulation was used in
\cite{yukieh} and it is necessary to say a few words about the
translation from that paper to this.

The definition of the zeta function used in \cite{yukieh} is stated
in Definition (2.10) of that paper. For our purposes we shall
always take the character $\omega$ appearing there to be the
trivial character. The domain of integration used in \cite{yukieh}
is $\R_{+}\times G_{\A}^{0}/G_{k}$, where
$G_{\A}^{0}=G_{1\A}^{0}\times G_{2\A}^{0}$ is the set of elements
of $G_{\A}$ both of whose entries have determinant of idele norm
$1$. We have $(\R_{+}\times G_{\A}^{0})/\ti T_{\A}^{1}\cong \ti
G_{\A}$ via the map which sends the class of $(\lambda,g^{0})$ to
the class of $(1,c(\underline{\lambda}))g^{0}$. In \cite{yukieh},
$\R_{+}\times G_{\A}^{0}$ is made to act on $V_{\A}$ by requiring
that $(\lambda,1)$ acts by multiplication by $\underline{\lambda}$
and the above isomorphism is compatible with this.

We must compare
the measure $d\ti g$ on $\ti G_{\A}$ with the measure $\md\lambda\,
dg^{0}$ which was used in \cite{yukieh}. We have $\ti G_{\A}\cong
(\R_{+}^2\times G_{\A}^{0})/(\R_{+}\times\ti T_{\A}^{1})$ where
$\R_{+}\times\ti T_{\A}^{1}$ is included in $\R_{+}^2\times
G_{\A}^{0}$ via $(\lambda,\ti t)\mapsto(\lambda,\lambda^{-1},\ti
t)$ and $\R_{+}^2\times G_{\A}^{0}$ maps onto $\ti G_{\A}$ via
$(\lambda_1,\lambda_2,g^0)\mapsto
(c(\ti{\underline{\lambda}}_1),c(\underline{\lambda}_2))g^{0}\cdot\ti
T_{\A}$. In this quotient we have chosen the measure $d\ti g$ to be
compatible with the measures $4\md\lambda_1\,\md\lambda_2\,dg^0$ on
$\R_{+}^2\times G_{\A}^0$ and $\md\lambda\,\md\ti t^1$ on
$\R_{+}\times\ti T_{\A}^1$, where the volume of $\ti T_{\A}^{1}/\ti
T_{k}$ under $\md\ti t^1$ is $1$ (as in section \ref{notation}).
From this it follows that the measures $4\md\lambda\, dg^0$ and
$\md\ti t^1$ are compatible with the measure $d\ti g$ in the
quotient $(\R_{+}\times G_{\A}^0)/\ti T_{\A}^{1}\cong \ti G_{\A}$.

Furthermore, $|\chi(1,c(\underline{\lambda}))|=\lambda^4$ and so if
$Z^*(\Phi,s)$ denotes the zeta function studied in \cite{yukieh}
then we have $Z(\Phi,s)=4 Z^*(\Phi,4s)$. In \cite{yukieh},
Corollary 8.16 it is
shown that $Z^{*}(\Phi,s)$ has a meromorphic continuation to the
region $\re(s)>6$ with a simple pole at $s=8$ with residue
$\gV_{k}\gV_{\kt}\widehat{\Phi}(0)$. Thus we arrive at
\begin{thm}\label{globalanalyticproperties}
The zeta function $Z(\Phi,s)$ has a meromorphic continuation to the
region $\re(s)>3/2$ with a simple pole at $s=2$ with residue
$\gV_{k}\gV_{\kt}\widehat{\Phi}(0)$.
\end{thm}

Note that $\hp(0)$ is the Fourier transform of $\Phi$ evaluated at 
the origin and so is simply the integral of $\Phi$ over 
the $V_{\A}$.  We define $\Sig(\Phi) = \hp(0)$ for 
$\Phi\in \cS(V_{\A})$.  For $v\in\gM$ and $\Phi_v\in 
\cS(V_{k_v})$ we can define the local version of 
the distribution $\Sig(\Phi)$ by 
\begin{equation} 
\Sig_v(\Phi_v) = \int_{V_{k_v}} \Phi_v(y) dy 
.\end{equation}
Since the coordinate system of $V$ consists of 
four coordinates in $k$ and two coordinates in $\ti k$, 
if $\Phi=\otimes_v \Phi_v$ then
\begin{equation}  \label{distribution}
\Sig(\Phi) = |\Delta_k|^{-2}|\Delta_{\ti k}|^{-1}\prod_v 
\Sig_v(\Phi_v)
.\end{equation}

This completes our review of the analytic properties of the global
zeta function. Before we can rewrite $Z(\Phi,s)$ in a form which
makes this analytic information bear on the problem at hand we must
return briefly to the local situation.

Let $v\in\gM_{\text{f}}$. If $F/k_v$ is a quadratic extension then
$F$ is generated over $k_v$ by either of the roots of some
irreducible polynomial $p(z)=z^2+a_1z+a_2\in k_v[z]$. In fact, this
polynomial may always be chosen to satisfy the more stringent
condition that $\calo_F$ is generated over $\calo_v$ by either of
the roots of $p(z)$. If this condition is satisfied then the
discriminant of $p(z)$ generates the ideal $\Del_{F/k_v}$. We wish
to recall how this may be achieved in each case.

Recall that $p(z)\in k_v[z]$ is called an \emph{Eisenstein
polynomial} if $a_1\in\gp_v$ and $a_2\in\gp_v\setminus\gp_v^2$. If
$F/k_v$ is a ramified extension then there is always an Eisenstein
polynomial whose roots generate $F$ over $k_v$ and any such
polynomial will satisfy the stronger condition stated above. For
each $v\in\gM_{\text{f}}$, $k_v$ has a unique unramified quadratic
extension. If $F$ is this extension and
$v\notin\gM_{\text{dy}}$ then we may satisfy the
stronger condition simply by taking $p(z)$ with $a_1=0$ and $-a_2$
any non-square unit in $k_v$. If $v\in\gM_{\text{dy}}$ then we must
instead take $p(z)$ to be an \emph{Artin-Schreier polynomial},
which means, by definition, that $p(z)$ is irreducible in $k_v[z]$,
$a_1=-1$ and $a_2$ is a unit.  Note that $p$ stays irreducible 
modulo $\gp_v$ in this case by Hensel's lemma.  

For each $v\in\gM_{\text{f}}$ we choose a list of representatives
$w_{v,1},\dots, w_{v,N_v}$, one for each of the $G_{k_v}$-orbits in
$V_{k_v}^{\sst}$, in such a way that $P(w_{v,i})$ generates the
ideal $\Del_{k(w_{v,i})/k_v}$ for $i=1,\dots, N_v$. This is
possible, in light of the previous paragraph, if we take each
$w_{v,i}$ to equal $w_p$ for a suitable $p(z)\in k_v[z]$. In the
special case where $k(w_{v,i})=k_v$ we take $w_{v,i}=w_p$ for
$p(z)=z^2-z$. For $v\in\gM_{\infty}$ we require instead that
$|P(w_{v,i})|_v=1$ for $i=1,\dots,N_v$, which is clearly possible.
In both cases we assume for convenience that $w_{v,1}$ represents
the orbit corresponding to $k_v$ itself. This done, if $F/k$ is a
quadratic extension then let $\om_{v,i_v(F)}$ represent the orbit
corresponding to $F_v/k_v$ (with $i_v(F)=1$ if $v$ splits in $F$).
Then we have
\begin{equation}\label{idealnormfactor}
\esN(\Del_{F/k})^{-1}=\prod_{v\in\gM_{\text{f}}}
\esN_{v}(\Del_{F/k,v})^{-1}
=\prod_{v\in\gM_{\text{f}}}|P(w_{v,i_v(F)})|_v
=\prod_{v\in\gM}|P(w_{v,i_v(F)})|_v\,.
\end{equation}

For $x\in L_0$ and $\Phi=\otimes\Phi_v\in\cS(V_{\A})$ we define the
\emph{orbital zeta function} of $x$ to be
$Z_{x}(\Phi,s)=\prod_{v\in\gM} Z_{x,v}(\Phi_v,s)$. If $x$ lies in
the orbit of $w_{v,i}$ in $V_{k_v}^{\sst}$ then we shall write
$\Xi_{x,v}(\Phi_v,s)= Z_{w_{v,i},v}(\Phi_v,s)$ and 
$\Xi_{x}(\Phi,s)=\prod_{v\in\gM}\Xi_{x,v}(\Phi_v,s)$. 
We call $\Xi_{x,v}(\Phi_v,s)$ the 
\emph{standard local zeta function} 
and  $\Xi_x(\Phi_v,s)$ the 
\emph{standard orbital zeta function}.  
\begin{prop}\label{ZetatoXiconversion}
For $x\in L_0$ and $\Phi=\otimes\Phi_v\in\cS(V_{\A})$
we have 
\begin{equation*}
Z_x(\Phi,s)=\esN(\Del_{k(x)/k})^{-s}\Xi_{x}(\Phi,s)\,.
\end{equation*}
\end{prop}
\begin{proof}
For each $v\in\gM$ let $i_v(x)$ be such that $x\in
G_{k_v}w_{v,i_v(x)}$. Then, from (\ref{lzeta}),
\begin{equation}
\begin{aligned}\label{ZetaXiequation}
Z_{x,v}(\Phi_v,s) &=
b_{x,v}|P(x)|_v^{-s}\int_{G_{k_v}x}|P(y)|_v^{s-2}\Phi_v(y)\,dy \\
&=\frac{|P(w_{v,i_v(x)})|_v^s}{|P(x)|_v^s}\cdot
\frac{b_{w_{v,i_v(x)},v}}{|P(w_{v,i_v(x)})|_v^{s}}
\int_{G_{k_v}w_{v,i_v(x)}}|P(y)|_v^{s-2}\Phi_v(y)\,dy \\
&=\frac{|P(w_{v,i_v(x)})|_v^s}{|P(x)|_v^s}\cdot
Z_{w_{v,i_v(x)},v}(\Phi_v,s) \\
&=\frac{|P(w_{v,i_v(x)})|_v^s}{|P(x)|_v^s}\cdot
\Xi_{x,v}(\Phi_v,s)
\end{aligned}
\end{equation}
where we have used Proposition \ref{bxindep} in passing from the first
line to the second. 
Applying (\ref{idealnormfactor}) to $F=k(x)$ we find
that $\prod_{v\in\gM}|P(w_{v,i_v(x)})|_v^s=\esN(\Del_{k(x)/k})^{-s}$. 
Since $x\in V_{k}^{\sst}$, $P(x)\in k^{\times}$ and so the Artin
product formula implies that $\prod_{v\in\gM}|P(x)|_v=1$. Now taking
the product over all $v\in\gM$ on both sides of (\ref{ZetaXiequation})
proves the identity. \end{proof}

For convenience, we introduce the abbreviation
\begin{equation}
\esR_1 = |\Delta_k|^{-{\frac 12}}
|\Delta_{\ti k}|^{-{\frac 12}}
\gC_k^{-2}\gC_{\ti k}^{-2}\,.
\end{equation}
\begin{prop}\label{ZetaSeries}
If $\Phi=\otimes\Phi_v\in\cS(V_{\A})$ then we have
\begin{equation*}
Z(\Phi,s)=\esR_1\sum_{x\in G_k\backslash L_0}
\esN(\Del_{k(x)/k})^{-s}\gC_{\kt(x)}\Xi_x(\Phi,s)\,.
\end{equation*}
\end{prop}
\begin{proof}
From Definition \ref{gzeta} we have
{\allowdisplaybreaks
\begin{align*}
Z(\Phi,s)&=\sum_{x\in G_k\backslash L_0}
\int_{G_{\A}/\ti T_{\A}G_k} |\chi(\ti g)|^s
\sum_{\gam\in G_k/G_{x\,k}}
\Phi(\ti g\gam x)\,d\ti g \\
&=\sum_{x\in G_k\backslash L_0}
\int_{G_{\A}/\ti T_{\A}G_{x\,k}}
|\chi(\ti g)|^s\Phi(\ti g x)\, d\ti g \\
&=\tfrac12\sum_{x\in G_k\backslash L_0}
\int_{G_{\A}/\ti T_{\A} G_{x\,k}^{\circ}}
|\chi(\ti g)|^s\Phi(\ti g x)\,d\ti g\text{\quad since
$[G_{x\,k}:G_{x\,k}^{\circ}]=2$} \\
&=\tfrac12\esR_1\gC_{\kt}
\sum_{x\in G_k\backslash L_0}
\int_{G_{\A}/\ti T_{\A}G_{x\,k}^{\circ}}|\chi(\ti g)|^s\Phi(\ti g
x)\,d_{\text{pr}}\ti g
\text{\quad by (\ref{tigmeasure})} \\
&=
\tfrac12\esR_1\gC_{\kt}
\sum_{x\in G_k\backslash L_0}
\bigl(\prod_v b_{x,v}\bigr)
\int_{G_{\A}/G_{x\,\A}^{\circ}}
|\chi(\ti g')|^s\Phi(\ti g'x)\, d_{\text{pr}}\ti g' \\
&\qquad\qquad\qquad\qquad\qquad\qquad\qquad\cdot
\int_{G_{x\,\A}^{\circ}/\ti T_{\A} G_{x\,k}^{\circ}}\,d_{\text{pr}}\ti
g'' \text{\quad by Definition \ref{bxdefinition}} \\
&=\tfrac12\esR_1\gC_{\kt}
\sum_{x\in G_k\backslash L_0}
\bigl(\prod_v Z_{x,v}(\Phi_v,s)\bigr)\cdot
\vol(G_{x\,\A}^{\circ}/\ti T_{\A}G_{x\,k}^{\circ})
\text{\quad by Definition \ref{lzeta}} \\
&=\esR_1
\sum_{x\in G_k\backslash L_0} Z_x(\Phi,s)\gC_{\kt(x)}
\text{\quad by Proposition \ref{object}} \\
&=\esR_1
\sum_{x\in G_k\backslash L_0}
\esN(\Del_{k(x)/k})^{-s}\gC_{\kt(x)}
\Xi_x(\Phi,s)
\text{\quad by Proposition \ref{ZetatoXiconversion}.}
\end{align*}}
\end{proof}

We are now ready to describe the filtering process. This process
was originally used in \cite{dawrb} and is described in
a general setting in \cite{yukiec}, \S0.5. Our discussion
will follow this latter reference, but with simplifications arising
from the fact that we know the residue of the global zeta function
explicitly (by Theorem \ref{globalanalyticproperties}).  

We set $S_0=\gM_{\infty}\cup\gM_{\text{rm}}\cup\gM_{\text{dy}}$ and
fix a finite set $S\supseteq S_0$ of places of $k$. For each finite
subset $T\supseteq S$ of $\gM$ we consider $T$-tuples
$\omega_T=(\omega_v)_{v\in T}$ where each $\omega_v$ is one of the
standard orbital representatives, $w_{v,i}$, for the orbits in
$V_{k_v}^{\sst}$ chosen above. If $x\in V_{k}^{\sst}$ and $x\in
G_{k_v}\omega_v$ then we write $x\approx \omega_v$ and if $x\approx
\omega_v$ for all $v\in T$ then we write $x\approx \omega_{T}$.

For later purposes, it is convenient to make the following definition.  
\begin{defn}
For any $v\in\gMf$,  $\Phi_{v,0}$ 
is the characteristic function of $V_{\calo_v}$. 
\end{defn}

Let $\Xi_{x,v}(s)=\Xi_{x,v}(\Phi_{v,0},s)$ and
$\Xi_{x,T}(s)=\prod_{v\notin T}\Xi_{x,v}(s)$. From the integral
defining $\Xi_{x,v}(s)$ it follows that for $v\notin S_0$ this
function may be expressed as $\Xi_{x,v}(s)=\sum_{n=-\infty}^{\infty}
a_{x,v,n}q_v^{-ns}$ for certain numerical coefficients $a_{x,v,n}$.
In section \ref{sec-omega} we shall establish the following
condition.
\begin{cond}\label{xicoeff}
For all $v\notin S_0$ and all $x\in V_{k_v}^{\text{\upshape ss}}$
we have $a_{x,v,n}=0$ for $n<0$, $a_{x,v,0}=1$ 
and $a_{x,v,n}\geq 0$ for all $n$.
\end{cond}

Suppose that we have Dirichlet series
$L_i(s)=\sum_{m=1}^{\infty}\ell_{i,m}m^{-s}$ for $i=1,2$. If
$\ell_{1,m}\leq\ell_{2,m}$ for all $m\geq1$ then we shall write
$L_1(s)\cleq L_2(s)$. In section \ref{sec-estimate} we shall
establish that for every $v\notin S_0$ there exists a
Dirichlet series $L_v(s)=\sum_{n=0}^{\infty}\ell_{v,n}q_v^{-ns}$
which satisfies the following condition.
\begin{cond}\label{aproperty}
\begin{itemize}
\item[(1)] For all $v\notin S_0$ and $x\in V_{k_v}^{\text{\upshape ss}}$,
$\Xi_{x,v}(s)\cleq L_v(s)$.
\item[(2)] The series defining $L_v(s)$ converges to a holomorphic
function in the region $\re(s)>1$ and the product $\prod_{v\notin
S_0}L_v(s)$ converges absolutely and locally uniformly in the
region $\re(s)>3/2$.
\item[(3)] For all $v\notin S_0$, $\ell_{v,0}=1$ and $\ell_{v,n}\geq 0$
for all $n$.
\end{itemize}
\end{cond}
For any $T\supseteq S$ we define $L_T(s)=\prod_{v\notin T}L_v(s)$.
Both $\Xi_{x,T}(s)$ and $L_T(s)$ are Dirichlet series and if we let
\begin{equation}
\Xi_{x,T}(s)=\sum_{m=1}^{\infty}a_{x,T,m}^*m^{-s}
\text{\qquad and \qquad}
L_T(s)=\sum_{m=1}^{\infty}\ell_{T,m}^*m^{-s}
\end{equation}
then $a_{x,T,m}^*$ (resp. $\ell_{T,m}^*$) is the sum of the terms
$\prod_{v\notin T}a_{x,v,n_v}$ (resp. $\prod_{v\notin
T}\ell_{v,n_v}$) over all possible factorizations $m=\prod_{v\notin
T}q_v^{n_v}$. Since only finitely-many places, $v$, of $k$ can have
$q_v$ equal to a power of a particular prime, the number of such
factorizations is finite. Also, in any such factorization,
$n_v=0$ for all but finitely-many $v$, and so this sum is well-defined.
It follows
from Conditions \ref{xicoeff} and \ref{aproperty} that $0\leq
a_{x,T,m}^*\leq\ell_{T,m}^*$ and $a_{x,T,1}^*=1$ for all $x\in
V_k^{\sst}$, all $T\supseteq S$ and all $m\geq1$. We shall use
these observations in the proof of Theorem \ref{prelim-main}
below.

We define
\begin{equation}\label{xisdefn}
\xi_{\omega_T}(s)=\sum_{x\in G_k\backslash L_0,x\approx\omega_T}
\esN(\Del_{k(x)/k})^{-s}\gC_{\kt(x)}\Xi_{x,T}(s)
\end{equation}
and
\begin{equation}
\xi_{\omega_S,T}(s)=\sum_{x\in G_k\backslash L_0,x\approx\omega_S}
\esN(\Del_{k(x)/k})^{-s}\gC_{\kt(x)}
\Xi_{x,T}(s)\,,
\end{equation}
which is the sum of $\xi_{\omega_T}(s)$ over all
$\omega_T=(\omega_v)_{v\in T}$ which extend the fixed $S$-tuple
$\omega_S$. In order to determine the analytic properties of these
Dirichlet series we require the following result.

\begin{lem}\label{testfunction}
Let $v\in\gM$, $x\in V_{k_v}^{\text{\upshape ss}}$ and $r\in\C$.
Then there exists $\Phi_v\in\cS(V_{k_v})$ such that the support of
$\Phi_v$ is contained in $G_{k_v}x$, $Z_{x,v}(\Phi_v,s)$ is an
entire function and $Z_{x,v}(\Phi_v,r)\neq0$.
\end{lem}
\begin{proof}
The set $G_{k_v}x$ is open and $y\mapsto |P(y)|_v^{r-2}$ is a
continuous function on it. We may therefore find an open set $U$
containing $x$, having compact closure $\bar{U}\subseteq G_{k_v}x$
and such that
\begin{equation}\label{easy-inequality}
\bigl| |P(y)|_v^{r-2}-|P(x)|_v^{r-2}\bigr|< \tfrac12|P(x)|_v^{r-2}
\end{equation}
for $y\in\bar{U}$.
We can then choose $\Phi_v\in\cS(V_{k_v})$ in such a way that
$\supp(\Phi_v)\subseteq\bar{U}$ and
$\int_{\bar{U}}\Phi_v(y)\,dy=1$. Now (\ref{easy-inequality})
implies that $|P(y)|_v$ doesn't vanish on $\bar{U}$ and hence it is
bounded both above and below by positive constants on this
compactum. Thus $Z_{x,v}(\Phi_v,s)$ is entire. The inequality
(\ref{easy-inequality}) also implies that
\begin{equation*}
\bigl|Z_{x,v}(\Phi_v,r)-b_{x,v}|P(x)|_v^{-2}\bigr|\leq
\tfrac12 b_{x,v}|P(x)|_v^{-2}
\end{equation*}
and hence $Z_{x,v}(\Phi_v,r)\neq0$.
\end{proof}

\begin{prop}\label{xisresidue}
Let $T\supseteq S$ be a finite set of places of $k$ and $\omega_T$
be a $T$-tuple, as above. The Dirichlet series $\xi_{\omega_T}(s)$
has a meromorphic continuation to the region $\re(s)>3/2$. Its only
possible singularity in this region is a simple pole at $s=2$ with
residue
\begin{equation*}
\esR_2\prod_{v\in T}b_{\omega_v,v}^{-1}|P(\omega_v)|^2_v
\end{equation*}
where
\begin{equation*}
\esR_2=\res_{s=1}\zeta_k(s)\cdot\res_{s=1}\zeta_{\kt}(s)\cdot
Z_{k}(2)Z_{\kt}(2)/|\Del_{k}|\,.
\end{equation*}
\end{prop}
\begin{proof}
For each $v\in T$ we choose $\Phi_v\in\cS(V_{k_v})$ such that
$\supp(\Phi_v)\subseteq G_{k_v}\omega_v$. Let
$\Phi=\bigotimes_{v\in T}\Phi_v\otimes\bigotimes_{v\notin
T}\Phi_{v,0}\in\cS(V_{\A})$. For $v\in T$ we have
$\Xi_{x,v}(\Phi_v,s)=0$ unless $x\approx\omega_v$ and hence
\begin{align*}
Z(\Phi,s)&=\esR_1\big(\prod_{v\in T}\Xi_{\omega_v,v}(\Phi_v,s)\big)
\sum_{x\in G_k\backslash L_0,x\approx\omega_T}
\esN(\Del_{k(x)/k})^{-s}\gC_{\kt(x)}\Xi_{x,T}(s) \\
&=\esR_1\big(\prod_{v\in
T}\Xi_{\omega_v,v}(\Phi_v,s)\big)\xi_{\omega_T}(s)
\end{align*}
by Proposition \ref{ZetaSeries}. Using Lemma \ref{testfunction}
and Theorem \ref{globalanalyticproperties}, this formula implies
the first statement.

Now choose $\Phi_v$ for $v\in T$ so that
$\Xi_{\omega_v,v}(\Phi_v,2)\neq0$. It follows directly from the
definition that
$\Xi_{\omega_v,v}(\Phi_v,2)=b_{\omega_v,v}|P(\omega_v)|_v^{-2}
\Sigma_v(\Phi_v)$ for all $v\in T$ and so the residue of
$\xi_{\omega_T}(s)$ at $s=2$ is
\begin{equation*}
\esR_1^{-1}\big(\prod_{v\in T}
b_{\omega_v,v}^{-1}|P(\omega_v)|^2_v\big)
\big(\prod_{v\in T}\Sigma_v(\Phi_v)\big)^{-1}\res_{s=2}Z(\Phi,s)\,.
\end{equation*}
We have $\Sigma_v(\Phi_{v,0})=1$ for $v\notin T$ and hence
\begin{equation*}
\res_{s=2}Z(\Phi,s)=\gV_{k}\gV_{\kt}|\Del_{k}|^{-2}|
\Del_{\kt}|^{-1}\prod_{v\in T}\Sigma_v(\Phi_v)\,.
\end{equation*}
Combining the last two equations shows that the residue of
$\xi_{\omega_T}(s)$ at $s=2$ is
\begin{equation*}
\esR_{1}^{-1}\gV_{k}\gV_{\kt}|\Del_{k}|^{-2}
|\Del_{\kt}|^{-1}
\big(\prod_{v\in T}b_{\omega_v,v}^{-1}|P(\omega_v)|^2_v\big)
\end{equation*}
and using the definition of $\esR_1$ and the values of $\gV_k$ and
$\gV_{\kt}$ (see the end of section \ref{measure}) gives the second
claim.
\end{proof}

\begin{cor}\label{xiSTresidue}
The Dirichlet series $\xi_{\omega_S,T}(s)$ has a meromorphic
continuation to the region $\re(s)>3/2$. Its only possible
singularity in this region is a simple pole at $s=2$ with residue
\begin{equation*}
\esR_2\big(\prod_{v\in S}
b_{\omega_v,v}^{-1}|P(\omega_v)|_v^2\big)\cdot
\prod_{v\in T\setminus S}
\sum_{x}\big(b_{x,v}^{-1}|P(x)|^2_v\big)
\end{equation*}
where the sum is over the complete set, $\{x\}$, of standard orbit
representatives for $G_{k_v}\backslash V_{k_v}^{\text{\upshape ss}}
$.
\end{cor}
\begin{proof}
We have $\xi_{\omega_S,T}(s)=\sum_{\omega_T}\xi_{\omega_T}(s)$
where the sum is over all $T$-tuples $\omega_T$ which extend the
$S$-tuple $\omega_S$. The claim follows immediately.
\end{proof}

We let $E_v=\sum_{x}b_{x,v}^{-1}|P(x)|_v^2$ for $v\notin S_0$ where
the sum is over all standard representatives, $x$, for orbits in
$G_{k_v}\backslash V_{k_v}^{\sst}$. 
In section \ref{final} we shall prove 
that the following condition holds.

\begin{cond}\label{evcond}
The product $\prod_{v\notin S_0}E_v$ converges to a positive
number.
\end{cond}

We are now ready to state and prove, subject to Conditions
\ref{xicoeff}, \ref{aproperty} and \ref{evcond}, the theorem
which is the goal of this section.

\begin{thm}\label{prelim-main}
Let $S\supseteq S_0$ be a finite set of places of $k$ and
$\omega_S$ be an $S$-tuple of standard orbital representatives.
Then
\begin{equation*}
\lim_{X\to\infty}X^{-2}
\sum_{\substack{x\in G_k\backslash L_0,x\approx\omega_S \\
\esN(\Del_{k(x)/k})\leq X}}\gC_{\kt(x)}=
\tfrac12\esR_2\prod_{v\in S}\big(b_{\omega_v,v}^{-1}
|P(\omega_v)|^2_v\big)\cdot
\prod_{v\notin S} E_v\,.
\end{equation*}
\end{thm}
\begin{proof}
In the following, sums over $x$ will be understood to include the
conditions $x\in G_k\backslash L_0$ and $x\approx \omega_S$ as well
as any further conditions which may be explicitly imposed. We have
$\xi_{\omega_S,T}(s)=\sum_{m=1}^{\infty}c_mm^{-s}$ where
\begin{equation*}
c_m=\sum_{x,n,\esN(\Del_{k(x)/k})n=m}
\gC_{\kt(x)}a_{x,T,n}^*\,.
\end{equation*}
Applying the Tauberian theorem (\cite{narkiewicz}, p. 464, Theorem
I) to $\xi_{\omega_S,T}(s)$ we obtain, in light of Corollary
\ref{xiSTresidue},
\begin{equation*}
\lim_{X\to\infty}X^{-2}
\sum_{x,n,\esN(\Del_{k(x)/k})n\leq X}
\gC_{\kt(x)}a_{x,T,n}^*=
\tfrac12\esR_2\big(
\prod_{v\in S}b_{\omega_v,v}^{-1}|P(\omega_v)|_v^{2}\big)\cdot
\prod_{v\in T\setminus S}E_v\,.
\end{equation*}
We shall denote the right-hand side of this equation by $\esL_T$.
Note that $\esL=\lim_{T\to\gM}\esL_T$ is the right-hand side of the
equation in the statement. Since $a_{x,T,n}^*\geq0$ for all $n$ and
$a_{x,T,1}^{*}=1$ we obtain
\begin{equation*}
\varlimsup_{X\to\infty} X^{-2}\sum_{\esN(\Del_{k(x)/k})\leq X}
\gC_{\kt(x)}
\leq\esL_T
\end{equation*}
for all $T$ and so $\varlimsup_{X\to\infty} X^{-2}
\sum_{\esN(\Del_{k(x)/k})\leq X}\gC_{\kt(x)}\leq\esL$. 
It follows that there is a constant $C$ such that
$\sum_{\esN(\Del_{k(x)/k})\leq X}\gC_{\kt(x)}\leq CX^2$ for all
$X>0$ (note that if $X<1$ then the sum is $0$).
Furthermore,
\begin{align*}
\sum_{\esN(\Del_{k(x)/k})\leq X}\gC_{\kt(x)}&=
\sum_{\esN(\Del_{k(x)/k})n\leq X}
\gC_{\kt(x)}a_{x,T,n}^*-
\sum_{\esN(\Del_{k(x)/k})n\leq X,n\geq2}
\gC_{\kt(x)}a_{x,T,n}^* \\
&\geq
\sum_{\esN(\Del_{k(x)/k})n\leq X}\gC_{\kt(x)}a_{x,T,n}^*-
\sum_{\esN(\Del_{k(x)/k})n\leq X,n\geq2}
\gC_{\kt(x)}\ell_{T,n}^* \\
&=\sum_{\esN(\Del_{k(x)/k})n\leq X}
\gC_{\kt(x)}a_{x,T,n}^*-
\sum_{n=2}^{\infty}\ell_{T,n}^*
\sum_{\esN(\Del_{k(x)/k})\leq X/n}\gC_{\kt(x)} \\
&\geq\sum_{\esN(\Del_{k(x)/k})n\leq X}
\gC_{\kt(x)}a_{x,T,n}^*-
CX^2\sum_{n=2}^{\infty}\ell_{T,n}^*n^{-2} \\
&=\sum_{\esN(\Del_{k(x)/k})n\leq X}\gC_{\kt(x)}a_{x,T,n}^*-
CX^2(L_T(2)-1)\,.
\end{align*}
It follows that, for all $T\supseteq S$,
\begin{equation*}
\varliminf_{X\to \infty}X^{-2}
\sum_{\esN(\Del_{k(x)/k})\leq X}\gC_{\kt(x)}\geq
\esL_T-C(L_T(2)-1)
\end{equation*}
and letting $T\to\gM$ we obtain
\begin{equation*}
\varliminf_{X\to\infty}X^{-2}\sum_{\esN(\Del_{k(x)/k})\leq X}
\gC_{\kt(x)}\geq\esL
\end{equation*}
since $\lim_{T\to\gM}L_T(2)=1$. The theorem follows. \end{proof}

The remainder of this paper and its companion
\cite{kable-yukie-pbh-II} are devoted to verifying the conditions
enunciated in this section and to evaluating the constants which
appear in Theorem \ref{prelim-main}. In the next section we make
use of the results of this work to state the theorem in a more
explicit form.

\section{The Mean Value Theorem}\label{final}

In this section we shall derive a more explicit and convenient mean
value theorem from Theorem \ref{prelim-main}. Throughout, $k$
will be a number field and $\kt$ a fixed quadratic extension of
$k$. If $F_1$ and $F_2$ are distinct quadratic extensions of $k$,
neither equal to $\kt$, then we shall call $F_1$ and $F_2$
\emph{paired} (with respect to $\kt$) if $F_2\subseteq F_1\cdot
\kt$. Since this condition uniquely determines $F_2$ from $F_1$, we
may write $F_2=F_1^*$ if $F_2$ and $F_1$ are paired. Our first
result will be used below to express $\gC_{\kt(x)}$ in terms of
$\gC_{k(x)}$ and $\gC_{k(x)^*}$ for $x\in L_0$.
\begin{prop}\label{C-factorization}
Suppose that $L/k$ is a biquadratic extension of number fields and
that $k_1$, $k_2$ and $k_3$ are the quadratic extensions of $k$
contained in $L$. Then
$\mathfrak{C}_{L}=\mathfrak{C}_{k}^{-2}
\mathfrak{C}_{k_1}
\mathfrak{C}_{k_2}
\mathfrak{C}_{k_3}$.
\end{prop}
\begin{proof}
This identity is perhaps the simplest instance of what is
known as a Brauer relation (see \cite{brauer-classrel}, p. 162, for
instance). 
For the reader's convenience we sketch the proof from the theory of
the Dedekind zeta function.
Using Theorem 1.1, Chapter XII of \cite{lang-algebraic-number-theory},
p. 230 we have the factorization
\begin{equation*}
\zeta_L(s)=\zeta_k(s)L(s,\chi_1)L(s,\chi_2)L(s,\chi_3)\,,
\end{equation*}
where $\chi_j$ is the idele class character of $k$ corresponding by
class field theory to $k_j$. Multiplying both sides of this
identity by $\zeta_k(s)^2$ we obtain
\begin{equation}\label{eq:zeta-identity}
\zeta_L(s)\zeta_k(s)^2=\zeta_{k_1}(s)\zeta_{k_2}(s)\zeta_{k_3}(s)\,.
\end{equation}
Since $\res_{s=1}\zeta_{F}(s)=\mathfrak{C}_F/|\Del_F|^{1/2}$, 
it follows that
\begin{equation*}
\mathfrak{C}_{L}\mathfrak{C}_{k}^2|\Del_L\Del_k^2|^{-1/2}
=\mathfrak{C}_{k_1}
\mathfrak{C}_{k_2}
\mathfrak{C}_{k_3}|\Del_{k_1}\Del_{k_2}\Del_{k_3}|^{-1/2}\,.
\end{equation*}
Recall that we have a functional equation
\begin{equation*}
\zeta_F(1-s)=(2^{-2r_2(F)}\pi^{-[F:\mathbb{Q}]}|\Del_F|)^{s-1/2}
\frac{\Gamma(s/2)^{r_1(F)}\Gamma(s)^{r_2(F)}}%
{\Gamma((1-s)/2)^{r_1(F)}\Gamma(1-s)^{r_2(F)}}\zeta_F(s)\,,
\end{equation*}
where $r_1(F)$ denotes the number of real places of $F$ and
$r_2(F)$ the number of complex places of $F$.
It is easy to check that $[L:\mathbb{Q}]+2[k:\mathbb{Q}]=
\sum_{j=1}^3[k_j:\mathbb{Q}]$ and $r_i(L)+2r_i(k)=\sum_{j=1}^3
r_i(k_j)$ for $i=1,2$. Comparing the factors in the functional
equation on both sides of (\ref{eq:zeta-identity}) now shows that
$|\Del_L\Del_k^2|=|\Del_{k_1}\Del_{k_2}\Del_{k_3}|$
and the identity follows. 
\end{proof}

For notational compactness we shall set $\vep_v(x)=
b_{x,v}^{-1}|P(x)|_v^2$ for all $v\in\gM$ and $x\in
V_{k_v}^{\sst}$. These constants are related to the quantities
calculated in later sections by the following result.
\begin{lem}
Let $v\in\gM_{\text{\upshape f}}$ and $x\in
V_{k_v}^{\text{\upshape ss}}$. Then
\begin{equation*}
\vep_v(x)=\vol(K_v\cap G_{x\,k_v}^{\circ})
\vol(K_vx)
\end{equation*}
where the first volume is evaluated with respect to the canonical
measure $dg_{x,v}''$ on $G_{x\, k_v}^{\circ}$ and the second with
respect to the measure on $V_{k_v}$ under which $V_{\calo_v}$ has
volume $1$.
\end{lem}
\begin{proof}
We have
\begin{align*}
1&=\int_{K_v}dg_v \\
&=b_{x,v}\int_{K_v G_{x\,k_v}^{\circ}/G_{x\,k_v}^{\circ}}
dg_{x,v}'\cdot
\int_{K_v\cap G_{x\,k_v}^{\circ}} dg_{x,v}''
\text{\quad by Definition \ref{bxdefinition}} \\
&=b_{x,v}\vol(K_v\cap G_{x\,k_v}^{\circ})
\int_{K_vx}|P(y)|_v^{-2}\,dy
\end{align*}
by (\ref{lint}) with $s=0$ and $\Phi$ the characteristic function
of $K_vx$. But $|P(y)|_v=|P(x)|_v$ for all $y\in K_vx$ and so
$1=b_{x,v}\vol(K_v\cap G_{x\,k_v}^{\circ})
|P(x)|_v^{-2}\vol(K_vx)$.
\end{proof}
Using this formula for $\vep_v(x)$ and the results of sections
\ref{sec-iv} and \ref{bxsif} we may determine the values of
$\vep_v(x)$ for all $v\notin\gM_{\text{dy}}\cup \gM_{\infty}$ and
all standard orbital representatives $x\in V_{k_v}^{\sst}$. We
record the results in Table \ref{table-nondyadic} below.

\newcommand{\tallstrut}{\rule[-7pt]{0cm}{22pt}}
\newcommand{\tinsert}{\quad\tallstrut}
\begin{table}
\begin{center}
\begin{tabular}{|c|c|}
\hline
\ Index\ &$\tinsert\vep_v(x)\tinsert$ \\
\hline\hline
(sp)&$\tinsert\tfrac12(1+q_v^{-1})(1-q_v^{-2})^2\tinsert$ \\
\hline
(in)&$\tinsert\tfrac12(1-q_v^{-1})(1-q_v^{-4})\tinsert$ \\
\hline
(rm)&$\tinsert\tfrac12(1-q_v^{-2})^2\tinsert$ \\
\hline
(sp ur)&$\tinsert\tfrac12(1-q_v^{-1})^3(1-q_v^{-2})\tinsert$ \\
\hline
(sp rm)&$\tinsert\tfrac12q_v^{-1}(1-q_v^{-1})(1-q_v^{-2})^3
\tinsert$ \\
\hline
(in ur)&$\tinsert\tfrac12(1-q_v^{-1})(1-q_v^{-4})\tinsert$ \\
\hline
(in rm)&$\tinsert\tfrac12q_v^{-1}(1-q_v^{-1})(1-q_v^{-2})
(1-q_v^{-4})\tinsert$ \\
\hline
(rm ur)&$\tinsert\tfrac12(1-q_v^{-1})^2(1-q_v^{-2})\tinsert$ \\
\hline
(rm rm)*&$\tinsert\tfrac12q_v^{-2}(1-q_v^{-2})^2\tinsert$ \\
\hline
(rm rm ur)&$\tinsert\tfrac12q_v^{-2}(1-q_v^{-1})^2(1-q_v^{-2})
\tinsert$ \\
\hline
\end{tabular}
\end{center}
\vspace*{5pt}
\caption{$\vep_v(x)$ for $v$ finite and non-dyadic}%
\label{table-nondyadic}
\end{table}

The first column
displays the index of the orbit and the second, $\vep_v(x)$,
where $x$ is the standard representative for the orbit. The values
of $\vol(K_v\cap G_{x\,k_v}^{\circ})$ which we use here are
contained in Propositions \ref{easyvolume},
\ref{iv-volume-(in_ur)}, \ref{iv-volume-(rm_ur)} and
\ref{iv-volume-(rm_rm_ur)} and the values of $\vol(K_vx)$ in
Propositions \ref{volume1}, \ref{volume2} and \ref{volume4}.

The infinite and dyadic places of $k$ both require special
treatment. We shall begin with the infinite places as the easier of
the two. We extend a classical notation ($r_1$ for the number of
real places and $r_2$ for the number of complex places) by letting
$r_{11}$ be the number of real places of $k$ which split in $\kt$
and $r_{12}$ the number of real places of $k$ which ramify in
$\kt$.
\begin{prop}\label{infvepprod}
For any $S$-tuple $\om_{S}$ we have
\begin{equation*}
\prod_{v\in\gM_{\infty}}\vep_v(\om_v)=
2^{2r_2-r_{11}}
\pi^{3r_{11}+2r_{12}+3r_2}\,.
\end{equation*}
In particular, the product does not depend on $\om_{S}$.
\end{prop}
\begin{proof}
For the standard orbital representatives, $x$, at the infinite
places we have required that $|P(x)|_v=1$ and so
$\vep_v(x)=b_{x,v}^{-1}$. If $v$ is a real place of $k$ which
splits in $\kt$ then $V_{k_v}^{\sst}$ is the union of two orbits
with indices (sp) and (sp rm) respectively. From Propositions
\ref{bxinf-(sp)} and \ref{bxinf-(sp_rm)} we see that
$\vep_v(\om_v)=\tfrac12\pi^3$ for both these orbits. In the
product, the total contribution from these places is thus
$2^{-r_{11}}\pi^{3r_{11}}$. If $v$ is a real place of $k$ which
ramifies in $\kt$ then $V_{k_v}^{\sst}$ is the union of two orbits
with indices (rm) and (rm~rm)* respectively. From Propositions
\ref{bxinf-(rm)} and \ref{bxinf-(rm_rm)*} we see that
$\vep_v(\om_v)=\pi^2$ for both these orbits. In the product, the
total contribution from these places is thus $\pi^{2r_{12}}$.
Finally, if $v$ is a complex place of $k$ then $V_{k_v}^{\sst}$
consists of a single orbit with index (sp) and, from Proposition
\ref{bxinf-(sp)}, $\vep_v(\om_v)=4\pi^3$ for this orbit. The
total contribution to the product from the complex places of $k$ is
thus $2^{2r_2}\pi^{3r_2}$ and the formula follows. \end{proof}

When $v\in\gM_{\text{dy}}$ we shall not calculate the constants
$\vep_{v}(x)$ individually in all cases. Rather we shall sometimes
calculate the sum of the $\vep_{v}(x)$ over a set of orbits with
similar arithmetical properties. This is because if $v\in\gMdy$, 
it is difficult to find the number of ramified quadratic extensions. 
This leads to a final version of
Theorem \ref{prelim-main} which contains no unevaluated constants, but
which employs an equivalence relation, denoted by $\asymp$, coarser
than the relation $\approx$. Our next task is to define this
relation.

Recall that, for $x,y\in V_{k_v}^{\sst}$, we write $x\approx y$ if
and only if $k_v(x)=k_v(y)$ (we have previously used this notation
only when $y$ was a standard orbital representative, but the
extension is convenient here). If $v\notin\gM_{\text{dy}}$ or if
$v\in\gM_{\text{dy}}$ but $k_v(x)/k_v$ is unramified (including the
case $k_v(x)=k_v$) then $x\asymp y$ will have the same meaning as
$x\approx y$. Suppose now that $v\in\gM_{\text{dy}}$ and that
$k_v(x)/k_v$ is ramified. If the type of $x$ is (sp~rm) or (in~rm)
then we shall write $x\asymp y$ if and only if
$\Del_{k_v(x)/k_v}=\Del_{k_v(y)/k_v}$. If the type of $x$ is 
(rm~rm)* or (rm~rm~ur) then we write $x\asymp y$ if and only if $y$
has the same type as $x$. Finally, if $x$ has type (rm~rm~rm) then
we write $x\asymp y$ if and only if
$\Del_{k_v(x)/k_v}=\Del_{k_v(y)/k_v}$ and $\Del_{\kt_v(x)/\kt_v} 
=\Del_{\kt_v(y)/\kt_v}$. This defines an equivalence relation on
$V_{k_v}^{\sst}$ for all places $v$ of $k$. If $\om_S$ is an
$S$-tuple of standard orbital representatives and $x\in L_0$ then
we write $x\asymp\om_S$ to mean that
$x\asymp\om_v$ for all $v\in S$.

The grouping of dyadic orbits is differently expressed in 
\cite{kable-yukie-pbh-II} and we must explain the connection
between the two formulations. 
For any $v\in\gM_{\text{dy}}$ we shall put $2\calo_v=\gp_v^{m_v}$.
If $x\in V_{k_v}^{\sst}$ then let
$\Del_{k_v(x)/k_v}=\gp_v^{\del_{x,v}}$ and, if
$v\notin\gM_{\text{sp}}$, also let
$\Del_{\kt_v/k_v}=\gp_v^{\ti\del_v}$ and
$\Del_{\kt_v(x)/\kt_v}=\ti\gp_v^{\ti\del_{x,v}}$. 
It is well-known that if $k_v(x)/k_v$ is ramified and $v$ is dyadic
then $\del_{x,v}$ takes one of the values $2,4,\dots,2m_v,2m_v+1$.
In \cite{kable-yukie-pbh-II} we introduce a natural number
$\lev(k_1,k_2)$, the \emph{level} of $k_1$ and $k_2$, which is
defined whenever $k_1$ and $k_2$ are ramified quadratic extensions
of a local field. Let us write $\lam_{x,v}=\lev(k_v(x),\kt_v)$ when
$v\in\gM_{\text{rm}}\cap\gM_{\text{dy}}$ and $k_v(x)/k_v$ is
ramified. If $\del_{x,v}\neq\ti\del_{v}$ then
\begin{equation}\label{level-minimum}
\lam_{x,v}=\min\{\lfloor\tfrac12(\del_{x,v}+1)\rfloor,
\lfloor\tfrac12(\ti\del_{v}+1)\rfloor\}\,,
\end{equation}
but if $\del_{x,v}=\ti\del_{v}$ then $\lam_{x,v}$ may take any
value from this minimum up to $\del_{x,v}$. We have the relation
\begin{equation}\label{disc-level-relation}
\ti\del_{x,v}=2(\del_{x,v}-\lam_{x,v})
\end{equation}
and so, with $\Del_{k_v(x)/k_v}$ fixed, $\Del_{\kt_v(x)/\kt_v}$ and
$\lev(k_v(x),\kt_v)$ determine one another. Thus the grouping of
dyadic orbits with index (rm~rm~rm) in \cite{kable-yukie-pbh-II},
by discriminant and level with $\kt_v$, coincides with the grouping
defined here.

If $x$ is a standard orbital representative in $V_{k_v}^{\sst}$ for
any $v\in\gM$ then let us write
\begin{equation*}
\bar{\vep}_v(x)=\sum_{y\asymp x}\vep_v(y)
\end{equation*}
where the sum is over standard orbital representatives which
satisfy $y\asymp x$. Thus $\bar{\vep}_v(x)=\vep_v(x)$ unless
$v\in\gM_{\text{dy}}$ and $k_v(x)/k_v$ is ramified. Also $y\asymp
x$ implies that $x$ and $y$ have the same type and since there is
only one orbit corresponding to each of the indices (rm~rm)* and
(rm~rm~ur), $\bar\vep_v(x)=\vep_v(x)$ if $x$ is the standard
representative for either of these orbits. In Table 
\ref{table-dyadic-ungrouped} we collect
the values of the constants $\vep_v(x)$ for those dyadic orbits
having $\bar\vep_v(x)=\vep_v(x)$ and in Table
\ref{table-dyadic-grouped} we collect the values of the constants
$\bar\vep_v(x)$ for the remaining dyadic orbits. 

\begin{table}
\begin{center}
\begin{tabular}{|c|c|}
\hline
\ Index\ & \tinsert$\vep_v(x)$\tinsert \\
\hline\hline
(sp) &\tinsert $\tfrac12(1+q_v^{-1})(1-q_v^{-2})^2$\tinsert \\
\hline
(in) &\tinsert $\tfrac12(1-q_v^{-1})(1-q_v^{-4})$\tinsert \\
\hline
(rm) &\tinsert $\tfrac12(1-q_v^{-2})^2$\tinsert  \\
\hline
(sp ur) &\tinsert $\tfrac12(1-q_v^{-1})^3(1-q_v^{-2})$\tinsert \\
\hline
(in ur) &\tinsert $\tfrac12(1-q_v^{-1})(1-q_v^{-4})$\tinsert \\
\hline
(rm ur) &\tinsert $\tfrac12(1-q_v^{-1})^2(1-q_v^{-2})$\tinsert \\
\hline
(rm rm)* &\tinsert $\tfrac12q_v^{-2\ti\del_v-2\lfloor\ti
\del_v/2\rfloor}(1-q_v^{-2})^2$
\tinsert \\
\hline
(rm rm ur) &\tinsert
$q_v^{-2\ti\del_v}(1-\tfrac12q_v^{-2\lfloor\ti\del_v/2\rfloor})
(1-q_v^{-1})^2(1-q_v^{-2})$\tinsert \\
\hline
\end{tabular}
\end{center}
\vspace*{5pt}
\caption{$\vep_v(x)$ for ungrouped dyadic orbits}
\label{table-dyadic-ungrouped}
\end{table}
\begin{table}
\begin{center}
\begin{tabular}{|c|c|c|}
\hline
\ Index\ &\ Conditions\ &\tinsert$\bar\vep_v(x)$\tinsert \\
\hline\hline
(sp rm) &$\del_{x,v}\leq 2m_v$& \tinsert 
$q_v^{-\del_{x,v}/2}(1-q_v^{-1})^2(1-q_v^{-2})^3$\tinsert \\
\hline
(sp rm) &$\del_{x,v}=2m_v+1$& \tinsert 
$q_v^{-(m_v+1)}(1-q_v^{-1})(1-q_v^{-2})^3$\tinsert \\
\hline
(in rm) &$\del_{x,v}\leq 2m_v$& \tinsert 
$q_v^{-\del_{x,v}/2}(1-q_v^{-1})^2(1-q_v^{-2})(1-q_v^{-4})$
\tinsert \\
\hline
(in rm) &$\del_{x,v}=2m_v+1$& \tinsert 
$q_v^{-(m_v+1)}(1-q_v^{-1})(1-q_v^{-2})(1-q_v^{-4})$
\tinsert \\
\hline
(rm rm rm) &$\del_{x,v}\neq\ti\del_v$, $\del_{x,v}\leq 2m_v$&\tinsert 
$q_v^{-(\del_{x,v}/2+\lam_{x,v})}
(1-q_v^{-1})^2(1-q_v^{-2})^2$
\tinsert \\
\hline
(rm rm rm) &$\del_{x,v}\neq\ti\del_v$, $\del_{x,v}=2m_v+1$&\tinsert
$q_v^{-(m_v+\lam_{x,v}+1)}(1-q_v^{-1})(1-q_v^{-2})^2$
\tinsert \\
\hline
(rm rm rm) &$\del_{x,v}=\ti\del_v\leq 2m_v$,
$\lam_{x,v}=\tfrac12\ti\del_v$ & \tinsert 
$q_v^{-2\lam_{x,v}}(1-q_v^{-1})(1-2q_v^{-1})(1-q_v^{-2})^2$
\tinsert \\
\hline
(rm rm rm) &$\del_{x,v}=\ti\del_v\leq 2m_v$,
$\lam_{x,v}>\tfrac12\ti\del_v$ & \tinsert 
$q_v^{-2\lam_{x,v}}(1-q_v^{-1})^2(1-q_v^{-2})^2$
\tinsert \\
\hline
(rm rm rm) &$\del_{x,v}=\ti\del_v=2m_v+1$ &\tinsert 
$q_v^{-2\lam_{x,v}}(1-q_v^{-1})^2(1-q_v^{-2})^2$
\tinsert \\
\hline
\end{tabular}
\end{center}
\vspace*{5pt}
\caption{$\bar\vep_v(x)$ for grouped dyadic orbits}
\label{table-dyadic-grouped}
\end{table}

The values of
$\vol(K_vx)$ and $\vol(K_v\cap G_{x\,k_v}^{\circ})$ used to
determine the entries in the two tables were drawn from
Propositions \ref{easyvolume}, \ref{iv-volume-(in_ur)}, 
\ref{iv-volume-(rm_ur)}, \ref{volume1}, \ref{volume3} of this paper
and Propositions 4.3, 5.21, 5.27 and Corollary 5.28
of \cite{kable-yukie-pbh-II}. 
In Table \ref{table-dyadic-grouped}, the second column records the
conditions on $\del_{x,v}$, $\ti\del_v$ and $\lam_{x,v}$ under
which the entry is valid. From (\ref{level-minimum}) and the
observations made in the previous paragraph it is easy to see that
the available conditions are exhaustive.

It will be convenient to extend the notation of section \ref{plan}
by writing
\begin{equation*}
E_v=\sum_x\vep_v(x)
\end{equation*}
for all $v\in\gM_{\text{f}}$, where the sum is taken over all
standard representatives, $x$, of orbits in $G_{k_v}\backslash
V_{k_v}^{\sst}$.  We call $E_v$ the \emph{local density} 
at the place $v$.  
\begin{prop}
Let {\rm $v\in\gM_{\text{f}}$}. Then $E_v=(1-q_v^{-2})E_v'$ where
\begin{equation}\label{Ev'value}
E_v'=\begin{cases} 1-3q_v^{-3}+2q_v^{-4}+q_v^{-5}-q_v^{-6}
&\text{\ if\ } v\in\gM_{\text{{\upshape sp}}}\,, \\
(1+q_v^{-2})(1-q_v^{-2}-q_v^{-3}+q_v^{-4})
&\text{\ if\ } v\in\gM_{\text{{\upshape in}}}\,, \\
(1-q_v^{-1})(1+q_v^{-2}-q_v^{-3}+
q_v^{-2\ti\del_v-2\lfloor\ti\del_v/2\rfloor-1})
&\text{\ if\ } v\in\gM_{\text{{\upshape rm}}}\,.
\end{cases}
\end{equation}
\end{prop}
\begin{proof}
First suppose that $v\notin\gM_{\text{dy}}$. Then every index
corresponds to a single orbit, with the exception of (sp~rm) and
(in~rm), which correspond to two orbits each. Using this and the
values of $\vep_v(x)$ given in Table \ref{table-nondyadic} it is
routine to check the given expressions.

Now suppose that $v\in\gM_{\text{dy}}$. We have
$E_v=\sum_x\bar\vep_v(x)$ where the sum now runs over a complete
set of representatives for the $\asymp$ equivalence classes. The
values of $\bar\vep_v(x)$ are given in Tables
\ref{table-dyadic-ungrouped} and \ref{table-dyadic-grouped} and
using them one can easily establish the claim when
$v\notin\gM_{\text{rm}}$. We carry out the case
$v\in\gM_{\text{rm}}$ explicitly since it is rather more elaborate.

First suppose that $\ti\del_v=2\ti\ell$ with $1\leq\ti\ell\leq
m_v$. The indices which are possible with our assumptions are (rm),
(rm~ur), (rm~rm)* and (rm~rm~ur), corresponding to one orbit each,
and (rm~rm~rm), which corresponds to many orbits. Using Table
\ref{table-dyadic-ungrouped}, the contribution to $E_v$ from the
first four of these indices is
\begin{equation}\label{1stcontrib}
\begin{aligned}
{} & \tfrac12(1-q_v^{-2})^2+\tfrac12(1-q_v^{-1})^2(1-q_v^{-2}) \\
& + \tfrac12 q_v^{-6\ti\ell}(1-q_v^{-2})^2
+q_v^{-4\ti\ell}(1-\tfrac12q_v^{-2\ti\ell})
(1-q_v^{-1})^2(1-q_v^{-2})\,.
\end{aligned}
\end{equation}
Recall that the orbits with index (rm~rm~rm) have been grouped
under $\asymp$ by $\del_{x,v}$ if $\del_{x,v}\neq\ti\del_v$ and by
level if $\del_{x,v}=\ti\del_v$. If $\del_{x,v}\neq\ti\del_v$ then
either $\del_{x,v}=2\ell$ with $\ell\neq\ti\ell$ or
$\del_{x,v}=2m_v+1$. Using Table \ref{table-dyadic-grouped} and 
(\ref{level-minimum}), we see that the contribution from these
equivalence classes is
\begin{equation}\label{2ndcontrib}
\begin{aligned}
\ &\sum_{\ell=1}^{\ti\ell-1}q_v^{-2\ell}(1-q_v^{-1})^2
(1-q_v^{-2})^2 +
\sum_{\ell=\ti\ell+1}^{m_v}q_v^{-(\ell+\ti\ell)}
(1-q_v^{-1})^2(1-q_v^{-2})^2 \\
\ & \quad +q_v^{-(m_v+\ti\ell+1)}(1-q_v^{-1})(1-q_v^{-2})^2 \\
\ & = (q_v^{-2}-q_v^{-2\ti\ell})(1-q_v^{-1})^2(1-q_v^{-2}) \\
& \quad +(q_v^{-(2\ti\ell+1)}-q_v^{-(m_v+\ti\ell+1)})(1-q_v^{-1})
(1-q_v^{-2})^2 \\
\ &\quad +q_v^{-(m_v+\ti\ell+1)}(1-q_v^{-1})(1-q_v^{-2})^2 \\
\ &= q_v^{-2}(1-q_v^{-1})^2(1-q_v^{-2})-
q_v^{-2\ti\ell}(1-q_v^{-1})^2(1-q_v^{-2}) \\
\ &\quad +q_v^{-(2\ti\ell+1)}(1-q_v^{-1})(1-q_v^{-2})^2\,.
\end{aligned}
\end{equation}
If $\del_{x,v}=\ti\del_v$ then the level, $\lam_{x,v}$, runs from
$\ti\ell$ up to $\ti\del_v-1$. By (\ref{disc-level-relation}), the
value $\lam_{x,v}=\ti\del_v=\del_{x,v}$, although possible,
corresponds to the orbit with index (rm~rm~ur) and so is excluded
here. The contribution from the equivalence classes with
$\del_{x,v}=\ti\del_v$ is thus
\begin{equation}\label{3rdcontrib}
\begin{aligned}
\ & q_v^{-2\ti\ell}(1-q_v^{-1})(1-2q_v^{-1})(1-q_v^{-2})^2 \\
& \quad +\sum_{i=\ti\ell+1}^{2\ti\ell-1}
q_v^{-2i}(1-q_v^{-1})^2(1-q_v^{-2})^2 \\
\ &=q_v^{-2\ti\ell}(1-q_v^{-1})(1-2q_v^{-1})(1-q_v^{-2})^2 \\
& \quad +(q_v^{-(2\ti\ell+2)}-q_v^{-4\ti\ell})(1-q_v^{-1})^2(1-q_v^{-2})\,.
\end{aligned}
\end{equation}
Let us now collect all the terms from (\ref{2ndcontrib}) and
(\ref{3rdcontrib}) which have $q_v^{-2\ti\ell}$ as a visible
factor. The result is
\begin{equation*}
\begin{aligned}
\ & q_v^{-2\ti\ell}\big[-(1-q_v^{-1})^2(1-q_v^{-2})+
q_v^{-1}(1-q_v^{-1})(1-q_v^{-2})^2 \\
\ & \hspace*{26pt} +(1-q_v^{-1})(1-2q_v^{-1})(1-q_v^{-2})^2+
q_v^{-2}(1-q_v^{-1})^2(1-q_v^{-2})\big] \\
\ &=q_v^{-2\ti\ell}(1-q_v^{-1})(1-q_v^{-2})\big[
-(1-q_v^{-1})+q_v^{-1}(1-q_v^{-2}) \\
\ &\hspace*{1.87in}+(1-2q_v^{-1})(1-q_v^{-2})
+q_v^{-2}(1-q_v^{-1})\big] \\
\ &=0
\end{aligned}
\end{equation*}
on expanding the factor in the square brackets. It remains to add
(\ref{1stcontrib}), the first term of (\ref{2ndcontrib}) and the
term $-q_v^{-4\ti\ell}(1-q_v^{-1})^2(1-q_v^{-2})$ from
(\ref{3rdcontrib}) to obtain $E_v$. This is easily done. The
situation where $\ti\del_v=2m_v+1$ is similar, but simpler, and we
leave it to the reader.
\end{proof}
In particular, this proposition verifies Condition \ref{evcond}
subject to the results of sections \ref{sec-iv} and \ref{bxsif}.

If $F/k$ is a quadratic extension distinct from $\kt/k$ then
$F=k(x)$ for some $x\in L_0$ and we shall write $F\approx\om_S$ if
$x\approx\om_S$ and $F\asymp\om_S$ if $x\asymp\om_S$.
\begin{thm}  \label{theorem-main}
Let $S\supseteq\gM_{\infty}$ be a finite set of places of $k$ and
$\om_S$ be an $S$-tuple of standard orbital representatives. Then
\begin{equation*}
\lim_{X\to\infty} X^{-2}
\sum_{\substack{[F:k]=2,F\asymp\om_S \\ \esN(\Del_{F/k})\leq X}}
\gC_{F}\gC_{F^*}
\end{equation*}
exists and has the value
\begin{equation*}
2^{-(r_1+r_2+1)}|\Del_{\kt}/\Del_k|^{1/2}\gC_{k}^3\zeta_{\kt}(2)
\prod_{v\in S\setminus\gM_{\infty}}(1-q_v^{-2})^{-1}
\bar\vep_v(\om_v)\cdot
\prod_{v\notin S}E_v'
\end{equation*}
where $\bar\vep_v(x)$ is given by Tables \ref{table-nondyadic},
\ref{table-dyadic-ungrouped} and \ref{table-dyadic-grouped} and
$E_v'$ by (\ref{Ev'value}).
\end{thm}
\begin{proof}
By Proposition \ref{C-factorization} we have
$\gC_{\kt(x)}=\gC_{k}^{-2}\gC_{\kt}\gC_{F}\gC_{F^*}$ if $F=k(x)$.
Recall, from Proposition \ref{xisresidue} and (\ref{zetaresidues}),
that
\begin{align*}
\esR_2&=|\Del_k|^{-1/2}\gC_{k}|\Del_{\kt}|^{-1/2}\gC_{\kt}\cdot
Z_k(2)Z_{\kt}(2)/|\Del_{k}| \\
&=\gC_k\gC_{\kt}|\Del_k|^{-3/2}|\Del_{\kt}|^{-1/2}Z_k(2)Z_{\kt}(2)
\end{align*}
and, from (\ref{Zkdefn}), that
\begin{align*}
Z_k(2)&=2^{-r_2}\pi^{-(r_1+r_2)}|\Del_k|\zeta_k(2)\,, \\
Z_{\kt}(2) &=2^{-\ti r_2}\pi^{-(\ti r_1+\ti r_2)}|\Del_{\kt}|
\zeta_{\kt}(2)\,,
\end{align*}
where $\ti r_1$ is the number of real places of $\kt$ and $\ti r_2$
the number of complex places of this field. Thus
\begin{equation}\label{newR2}
\esR_2=2^{-(r_2+\ti r_2)}\pi^{-(r_1+\ti r_1+r_2+\ti r_2)}
|\Del_{\kt}/\Del_{k}|^{1/2}\gC_k\gC_{\kt}\zeta_k(2)\zeta_{\kt}(2)
\,.
\end{equation}

Let $T=S\cup S_0$ and choose a $T$-tuple, $\om_T'=(\om_v')$. According to
Theorem \ref{prelim-main},
\begin{equation}\label{Climit}
\lim_{X\to\infty}\sum_{
\substack{[F:k]=2,F\approx\om_T' \\ \esN(\Del_{F/k})\leq X}}
\gC_{F}\gC_{F^*}
\end{equation}
exists and equals
\begin{equation*}
\tfrac12\gC_{k}^2\gC_{\kt}^{-1}\esR_2
\prod_{v\in T}\vep_v(\om_v')\cdot
\prod_{v\notin T}E_v\,.
\end{equation*}
Making use of (\ref{newR2}) and Proposition \ref{infvepprod} this
quantity equals
\begin{equation*}
2^{r_2-r_{11}-\ti r_2-1}
\pi^{3r_{11}+2r_{12}+2r_2-r_1-\ti r_1-\ti r_2}
|\Del_{\kt}/\Del_{k}|^{1/2}\gC_{k}^3\zeta_k(2)\zeta_{\kt}(2)
\prod_{v\in T\setminus\gM_{\infty}}\vep_v(\om_v')\cdot
\prod_{v\notin T}E_v\,.
\end{equation*}
But $\ti r_1=2r_{11}$, $\ti r_2=r_{12}+2r_2$ and
$r_1=r_{11}+r_{12}$. Thus
\begin{equation*}
\begin{aligned}
r_2-r_{11}-\ti r_2-1 & =r_2-r_{11}-r_{12}-2r_2-1 \\
& = -(r_2+r_{11}+r_{12}+1)=-(r_1+r_2+1)
\end{aligned}
\end{equation*}
and
\begin{equation*}
\begin{aligned}
3r_{11}+2r_{12}+2r_2-r_1-\ti r_1-\ti r_2 
& = 3r_{11}+2r_{12}+2r_2-r_{11}-r_{12}-2r_{11}-r_{12}-2r_2 \\
& =0
\end{aligned}
\end{equation*}
and we have evaluated (\ref{Climit}) as
\begin{equation}\label{Tfinallim}
2^{-(r_1+r_2+1)}|\Del_{\kt}/\Del_k|^{1/2}\gC_{k}^3\zeta_k(2)
\zeta_{\kt}(2)\prod_{v\in T\setminus\gM_{\infty}}\vep_v(\om_v')\cdot
\prod_{v\notin T}E_v\,.
\end{equation}
Now
\begin{align*}
\prod_{v\notin T}E_v &= \prod_{v\notin T}(1-q_v^{-2})\cdot
\prod_{v\notin T} E_v' \\
&=\zeta_k(2)^{-1}\prod_{v\in
T\setminus\gM_{\infty}}(1-q_v^{-2})^{-1}\cdot
\prod_{v\notin T}E_v'
\end{align*}
and so (\ref{Tfinallim}) equals
\begin{equation}\label{Tsimplim}
2^{-(r_1+r_2+1)}|\Del_{\kt}/\Del_k|^{1/2}\gC_k^3\zeta_{\kt}(2)
\prod_{v\in T\setminus\gM_{\infty}}(1-q_v^{-2})^{-1}\vep_v(\om_v')
\cdot\prod_{v\notin T} E_v'\,.
\end{equation}
Now we sum (\ref{Climit}) and (\ref{Tsimplim}) over all $T$-tuples
$\om_T'=(\om_v')$ which satisfy $\om_v'\asymp\om_v$ for all $v\in
S$ to obtain the statement of the theorem.
\end{proof}
Note that in Theorem \ref{theorem-main}, 
$S$ does not have to contain $S_0$.  

Given an $S$-tuple, $\om_S$, with $S\supseteq\gM_{\infty}$ let us
define
\begin{align*}
n_{++}&=\#\{v\in\gM_{\R}\mid
v\in\gM_{\text{sp}}\text{\ and\ }k_v(\om_v)=k_v\}\,, \\
n_{+-}&=\#\{v\in\gM_{\R}\mid
v\in\gM_{\text{sp}}\text{\ and\ }k_v(\om_v)\neq k_v\}\,, \\
n_{-+}&=\#\{v\in\gM_{\R}\mid
v\in\gM_{\text{rm}}\text{\ and\ }k_v(\om_v)=k_v\}\,, \\
n_{--}&=\#\{v\in\gM_{\R}\mid
v\in\gM_{\text{rm}}\text{\ and\ }k_v(\om_v)\neq k_v\}\,.
\end{align*}
If $F$ is a quadratic extension of $k$ and 
$F\asymp\om_S$ then we denote the composition of 
$F$ and $\ti k$ by $\ti F$ (which 
corresponds to $L$ in Proposition \ref{C-factorization}). 
Then it is easy to see that
\begin{equation*}
\begin{aligned}
r_1(F) & = 2(n_{++}+n_{-+})\,,\quad 
\hskip 5pt r_2(F)=n_{--}+n_{+-}+2r_2\,, \\
r_1(F^*) & = 2(n_{++}+n_{--})\,,\quad r_2(F^*)=n_{+-}+n_{-+}+2r_2\,, \\
r_1(\ti F) & = 4n_{++}\,, 
\hskip 0.85 in r_2(\ti F) = 2(n_{+ -}  +n_{- +} +n_{- -})+4r_2\,, 
\end{aligned}
\end{equation*}
and so $r_1(F)$, $r_1(F^*)$, $r_1(\ti F)$, 
$r_2(F)$, $r_2(F^*)$, and $r_2(\ti F)$ depend only
upon $\om_S$. This allows us to define
\begin{equation*}
\begin{aligned} 
\gc (\om_S) & = 2^{r_1(F)+r_1(F^*)}
(2\pi)^{r_2(F)+r_2(F^*)}, \\
\ti \gc (\om_S) & =2^{r_1(\ti F)} (2\pi)^{r_2(\ti F)}
\end{aligned}
\end{equation*}
where $F\neq\kt$ is any quadratic extension of $k$ satisfying
$F\asymp\om_S$.
\begin{cor}  \label{cor-main1}
Let $S\supseteq\gM_{\infty}$ be a finite set of places of $k$ and
$\om_S$ be an $S$-tuple of standard orbital representatives. Then
\begin{equation*}
\lim_{X\to\infty}X^{-2}
\sum_{\substack{[F:k]=2,F\asymp\om_S \\ \esN(\Del_{F/k})\leq X}}
h_FR_Fh_{F^*}R_{F^*}
\end{equation*}
exists and equals
\begin{equation*}
2^{-(r_1+r_2+1)}\mathfrak{c}(\om_S)^{-1}e_k^2
|\Del_{\kt}/\Del_k|^{1/2}\gC_k^3\zeta_{\kt}(2)
\prod_{v\in
S\setminus\gM_{\infty}}(1-q_v^{-2})^{-1}\bar\vep_v(\om_v)\cdot
\prod_{v\notin S}E_v'\,.
\end{equation*}
\end{cor}
\begin{proof}
Let $F/k$ be a quadratic extension and suppose that $F$ contains a
primitive $n^{\text{th}}$ root of unity, $\zeta_n$, for some
$n$. Since $[\Q(\zeta_n):\Q]=\varphi(n)$, it follows that
$\varphi(n)\leq[F:\Q]=2[k:\Q]$. But it is well-known that
$\varphi(n)\to\infty$ as $n\to\infty$ and so there is some constant
$N$, independent of $F$, such that $n\leq N$. We conclude that, for
all but finitely-many quadratic extensions $F$ of $k$, 
$e_F=e_{F^*}=e_k$.
This finite list of exceptions may be ignored in the limit. Since
\begin{equation*}
\gC_F=2^{r_1(F)}(2\pi)^{r_2(F)}h_FR_Fe_F^{-1}\,,
\end{equation*}
the corollary is now an immediate consequence of the theorem and
the definition of $\mathfrak{c}(\om_S)$.
\end{proof}
\begin{cor}  \label{cor-main2}  With the same assumptions
as in Corollary \ref{cor-main1}, 
\begin{equation*}
\lim_{X\to\infty}X^{-2}
\sum_{\substack{[F:k]=2,F\asymp\om_S \\ \esN(\Del_{F/k})\leq X}}
h_{\ti F}R_{\ti F}
\end{equation*}
exists and equals
\begin{equation*}
\begin{aligned}
{} & 2^{-(r_1+r_2+1)} 2^{r_1(\ti k)}(2\pi)^{r_2(\ti k)}
\ti\gc(\om_S)^{-1}
|\Del_{\kt}/\Del_k|^{1/2}\gC_kh_{\ti k}R_{\ti k}\zeta_{\kt}(2) \\
& \times \prod_{v\in S\setminus\gM_{\infty}}
(1-q_v^{-2})^{-1}\bar\vep_v(\om_v)\cdot
\prod_{v\notin S}E_v'\,.
\end{aligned}
\end{equation*}
\end{cor}
\begin{proof}
By Proposition \ref{C-factorization}, 
$\gC_{\ti F} = \gC_k^{-2}\gC_F\gC_F\gC_{\ti k}$.   
So 
\begin{equation*}
\begin{aligned} 
h_{\ti F}R_{\ti F} & = 2^{-r_1(\ti F)}(2\pi)^{-r_2(\ti F)}
e_{\ti F}\gC_{\ti F} \\
& = \ti\gc(\om_S)^{-1}e_{\ti F}\gC_k^{-2}\gC_F\gC_{F^*}\gC_{\ti k} \\
& = 2^{r_1(\ti k)}(2\pi)^{r_2(\ti k)}
\ti\gc (\om_S)^{-1}\gc (\om_S)
e_{\ti F}e_{\ti k}^{-1}e_F^{-1}e_{F^*}^{-1}
h_{\ti k}R_{\ti k}\gC_k^{-2} h_FR_Fh_{F^*}R_{F^*}
.\end{aligned}
\end{equation*}
As in the proof of Corollary \ref{cor-main1}, $e_F=e_{F^*}=e_k$
and $e_{\ti F} = e_{\ti k}$ except for a finite number of 
quadratic extensions $F$.  
Therefore Corollary \ref{cor-main2} follows from 
Corollary \ref{cor-main1}.
\end{proof}  

We now specialize to the case  $k=\Q$ and $S=\gM_{\infty}$.
Suppose $\ti k = \Q(\sqrt{d_0})$ where $d_0\not=1$ is a square free
integer.  Then $r_1=1,r_2=0,h_k=1,e_k=2$ and $\gC_k=1$.
It is easy to verify that $2^{-(r_1+r_2+1)}\gc(\om_S)^{-1}e_k^2
= \gc(\om_S)^{-1}$ and $2^{-(r_1+r_2+1)} 
2^{r_1(\ti k)}(2\pi)^{r_2(\ti k)}\ti\gc(\om_S)^{-1}$ both coincide
with $c_{\pm}(d_0)^{-1}$ as defined in the introduction.
Therefore Theorems \ref{simple-mainthm1} and \ref{simple-mainthm2}  
are special cases of Corollaries \ref{cor-main1} and \ref{cor-main2}.

\section{The omega sets and their properties}\label{sec-omega}
 
The main purpose of this section is to verify Condition \ref{xicoeff}.
Let $v\in\gMf$ and $x\in V^{\sst}_{k_v}$.  
The function $\Xi_{x,v}(s)$ is defined as an integral over
$G_{k_v}/G_{x\,k_v}^{\circ}$ and our strategy is to replace this by
an integral over a carefully chosen set $\Om_{x,v}\subseteq
G_{k_v}$ called the omega set.  
We impose on the omega set, $\Om_{x,v}$, several
conditions derived from an analysis of Datskovsky's calculations of
standard local zeta functions in \cite{dats}. Once we show that
these conditions can be satisfied, Condition \ref{xicoeff} is an
almost immediate consequence. Thus the bulk of the work in this
section is devoted to finding the omega sets and verifying their
properties.

For the sake of Condition \ref{xicoeff}, 
it is enough to assume that $v\in \gM\setminus S_0$.  
However, verifying Condition \ref{xicoeff} 
will not be our only application of the
existence of omega sets. We shall also require them in certain
proofs in section \ref{bxsif}  and, for this, greater generality
will be needed. Thus we shall allow $v$ to be any finite place of
$k$ and consider orbits of types other than 
three types (rm rm)*, (rm rm ur), and (rm rm rm) at dyadic 
places $v\in\gMdy$.  

Before we begin, we shall record as a lemma a simple observation
which will be useful both later in this section and in the next.
\begin{lem}\label{differentreps}
Suppose that $v\in\gM$, $x\in V_{k_v}^{\text{\upshape ss}}$ and $y\in
G_{k_v}x$. If $|P(x)|_v=|P(y)|_v$ then
$Z_{x,v}(\Phi,s)=Z_{y,v}(\Phi,s)$ for all $\Phi\in\cS(V_{k_v})$.
\end{lem}
\begin{proof}
Examining the second equation in Definition \ref{lzeta} we see, in
light of Proposition \ref{bxindep} and the hypotheses, that every
factor in the definition of the local orbital zeta function remains
unchanged when we replace $x$ by $y$.
\end{proof}
For each $x\in V_{k_v}^{\sst}$ we choose an element $g_x\in
G_{k_v(x)}$ such that $g_xw=x$ and $g_x$ satisfies Condition 
\ref{ghomcond} if $k_v(x)\neq k_v$. From this choice we obtain an
isomorphism $\theta_{g_x}:G_{x\,k_v}^{\circ}\to H_{x\,k_v}$ where
$H_{x\,k_v}$ is defined by (\ref{Hxdefn}).
\begin{defn}\label{omegadefn} A set $\Om_{x,v}\subseteq G_{k_v}$ is
called an omega set for $x$ if it has the following properties:
\begin{itemize}
\item[(1)] $\Om_{x,v}x=(G_{k_v}x)\cap V_{\calo_v}$.
\item[(2)]
$K_v\Om_{x,v}\theta_{g_x}^{-1}(H_{x\,\calo_v})=\Om_{x,v}$.
\item[(3)] If $g_1,g_2\in\Om_{x,v}$, $h\in G_{x\,k_v}^{\circ}$ and
$g_1=g_2h$ then $h\in\theta_{g_x}^{-1}(H_{x\,\calo_v})$.
\item[(4)] If $g\in\Om_{x,v}$ then $|\chi(g)|_v\leq1$ with equality
only if $g\in K_v$.
\end{itemize}
\end{defn}

Below we give omega sets for representatives of each of the orbit
types which we require. These include the six orbit types possible
under the restriction that $v\notin S_0$, as well as the orbits of
type (rm) and (rm~ur). For
the orbits of type (sp), (in) and (rm)
it will be convenient to use $x=w$
as the orbital representative instead of the standard $w_p$. This
is permissible for the purpose at hand 
by Lemma \ref{differentreps}. For the orbits of
types (sp~ur), (sp~rm), (in~ur), (in~rm) and (rm~ur) we shall use the
standard representatives.

If $p(z)=z^2+a_1z+a_2\in k_v[z]$ then we shall let
$\al=\{\al_1,\al_2\}$ be the set of roots of $p$ and write
$e(\al)={}^{t}(1\ -\al_1)$ (a column vector in $k_v(w_p)^2$). 
If $\ell={}^{t}(\ell_1\ \ell_2)$ is any such column vector then
we set
$\Vert\ell\Vert=\max\{|\ell_1|_{k_v(w_p)},|\ell_2|_{k_v(w_p)}\}$.
Let $t$ be as in (\ref{standardtoruselement}) for the field $k_v$
and $n(u)=\big(n(u_1),n(u_2),n(u_3)\big)$ for $u=(u_1,u_2,u_3)\in
k_v^3$ or $n(u)=\big(n(u_1),n(u_2)\big)$ for
$u=(u_1,u_2)\in\kt_v\times k_v$. Let $g=\k tn(u)$ be the Iwasawa
decomposition of $g\in G_{k_v}$. In section \ref{final} we
described the form of the polynomial $p(z)$ for each of the
standard orbital representatives. It will be convenient here to add
the assumption that $a_1=0$ whenever $v$ is not dyadic, as we may.

For the index (sp) with orbital representative $x=w$ we define
\begin{equation}\label{omega1}
\Om_{x,v}=\{g=\k tn(u)\mid t_{ij}=1\text{\ for\ }i,j=1,2
\text{\ and\ }gx\in V_{\calo_v}\}\,.
\end{equation}
For the indices (in) and (rm)
with orbital representative $x=w$ we define
\begin{equation}\label{omega2}
\Om_{x,v}=\{g=\k tn(u)\mid t_{11}=t_{12}=1\text{\ and\ }gx\in
V_{\calo_v}\}\,.
\end{equation}
For the index (sp~ur) with orbital representative $x=w_p$
we define
\begin{equation}\label{omega3}
\Om_{x,v}=\{g=(g_1,g_2,g_3)\big| 
|\det(g_1)|_v=1\text{ or }q_v^{-1}, |\det(g_2)|=1\text{\ or\ }
q_v,\ gx\in V_{\calo_v}\}\,.
\end{equation}
For the index (sp~rm) with orbital representative $x=w_p$
we define
\begin{equation}\label{omega4}
\Om_{x,v}=\{g=(g_1,g_2,g_3)\big|
|\det(g_i)|_v=1\text{\ for\ }i=1,2,\ gx\in V_{\calo_v}\}\,.
\end{equation}
For the index (in~ur) with orbital representative $x=w_p$
we define
\begin{equation}\label{omega5}
\Om_{x,v}=\{g=(g_1,g_2)\big|
|\det(g_1)|_{\kt_v}=1,\Vert g_1e(\al)\Vert=1,\ gx\in
V_{\calo_v}\}\,.
\end{equation}
For the index (in~rm) with orbital representative $x=w_p$
we define
\begin{equation}\label{omega6}
\Om_{x,v}=\{g=(g_1,g_2)\big|
|\det(g_1)|_{\kt_v}=1,\ gx\in V_{\calo_v}\}\,.
\end{equation}
Finally, for the index (rm~ur) with orbital representative $x=w_p$
we define
\begin{equation}\label{omega7}
\Om_{x,v}=\{g=(g_1,g_2)\big|
|\det(g_1)|_{\kt_v}=1\text{\ or\ }q_v^{-1},\ gx\in V_{\calo_v}\}\,.
\end{equation}
In every case we shall write
\begin{equation}\label{omega8}
\Om_{x,v}^1=\{g\in\Omega_{x,v}\big||\chi(g)|_v=1\}\,.
\end{equation}
\begin{prop} The sets defined by 
(\ref{omega1})--(\ref{omega7})
have properties (1), (2) and (3) 
of Definition \ref{omegadefn}.
\end{prop}
\begin{proof}
If $\k\in K_v$ then $\k V_{\calo_v}=V_{\calo_v}$, $|\det(\k)|_v=1$
and $\Vert\k e\Vert=\Vert e\Vert$ for any vector $e$. This makes it
clear that $K_v\Om_{x,v}=\Om_{x,v}$ in all cases. The rest of the
argument will be case by case, but we make two observations which
will be used repeatedly. First, it follows at once from the
definition in every case that $\Om_{x,v}x\subseteq G_{k_v}x\cap
V_{\calo_v}$ and so to establish (1) we need only prove the reverse
inclusion. This will be done if we can show that given $g\in
G_{k_v}$ with $gx\in V_{\calo_v}$
we can find $h\in G_{x\,k_v}^{\circ}$ such that
$gh\in\Om_{x,v}$. Secondly, any $h\in G_{x\,k_v}^{\circ}$ may be
expressed as $h=g_xs_x(t_x)g_x^{-1}$, in the notation of 
(\ref{stabelm1})--(\ref{stabelm4}), and
$h\in\theta_{g_x}^{-1}(H_{x\,\calo_v})$ if and only if all the
components of $t_x$ are units.

Consider the cases (sp), (in) and (rm). We may assume, for simplicity,
that $g_x$ has been chosen to be the identity. Take $g\in G_{k_v}$
with $gx\in V_{\calo_v}$
and let $g=\k(g)t(g)n(u(g))$ be its Iwasawa decomposition. Let
$s_x(t_x)$ be as in (\ref{stabelm1}) or (\ref{stabelm5}). By
choosing
$t_x=(t_{11}(g)^{-1},t_{12}(g)^{-1},t_{21}(g)^{-1},t_{22}(g)^{-1})$ in
the first case and $t_x=(t_{11}(g)^{-1},t_{12}(g)^{-1})$ in the
second, we may arrange that $gs_x(t_x)\in\Om_{x,v}$. This proves
property (1). Moreover, if $g\in\Om_{x,v}$ and all the components
of $t_x$ are units then commuting $s_x(t_x)$ past the $T_{k_v}$ and
$N_{k_v}$ factors in the Iwasawa decomposition and absorbing it
into the $K_v$ factor shows that $gs_x(t_x)\in\Om_{x,v}$ also,
which proves property (2). For property (3), observe that in the
Iwasawa decomposition, the $T_{k_v}$ factor is unique up to
multiplication of its diagonal elements by units. Thus if
$g_1,g_2\in\Om_{x,v}$ and $g_1=g_2h$ with $h=s_x(t_x)$ then
$s_x(t_x)\in H_{x\,\calo_v}$. This proves property (3).

We next turn to case (sp~ur). Let $s_x(t_x)$ be as in
(\ref{stabelm2}) and $g\in G_{k_v}$ with $gx\in V_{\calo_v}$.
Note that
\begin{equation*}
|\det s_{x1}(t_x)|_v=|\n_{k_v(x)/k_v}(t_{11})|_v
\end{equation*}
and since $k_v(x)/k_v$ is unramified, this may be any even power of
$q_v$. The same holds for $|\det s_{x2}(t_x)|_v$ and the
determinants of the components of $g_xs_x(t_x)g_x^{-1}$ are the
same as those of the corresponding components of $s_x(t_x)$. It
follows that we can arrange
$g\big(g_xs_x(t_x)g_x^{-1}\big)\in\Om_{x,v}$ for a suitable choice
of $t_x$ and this proves (1). If $s_x(t_x)\in H_{x\,\calo_v}$ then
the determinants of each of its components are units and this makes
(2) obvious. Also, this argument shows that if
$g_1,g_2\in\Om_{x,v}$, $h=g_xs_x(t_x)g_x^{-1}$ and $g_1=g_2h$ then
$t_{11}$ and $t_{21}$ are units, which implies that $s_x(t_x)\in
H_{x\,\calo_v}$; hence (3).

The case (sp~rm) is very similar, with the one difference that
since $k_v(x)/k_v$ is ramified, $|\det s_{xj}(t_x)|_v$ can be any
integer power of $q_v$.

Next we treat (in~ur). Let $g=(g_1,g_2)\in G_{k_v}$ with $gx\in
V_{\calo_v}$ and $s_x(t_x)$ be as in (\ref{stabelm3}). Note that
$e(\al)$ is an eigenvector for the first component of
$h=g_xs_x(t_x)g_x^{-1}$ with eigenvalue $t_{11}$. So if
$h=(h_1,h_2)$ then $\Vert g_1h_1e(\al)\Vert=
|t_{11}|_{\kt_v}\Vert g_1e(\al)\Vert$. Also,
$|\det(g_1h_1)|_{\kt_v}=|\det(g_1)|_{\kt_v}|t_{11}t_{12}|_{\kt_v}$.
We are free to choose the pair $(t_{11},t_{11}t_{12})\in\kt_v^2$
arbitrarily and so there exists $h\in G_{x\,k_v}^{\circ}$ with
$gh\in\Om_{x,v}$, proving (1). If $g\in\Om_{x,v}$ and
$h\in\theta_{g_x}^{-1}(H_{x\,\calo_v})$ then $t_{11}$ and $t_{12}$
are units and so $\Vert ghe(\al)\Vert=\Vert ge(\al)\Vert$ and
$|\det(gh)|_{\kt_v}=|\det(g)|_{\kt_v}$, which proves (2). Also, if
$g_1,g_2\in\Om_{x,v}$, $h=g_xs_x(t_x)g_x^{-1}$ and $g_1=g_2h$ then
$|t_{11}|_{\kt_v}=|t_{11}t_{12}|_{\kt_v}=1$, which implies that
$h\in\theta_{g_x}^{-1}(H_{x\,\calo_v})$ and (3) follows.

Finally, cases (in~rm) and (rm~ur) are very similar to cases
(sp~rm) and (sp~ur). Note that if $s_x(t_x)$ is as in 
(\ref{stabelm4}) then $|\det s_{x1}(t_x)|_{\ti k_v}=|\n_{\ti
k_v(x)/\ti k_v}(t_{11})|_{\ti k_v}$. In case (in~rm), $\ti
k_v(x)/\ti k_v$ is ramified and so this takes every value in $|\ti
k_v^{\times}|_{\ti k_v}$. In case (rm~ur), $\ti k_v(x)/\ti k_v$ is
unramified and so $|\det s_{x1}(t_x)|_{\ti k_v}$ takes every value
in $|(\ti k_v^{\times})^2|_{\ti k_v}$. The rest of the argument is
identical to that in the cases already mentioned.
\end{proof}

Using only parts (1), (2) and (3) of Definition \ref{omegadefn} we
can prove the following.
\begin{prop}\label{omegaint}
Let $\Psi_{x,v}$ be the characteristic function of 
$\Om_{x,v}$.  Then 
\begin{equation*}
Z_{x,v}(\Phi_{v,0},s) = \int_{G_{k_v}}|\chi(g)|_v^s 
\Psi_{x,v}(g)\,dg_v\,.\end{equation*}
\end{prop} 
\begin{proof}
Since 
\begin{equation*}
dg_v = d\ti g_v \md \ti t_v,\; dg_{x,v}'' = d\ti g_{x,v}'' 
\md \ti t_v,\; dg_v = b_{x,v}dg_{x,v}'dg_{x,v}''
,\end{equation*}
$d\ti g_v = b_{x,v} dg_{x,v}'d\ti g_{x,v}''$.  
So the right hand side of the above identity is 
\begin{equation} \label{f1}
b_{x,v}\int_{G_{k_v}/G^{\circ}_{x\, k_v}}|\chi(g_{x,v}')|_v^s 
\left( \int_{G^{\circ}_{k_v}} \Psi_{x,v}(g_{x,v}'g_{x,v}'')\,
dg_{x,v}'' \right)\,dg_{x,v}'\,.\end{equation} 

By (2) and (3) of Definition \ref{omegadefn},
$\Psi_{x,v}(g_{x,v}'g_{x,v}'')$ is non-zero 
if and only if $g_{x,v}'\in \Om_{x,v}$ and 
$g_{x,v}''\in \theta_{g_x}^{-1}(H_{x\, \co_v})$.
Since we chose the measure
$dg_{x,v}''$ so that the volume of this set is one, 
\begin{equation*}
\int_{G^{\circ}_{k_v}} \Psi_{x,v}(g_{x,v}'g_{x,v}'') dg_{x,v}'' 
\end{equation*}
is the  characteristic function  of 
$\Om_{x,v}G^{\circ}_{x\,k_v}/G^{\circ}_{x\,k_v} \cong G_{k_v} x\cap V_{\co_v}$. 
Therefore, (\ref{f1}) is 
\begin{equation*}
b_{x,v}\int_{G_{k_v}/G^{\circ}_{x\, k_v}}|\chi(g_{x,v}')|_v^s 
\Phi_{v,0}(g_{x,v}'x)dg_{x,v}' 
,\end{equation*}
which is the definition of $Z_{x,v}(\Phi_{v,0},s)$.  
\end{proof} 

Before we verify part (4) of Definition \ref{omegadefn} it will be
convenient to prove three lemmas. First note that we may let
$\gl(2)_{k_v}$ act on the space of quadratic polynomials in
$k_v[z]$ by regarding such polynomials as the inhomogeneous forms
of binary quadratic forms. With this convention, if
$p(z)=z^2+a_1z+a_2\in k_v[z]$ and $g=a(t_1,t_2)n(u)$ then
\begin{equation*}
gp(z)=t_1^2z^2+t_1t_2(2u+a_1)z+t_2^2(u^2+a_1u+a_2)\,.
\end{equation*}
\begin{lem}\label{ramlem}
Suppose that $p(z)$ is an Eisenstein polynomial. Let
$t\in k_v^{\times}$, $u\in k_v$, $i=0$ or $1$ and suppose that
$\pi_v^ia(t,t^{-1}\pi_v^{-i})n(u)p(z)\in\calo_v[z]$. Then
$t\in\calo_v^{\times}$ and $u\in\calo_v$. Moreover, if $i=1$ then
$u\in\gp_v$.
\end{lem}
\begin{proof}
We have $\pi_v^i t^2\in\calo_v$, which implies that $t\in\calo_v$
since $i=0$ or $1$. Since
$t^{-2}\pi_v^{-i}(u^2+a_1u+a_2)\in\calo_v$, $(u^2+a_1u+a_2)\in
t^2\pi_v^i\calo_v$. In particular, $u^2+a_1u+a_2\in\calo_v$ and so
$u(u+a_1)\in\calo_v$. If $u\notin\calo_v$ then
$\ord(u+a_1)=\ord(u)$ and we reach a contradiction. Hence
$u\in\calo_v$. The order of $u^2+a_1u+a_2$ is either $0$ (if
$u\in\calo_v^{\times}$) or $1$ (if $u\in\gp_v$). If $i=0$ this
forces $t\in\calo_v^{\times}$ and if $i=1$ it forces first
$t\in\calo_v^{\times}$ and then $u\in\gp_v$.
\end{proof}
\begin{lem}\label{urlem} Suppose that $p(z)=z^2+a_2$ with
$-a_2\in\calo_v^{\times}\setminus (\calo_v^{\times})^2$, if
{\rm $v\notin\gM_{\text{dy}}$}, or that $p(z)$ is an Artin-Schreier
polynomial, if {\rm $v\in\gM_{\text{dy}}$}.
Let $t\in
k_v^{\times}$, $u\in k_v$, $i=-1$, $0$ or $1$ and suppose that
$\pi_v^ia(t,t^{-1}\pi_v^{-i})n(u)p(z)\in\calo_v[z]$. Then $i=0$,
$t\in\calo_v^{\times}$ and $u\in\calo_v$.
\end{lem}
\begin{proof}
The conditions imply that $\pi_v^i t^2$ and $t^{-2}\pi_v^{-i}p(u)$
are integral. Since $-1\leq i\leq 1$, $t\in\calo_v$. Thus
$p(u)\in\pi_v^{-1}\calo_v$, which implies that
$u(u+a_1)\in\pi_v^{-1}\calo_v$. If $u\notin\calo_v$ then
$\ord(u)=\ord(u+a_1)$ and so $\ord(u(u+a_1))$ is a negative, even
integer. This is a contradiction and so $u\in\calo_v$. The
reduction of the polynomial $p(z)$ has no roots in
$\calo_v/\gp_v$ and thus $p(u)\in\calo_v^{\times}$ for all
$u\in\calo_v$. It follows that $t^2\pi_v^i\in\calo_v^{\times}$.
This gives $i=0$ and $t\in\calo_v^{\times}$, as required.
\end{proof}
\begin{lem}\label{lessthanlem}
Let $x$ be a standard orbital representative and
suppose that $y\in V_{\calo_v}$ lies in the orbit of $x$ under
$G_{k_v}$. Then $|P(y)|_v\leq |P(x)|_v$.
\end{lem}
\begin{proof}
If $k_v(x)=k_v$ then $|P(x)|_v=1$ and $P(y)\in\calo_v$ since $y\in
V_{\calo_v}$. The statement follows in this case. We now assume
that $k_v(x)\neq k_v$. Let
$F_y(v_1,v_2)=b_0v_1^2+b_1v_1v_2+b_2v_2^2$ and consider the
polynomial $r(z)=z^2+b_1z+b_0b_2$. Since $y\in V_{\calo_v}$,
$b_0,b_1,b_2\in\calo_v$ and so $r(z)\in\calo_v[z]$. The
discriminant of $r(z)$ is equal to the discriminant of $F_y$ and so
if $\beta$ is a root of $r(z)$ then $\beta\in k_v(y)=k_v(x)$. It
follows that $\calo_v[\beta]\subseteq\calo_{k_v(x)}$ and hence that
$P(y)\calo_v\subseteq\Del_{k_v(x)/k_v}$. But the standard orbital
representative was chosen so that $\Del_{k_v(x)/k_v}=P(x)\calo_v$
and the statement follows in this case also.
\end{proof}

\begin{prop}\label{lzct1}
The sets defined by (\ref{omega1})--(\ref{omega7})
have property (4) of Definition \ref{omegadefn}. 
Consequently, they are omega sets.
\end{prop} 
\begin{proof} 
If $g\in\Om_{x,v}$ then $gx\in V_{\calo_v}$ and so
$|P(gx)|_v\leq|P(x)|_v$ by Lemma \ref{lessthanlem}. But
$|P(gx)|_v=|\chi(g)|_v|P(x)|_v$ and it follows that
$|\chi(g)|_v\leq 1$. This establishes the first part of (4) in
Definition \ref{omegadefn}.

We now have to show that if $g\in\Om_{x,v}^{1}$ then $g\in K_v$.
The orbital representatives have already been fixed in
(\ref{omega1})--(\ref{omega7}) and the notation introduced there
will be used without comment below.

We begin with the cases (sp), (in) and (rm). Let $g\in\Om_{x,v}^{1}$; we
have to show that $g\in K_v$. By (2) of Definition \ref{omegadefn},
$\Om_{x,v}^{1}$ is left $K_v$-invariant and so we may assume that
$g=tn(u)$. Since $g\in\Om_{x,v}$ we have
$t_{11}=t_{12}=t_{21}=t_{22}=1$ in case (sp) and $t_{11}=t_{12}=1$
in cases (in) and (rm). The assumption that $|\chi(g)|_v=1$ implies that
$|t_{31}t_{32}|_v=1$ in case (sp) and that $|t_{21}t_{22}|_v=1$ in
cases (in) and (rm). In case (sp) we have
\begin{equation}
gw = \left(t_{31} \pmatrix 1 & u_2\\ u_1 & u_1u_2\endpmatrix,
t_{32}\pmatrix u_3 & u_2u_3\\ u_1u_3 & 1+u_1u_2u_3\endpmatrix\right)
.\end{equation} 
and in cases (in) and (rm) we have
\begin{equation}
gw = \left(t_{21} \pmatrix 1 & u_1^{\sig}\\ 
u_1 & \n_{\ti k_v/k_v}(u_1)\endpmatrix,
t_{22}\pmatrix u_3 & u_1^{\sig}u_3\\ 
u_1u_3 & 1+\n_{\ti k_v/k_v}(u_1)u_3\endpmatrix\right)
.\end{equation} 
Let $a= \text{ord}_{k_v} (t_{31})$ or $\text{ord}_{k_v} (t_{21})$. 
Then, by assumption, 
$\text{ord}_{k_v} (t_{32})=-a$ or 
$\text{ord}_{k_v} (t_{22})=-a$. 
Consider case (sp).  
Let $\bar u_i=\pi_v^a u_i$ for $i=1,2$, and $\bar u_3=\pi_v^{-a}u_3$.
Then $gw\in V_{\co_v}$ if and only if 
\begin{equation*}
\begin{aligned}
{}& \pi_v^a,\; \bar u_1,\bar u_2,\bar u_3,\;  
\pi^{-a}\bar u_1\bar u_2,\; \pi^{-a}\bar u_1\bar u_3,\; 
\pi^{-a}\bar u_2\bar u_3, \\
& \pi_v^{-a}(1+ \pi_v^{-a}\bar u_1\bar u_2\bar u_3)\,.
\end{aligned}
\end{equation*}
are integral. So $a\geq 0$.  
We assume $a>0$ and  deduce a contradiction. 
Suppose $\bar u_1$ is not a unit.  Then 
\begin{equation*}
\pi_v^{-a}\bar u_1\bar u_2\bar u_3
= (\pi_v^{-a}\bar u_2\bar u_3)\bar u_1
\equiv 0 \;(\gp_v)\,.
\end{equation*}
Then  $1+ \pi_v^{-a}\bar u_1\bar u_2\bar u_3$ is a 
unit. This implies $\pi_v^{-a}(1+ \pi_v^{-a}\bar u_1\bar u_2\bar u_3)
\notin \co_v$, which is a contradiction.  
So $\bar u_1$ is a unit and similarly $\bar u_2,\bar u_3$
are units also.    
Then the order of $\pi_v^{-a}(1+ \pi_v^{-a}\bar u_1\bar u_2\bar u_3)$
is $-2a$, which is a contradiction.  This implies $a=0$. 
Then $u_i\in \co_v$ for $i=1,2,3$.  
Cases (in) and (rm) are similar using $u_1,u_1^{\sig},u_2$
in the places of $u_1,u_2,u_3$ above.  The only difference is 
that we consider elements in $\ti \co_v$. 

Next we treat the case (sp~rm).
Suppose $g=(g_1,g_2,g_3)\in \Om^1_{x,v}$. 
Then $|\det g_i|_v=1$ for $i=1,2,3$.  
We may assume that $g_1,g_2,g_3$ are lower triangular.  
Note that $F_{w_p}(z,1)=p(z)$.  
So $F_{gw_p}(z,1)= (\det g_1\det g_2)g_3p(z)$ is integral. 
Since $\det g_1,\det g_2\in \co_v^{\times}$, 
$g_3\in \gl(2)_{\co_v}$ by Lemma \ref{ramlem}.  

In this case, we can regard $V$ as 
$\aff^2\otimes \aff^2\otimes \aff^2$.  
Instead of the third factor, we can use the first and the second factors
to make equivariant maps similar to $F_x$.
Then because of the symmetry of
our element $w_p$, $g_1,g_2\in \gl(2)_{\co_v}$ by Lemma \ref{ramlem}
again. This concludes the verification in this case.

Now we consider the case (sp ur).
Let $g=(g_1,g_2,g_3)\in \Om_{x,v}$ and $|\chi(g)|_v=1$.  
In this case there are four possibilities
as follows:
\begin{itemize}
\item[(A)] $|\det g_1|_v=|\det g_2|_v=1$,
\item[(B)] $|\det g_1|_v=1,\; |\det g_2|_v=q_v^{-1}$,  
\item[(C)] $|\det g_1|_v=q_v^{-1},\; |\det g_2|_v=1$,  
\item[(D)] $|\det g_1|_v=q_v^{-1},\; |\det g_2|_v=q_v$.  
\end{itemize}
For these cases, $|\det g_3|_v=1,q_v,q_v,1$ respectively.  
The argument in case (A) is similar to that used in
case (sp rm).
In case (B), $F_{gw_p}(z,1) = \pi_v g_3 p(z)$, and 
$\det g_3=\pi_v^{-1}$.  Since $F_{gw_p}(z,1)$ is integral, 
this corresponds to the case $i=1$ in Lemma \ref{urlem}. 
Therefore this cannot happen. Cases (C), (D) are similar to 
case (B) because of the symmetry (considering an equivariant 
map using the second $\aff^2$ factor  in case (D)).    

Now we consider the case (in~ur). Suppose that
$g=(g_1,g_2)\in\Om_{x,v}^1$. This implies that
$|\det(g_1)|_{\kt_v}=|\det(g_2)|_v=1$. We have
\begin{equation*}
F_{gx}(z,1)= \n_{\ti k_v/k_v}(\det g_1)g_2p(z)
\end{equation*}
and, since $\n_{\kt_v/k_v}(\det g_1)$ is a unit by assumption,
$g_2\in\gl(2)_{\calo_v}$ by Lemma \ref{urlem}. Since $\Om_{x,v}$ is
left $K_v$-invariant we may assume that $g_2=1$ and that $g_1$ is
lower triangular, say $g_1=a(t_{11},t_{12})n(u_1)$. Note that
\begin{equation}\label{eigen}
g_1e(\al) = \pmatrix t_{11}\\ t_{12}(u_1-\al_1)\endpmatrix\,.
\end{equation}
and this is a primitive integral vector. Computation gives
$(g_1,1)w_p=(M_1,M_2)$ where
\begin{equation}\label{expression}
\begin{aligned}
M_1 & = \pmatrix 0 & t_{11}t_{12}^{\sig}\\ t_{11}^{\sig}t_{12} & 
\n_{\ti k_v/k_v}(t_{12})[\tr_{\ti k_v/k_v}(u_1)+a_1]\endpmatrix\,, \\
M_2 & = \pmatrix
\n_{\ti k_v/k_v}(t_{11}) & t_{11}t_{12}^{\sig}(u_1^{\sig}+a_1)\\ 
t_{11}^{\sig}t_{12}(u_1+a_1)&  
\n_{\ti k_v/k_v}(t_{12})m(u_1,p)\endpmatrix
\end{aligned}
\end{equation}
with
\begin{equation*}
m(u_1,p)=a_1^2-a_2+a_1\tr_{\ti k_v/k_v}(u_1)+\n_{\ti k_v/k_v}(u_1)
\end{equation*}
and both these matrices must be integral.
Let $\bar u_1=u_1-\al_1$. Then $\tr_{\ti k_v/k_v}(u_1)+a_1=
\tr_{\ti k_v/k_v}(\bar u_1)$ and
\begin{equation*}
m(u_1,p)=\n_{\ti k_v/k_v}(\bar u_1)-\tr_{\ti k_v/k_v}(\al_1\bar
u_1)
\end{equation*}
and so $M_1$ and $M_2$ are integral if and only if
\begin{equation*}
\begin{aligned}
\  & t_{11},\ \bar u_1,\ \n_{\kt_v/k_v}(t_{12})\tr_{\kt_v/k_v}(\bar
u_1),  \\
\  & \n_{\kt_v/k_v}(t_{12})
[\n_{\kt_v/k_v}(\bar u_1)-\tr_{\kt_v/k_v}(\al_1\bar u_1)]
\end{aligned}
\end{equation*}
are integral. Since $\al_1\in\ti\calo_v$, it follows that
$u_1\in\ti\calo_v$. Also $t_{11}\in\ti\calo_v$ and it remains to
show that $t_{11}$ and $t_{12}$ are units. From the definition of
$\Om_{x,v}$ we know that $|t_{11}t_{12}|_{\kt_v}=1$. Let
$\ord_{\kt_v}(t_{11})=i$; we assume that $i>0$ and deduce a
contradiction. We have $\ord_{\kt_v}(t_{12})=-i$ and, from
(\ref{eigen}), we conclude that $\ord_{\kt_v}(\bar u_1)=i$. Thus we
may write $\bar u_1=\pi_v^i(\bar u_{11}+\bar u_{12}\al_1)$ where
$\bar u_{11},\bar u_{12}\in\calo_v$ and $\bar u_{11}+\bar
u_{12}\al_1\in\ti\calo_v^{\times}$. Then
\begin{align*}
\n_{\kt_v/k_v}(\bar u_1) &= \pi_v^{2i}[\bar u_{11}^2-a_1\bar
u_{11}\bar u_{12}+a_2\bar u_{12}^2], \\
\tr_{\kt_v/k_v}(\bar u_1) &= \pi_v^i[2\bar u_{11}-a_1\bar u_{12}],
\\
\tr_{\kt_v/k_v}(\al_1\bar u_1) &= \pi_v^i[-a_1\bar
u_{11}+(a_1^2-2a_2)\bar u_{12}]
\end{align*}
and, since $\ord_{k_v}(\n_{\kt_v/k_v}(t_{12}))=-2i$, it follows
that
\begin{equation}\label{linsystem1}
\begin{aligned}
-a_1\bar u_{11}+(a_1^2-2a_2)\bar u_{12} &\equiv 0\;(\gp_v^i),
\\ 2\bar u_{11}-a_1\bar u_{12} &\equiv 0\; (\gp_v^i)\,.
\end{aligned}
\end{equation}
Regarding this as a linear system for $(\bar u_{11},\bar u_{12})$,
the determinant of the coefficient matrix is $-a_1^2+4a_2=-P(x)$.
This is a unit by the choice of $x$ and so (\ref{linsystem1})
implies that $(\bar u_{11},\bar u_{12})\equiv (0,0)\;(\gp_v)$.
This contradicts $\bar u_{11}+\bar
u_{12}\al_1\in\ti\calo_v^{\times}$ and so $i=0$. This completes the
case (in~ur).

Next we must deal with the case (in~rm). Suppose that
$g=(g_1,g_2)\in\Om_{x,v}^1$. By similar arguments to those of the
previous case, using Lemma \ref{ramlem} in place of Lemma
\ref{urlem}, we see that $g_2\in\gl(2)_{\calo_v}$. Hence we may
assume that $g_2=1$ and that $g_1=a(t_{11},t_{12})n(u_1)$ is lower
triangular. Then $(g_1,1)w_p=(M_1,M_2)$ where $M_1$ and $M_2$ are
given by (\ref{expression}). Since
$t_{11}t_{12}^{\sig}\in\ti\calo_v^{\times}$, $M_1$ and $M_2$ are
integral if and only if
\begin{equation}\label{intcond-inrm}
t_{11},\ u_1,\ \n_{\kt_v/k_v}(t_{12})[\tr_{\kt_v/k_v}(u_1)+a_1],\ 
\n_{\kt_v/k_v}(t_{12})m(u_1,p)
\end{equation}
are integral. Let $\ord_{\kt_v}(t_{11})=i$; we shall again assume
that $i>0$ and derive a contradiction. We have
$\ord_{\kt_v}(t_{12})=-i$, so that
$\ord_{k_v}(\n_{\kt_v/k_v}(t_{12}))=-2i$. Thus
$\tr_{\kt_v/k_v}(u_1)\equiv -a_1\;(\gp_v^{2i})$ and, since $p(z)$
is an Eisenstein polynomial, it follows that
$\tr_{\kt_v/k_v}(u_1)\equiv 0\;(\gp_v)$. Also, $m(u_1,p)\equiv
0\;(\gp_v^{2i})$ and, using our conclusion about
$\tr_{\kt_v/k_v}(u_1)$ together with the fact that $p(z)$ is an
Eisenstein polynomial, we deduce that $\n_{\kt_v/k_v}(u_1)\equiv
a_2\;(\gp_v^2)$. But $\ord_{k_v}(a_2)=1$ and
$\ord_{k_v}(\n_{\kt_v/k_v}(u_1))=2\ord_{\ti k_v}(u_1)$ is always
even, so this last congruence is impossible. This contradiction
completes the case (in~rm).

Finally we must deal with the case (rm~ur). Suppose that
$g=(g_1,g_2)\in\Om_{x,v}^1$. There are apparently two
possibilities: either $|\det(g_1)|_{\kt_v}=|\det(g_2)|_v=1$ or
$|\det(g_1)|_{\kt_v}=q_v^{-1}$ and $|\det(g_2)|_v=q_v$. However,
Lemma \ref{urlem} shows that the second possibility cannot occur
and, moreover, that $g_2\in\gl(2)_{\calo_v}$. Thus we may assume,
as usual, that $g_1=a(t_{11},t_{12})n(u_1)$ and $g_2=1$. The
matrices $M_1$ and $M_2$ given by (\ref{expression}) must be
integral and, since $t_{11}t_{12}^{\sig}$ is a unit, this happens
if and only if the quantities enumerated in (\ref{intcond-inrm}) are
all integral. Again assume that $\ord_{\kt_v}(t_{11})=i$ and that
$i>0$. Then $\tr_{\kt_v/k_v}(u_1)+a_1\equiv 0\;(\gp_v^i)$. If $v$
is dyadic then $a_1=-1$ and so this congruence forces
$\tr_{\kt_v/k_v}(u_1)$ to be a unit. However, since $\kt_v/k_v$ is
ramified, $u_1^{\sig}\equiv u_1\;(\ti\gp_v)$ and so
$\tr_{\kt_v/k_v}(u_1)\equiv 2u_1\equiv 0\;(\ti\gp_v)$, which
implies that $\tr_{\kt_v/k_v}(u_1)$ is not a unit. This
contradiction completes that proof in the dyadic case. Now assume
that $v$ is not dyadic. Then $a_1=0$ and so $\tr_{\kt_v/k_v}(u_1)$
is not a unit. We can write $u_1=u_{11}+u_{12}\sqrt{\pi_v}$ with
$u_{11},u_{12}\in\calo_v$ and a suitable choice of uniformizer
$\pi_v$. Since $\tr_{\kt_v/k_v}(u_1)=2u_{11}$, we conclude that
$u_{11}$ is not a unit and hence that $u_1$ is not a unit. However,
$m(u_1,p)=-a_2+\n_{\kt_v/k_v}(u_1)\equiv0\;(\gp_v^i)$ and $a_2$ is
a unit. This contradiction completes the proof in the non-dyadic
case.
\end{proof}
Having completed the verification that $\Om_{x,v}$ is an omega set
in every case, we can now quickly achieve the aim of this section.
\begin{cor}\label{611fulfilled}
Condition \ref{xicoeff} holds. Moreover, $a_{x,v,n}=0$ if $n$ is
odd.
\end{cor}
\begin{proof}
Let $v\in\gM\setminus S_0$ and $y\in V_{k_v}^{\sst}$. From Lemma
\ref{differentreps}, Proposition \ref{omegaint} and the choices
made above we have
\begin{equation}\label{xiasint}
\Xi_{y,v}(s)=Z_{x,v}(\Phi_{v,0},s)=
\int_{\Om_{x,v}}|\chi(g)|_v^s\,dg_v
\end{equation}
where $x$ is the representative we have chosen here to represent
the orbit of $y$. Let
$V_j=\{g\in\Om_{x,v}\big||\chi(g)|_v=q_v^{-j}\}$. From
(\ref{xiasint}) we obtain
\begin{equation*}
\Xi_{y,v}(s)=\sum_{j=-\infty}^{\infty}\vol(V_j)q_v^{-js}\,.
\end{equation*}
However, from (4) in the definition of an omega set,
$V_j=\emptyset$ if $j<0$. Thus the sum really only extends from $0$
to $\infty$ and $a_{y,v,n}=\vol(V_n)$ for $n\geq0$. This makes it
clear that $a_{y,v,n}\geq0$ for all $n$. Since $\chi$ is the square
of a rational character, $V_n=\emptyset$ if $n$ is odd and this
gives the last statement. Finally, again by (4) of the definition,
$V_0=\Om_{x,v}^1=K_v$ and so $a_{y,v,0}=\vol(K_v)=1$.
\end{proof}

\section{The estimate of the local zeta functions}%
\label{sec-estimate}

The purpose of this section is to verify Condition \ref{aproperty}.
So  we assume that $v\in\gM\setminus
S_0$ and $x\in V_{k_v}^{\sst}$. Our method will be to estimate
$\Xi_{x,v}(s)$ by expressing it as an integral over a domain,
$\Gamma_v$, adapted to the purposes of this section as the omega sets
were to those of section \ref{sec-omega}. Throughout this section,
if $T_{x1}$ and $T_{x2}$ are distributions depending on $x$ and
$T_{x1}=C_xT_{x2}$ for some constant $C_x\neq0$ then we shall write
$T_{x1}\propto T_{x2}$. After working with such proportionality
statements, we shall appeal to the results of section
\ref{sec-omega} to strengthen them to inequalities. Thus the
results of this section depend logically on those of the last.

We introduce the following objects ($j\geq 0$ in the last
equation).
\begin{equation}  \label{gamdefn}
\begin{aligned}
\gam & = \begin{cases} 
(a(1,t_1)n(u_1),a(1,t_2)n(u_2),n(u_3)a(t_3,t_4)) &  v\in \gMsp, \\
(a(1,t_1)n(u_1),n(u_2)a(t_2,t_3)) &  v\notin \gMsp, \end{cases} \\
d\gam & = \begin{cases} 
\md t_1 \md t_2 \md t_3 \md t_4 du_1 du_2 du_3 & v\in \gMsp, \\
\md t_1 \md t_2 \md t_3 du_1 du_2  & v\notin \gMsp, 
\end{cases} \\
\Gamma_v & = \begin{cases} 
\{\gam \mid  t_1,t_2,t_3,t_4\in k_v^{\times},\; 
u_1,u_2,u_3\in k_v \} &  v\in \gMsp, \\
\{\gam\mid  t_1\in \ti k_v^{\times},t_2,t_3\in k_v^{\times},\; 
u_1\in \ti k_v,u_2\in k_v \} & v\notin \gMsp, \end{cases} \\
\Gamma_{v}^j & = \begin{cases} 
\{\gam\in\Gamma_v\big| |t_1t_2t_3t_4|_v=q_v^{-j}\} & v\in \gMsp,\\
\{\gam\in\Gamma_v\big| |\normv{t_1}t_2t_3|_v=q_v^{-j}\} & v\notin \gMsp
.\end{cases}
\end{aligned}
\end{equation}
In the above definition, $\md t_1, du_1$, etc., are 
the standard measures on 
$k_v^{\times},\ti k_v^{\times},k_v$, or $\ti k_v$, and $d\gam$ is 
thus a  measure on $\Gam_v$ right invariant with respect to 
the last entry and left invariant with respect to the 
other  entries.   
\begin{lem} \label{lemma:gamma_set}
If {\rm $x\in V^{\sst}_{k_v}$} is a standard orbital representative
then
\begin{equation*}
\int_{G_{k_v}/G^{\circ}_{x\,k_v}}f(g'_{x,v} x)\,dg'_{x,v}
\propto \int_{\Gamma_v}f(\gam x) d\gam 
\end{equation*}
for every $f\in L^1(G_{k_v}x)$ which is invariant on the left by the
action of elements of the form $(1,1,\ka)$ or $(1,\ka)$ 
with $\ka\in\GL{2}_{\co_v}$. 
\end{lem}
\begin{proof}  
We begin with the case $v\in \gMsp$. Define
\begin{equation}
\bar{\Gamma}_v=\left\{\bar \gam\in G_{k_v}\left| \begin{matrix} 
\bar \gam=(a(1,t_1)n(u_1),a(1,t_2)n(u_2),g_3) \\ 
t_1,t_2\in k_v^{\times},\; u_1,u_2\in  k_v \end{matrix}\right. \right\}
.\end{equation}
Suppose that $x=w_p$ where $p(z)=z^2+a_1z+a_2$ (recall that all the
standard orbital representatives have this form). We claim that
$\bar{\Gamma}_v\cap G^{\circ}_{x\, k_v}=\{1\}$ and
that 
\begin{equation*}
\bar{\Gamma}_v G^{\circ}_{x\,k_{v}}
= \{(g_1,g_2,g_3) \mid 
g_{i11}^2+a_1g_{i11}g_{i12}+ a_2g_{i12}^2\neq 0,\,i=1,2\}
.\end{equation*}

The elements of the group $G^{\circ}_{x\, k_v}$ have the form
described in Lemma \ref{stabexpl(1)}.
This makes it clear that $\bar{\Gamma}_v\cap G^{\circ}_{x\,k_v}=\{1\}$. Since
the last entry in elements of $\bar{\Gamma}_v$ is unrestricted, we need
only show that the equation
\begin{equation}  \label{lineq}
\begin{pmatrix} 1 & 0\\ u' & t\end{pmatrix}
\begin{pmatrix} m_{11} & m_{12}\\ m_{21} & m_{22}\end{pmatrix}=
\begin{pmatrix} c & -d\\a_2d & c-a_1d\end{pmatrix}
\end{equation}
is always solvable for $t\not= 0$, $u'$ and $c$ and $d$ satisfying 
$c^2-a_1cd+a_2d^2\not= 0$ provided that
$m_{11}^2+a_1m_{11}m_{12}+a_2m_{12}^2\not= 0$ and the matrix
$(m_{ij})$ is non-singular.

If (\ref{lineq}) holds, 
we must take $c=m_{11}$ and
$d=-m_{12}$ and then the equation is equivalent to
\begin{equation*}
\begin{pmatrix} m_{11} & m_{21}\\ m_{12} & m_{22}\end{pmatrix}
\begin{pmatrix} u'\\ t\end{pmatrix}=
\begin{pmatrix} -a_2m_{12}\\ m_{11}+a_1m_{12}\end{pmatrix}\,,
\end{equation*}
which is solvable for $t$ and $u'$ since the coefficient matrix is
non-singular by hypothesis. If $t=0$ then we have $u'm_{11}=-a_2m_{12}$ and
$u'm_{12}=m_{11}+a_1m_{12}$.  Multiplying the first equation by $m_{12}$,
the second by $m_{11}$ and subtracting, we obtain $m_{11}^2+a_1m_{11}m_{12}
+a_2m_{12}^2=0$, contrary to hypothesis. This proves the second claim.

Let $d_{\ell}\bar\gam=\md t_1\,\md t_2\,du_1\,du_2 dg_3$.  
Then $d_{\ell}\bar\gam$ is a left Haar
measure on the (non-unimodular) group $\bar{\Gamma}_v$.
From what we have just shown, it follows that $G_{k_v}\setminus
\bar{\Gamma}_v\cdot G^{\circ}_{x\,k_v}$ always has measure zero. Thus 
we have
\begin{equation}
\begin{aligned}
\label{eq:measure_theory1}
\int_{G_{k_v}/G^{\circ}_{x\, k_v}}f(g'_{x,v} x)\,dg'_{x,v}
&=\int_{\bar{\Gamma}_v\cdot G^{\circ}_{x\, k_v}/G^{\circ}_{x\, k_v}}
f(g'_{x,v} x)\,dg'_{x,v} \\
&\propto\int_{\bar{\Gamma}_v}f(\gam x)d_{\ell}\gam
\end{aligned}
\end{equation}
for all $f\in L^1(G_{k_v}x)$. 
Now if $\varphi\in
L^1(\GL{2}_{k_v})$ is left invariant under $\GL{2}_{\co_v}$ then
the Iwasawa decomposition implies that
\begin{equation*}
\int_{\GL{2}_{k_v}}\varphi(h)\,dh\propto
\int_{B}\varphi(b)\,d_{r}b
\end{equation*}
where $B=\{n(u_3)a(t_3,t_4)\mid t_3,t_4\in k^{\times},u_3\in k\}$ and
$d_{r}b$ denotes the right Haar measure on the group $B$. It is
easy to check that $d_{r}b=\md t_3\,\md t_4\,du_3$ and applying this
in (\ref{eq:measure_theory1}) we obtain the conclusion.

Finally, almost identical arguments apply in the case when
$(G_{k_v},V_{k_v})$ is not split and we shall not repeat them.
\end{proof}
\begin{prop}
If $p(z)=z^2-z$ then we have
\begin{equation*}
\Xi_{w_p,v}(s)=(1-q_v^{-(2s-1)})^{-1}
(1-q_v^{-(2s-2)})^{-1}\,.
\end{equation*}
\end{prop}
\begin{proof}
Our work will be simplified if we compute 
with the element $x=n_0w_p$ with 
$n_0=(1,1,{}^{t}n(1))$ or $(1,{}^{t}n(1))$ instead of with
the element $w_p$. By Lemma \ref{differentreps},
$Z_{w_p,v}(\Phi_{v,0},s)=Z_{x,v}(\Phi_{v,0},s)$ and so this is
permissible. 
 
Suppose that $v\in\gM_{\text{sp}}$. Then, by Lemma \ref{stabexpl(1)}, 
elements of $G_{x\,k_v}^{\circ}$ have the form
\begin{equation}    \label{xstabform}
\left(\begin{pmatrix} c_{11}  &c_{11}-c_{12}\\
0 & c_{12}\end{pmatrix},
\begin{pmatrix} c_{21} & c_{21}-c_{22}\\
0 & c_{22}\end{pmatrix},*\right)
\end{equation}
where $*$ is determined by the other two entries. 
Note that the  conjugation by $n_0$ does not change 
the first two components.   Let 
\begin{equation}
\begin{aligned}
\mu & = ({}^{t}n(u_1),{}^{t}n(u_2),a(t_1,t_2)n(u_3)), \\
d\mu & = |t_1^{-1}t_2|_v\md t_1\,\md t_2\,du_1\,du_2\,du_3, \\
S & = \{\mu\mid t_1,t_2\in k_v^{\times},\; u_1,u_2,u_3\in k_v\}
.\end{aligned}
\end{equation}
From (\ref{xstabform}) and the Iwasawa decomposition 
it follows that $K_vSG_{x\,k_v}^{\circ}=G_{k_v}$ and
$dg\propto d\ka\,d\mu\, dg''_{x,v}$. 
Since $\Phi_{v,0}$ is $K_v$-invariant, 
\begin{align*}
\Xi_{x,v}(s) & = b_{x,v}\int_{G_{k_v}/G_{x\,k_v}^{\circ}}
|\chi(g'_{x,v})|_v^s\Phi_{v,0}(g'_{x,v} x)\,dg'_{x,v} \\
& \propto\int_{S}|\chi(\mu)|_v^{s}\Phi_{v,0}(\mu x)\,d\mu.
\end{align*}
Computation gives 
\begin{equation*}
\mu x =\left( \begin{pmatrix} t_1 & 0\\ 0 & 0\end{pmatrix},
t_2\begin{pmatrix} u_3-u_1-u_2+u_1u_2+1 & u_1-1\\ 
u_2-1 & 1\end{pmatrix} \right)
.\end{equation*}

Introducing the variables
\begin{equation*}
\bar u_1 = t_2(u_1-1),\; 
\bar u_2 = t_2(u_2-1),\; 
\bar u_3 = t_2(u_3-u_1-u_2+u_1u_2+1)
\end{equation*}
we have $d\bar u_1\,d\bar u_2\,d\bar u_3=|t_2|_v^3du_1\,du_2\,du_3$.
So
\begin{equation*}
\begin{aligned}
\Xi_{x,v}(s) & \propto
\int|t_1|_v^{2s-1}|t_2|_v^{2s-2}\Phi_{v,0}\left(\left(
\begin{smallmatrix} t_1 & 0\\ 0 & 0\end{smallmatrix}\right),
\left( \begin{smallmatrix}\bar u_3 & \bar u_1\\ 
\bar u_2 & t_2\end{smallmatrix} \right)\right)\, 
\md t_1\md t_2d\bar u_1 d\bar u_2 d\bar u_3 \\
&=\int_{|t_1|_v,|t_2|_v\leq 1}
|t_1|_v^{2s-1}|t_2|_v^{2s-2}\,\md t_1 \md t_2\\
& = (1-q_v^{-(2s-1)})^{-1}(1-q_v^{-(2s-2)})^{-1}
.\end{aligned}
\end{equation*}
But we know from Condition \ref{xicoeff} that the constant term 
in $\Xi_{x,v}(s) $ is
$1$ and so $\Xi_{x,v}(s) $ has the stated
value. When $v\in\gM_{\text{in}}$ the
calculation is a simple variation on the above and we shall not
reproduce it here. 
\end{proof}
\begin{prop}\label{splv}
Let {\rm $v\in \gMsp$} 
and suppose that
$x$ is the standard orbital representative for an orbit with
$k_v(x)\neq k_v$. If
\begin{equation*}
L_v(s)=1+ 8(1-q_v^{-2(s-1)})^{-3}q_v^{-2(s-1)}
(4-3q_v^{-2(s-1)}+q_v^{-4(s-1)})
\end{equation*}
then $\Xi_{x,v}(s) \cleq L_v(s)$.  
\end{prop}
\begin{proof}
The standard orbital representative is $x=w_p$ for some quadratic
polynomial $p(z)=z^2+a_2$ which is irreducible over $k_v$ (we may
assume that $a_1=0$ since $v\notin\gM_{\text{dy}}$). 
Let $\gam$, $d\gam$, $\Gam_v$ and $\Gam^j_v$ be as in (\ref{gamdefn}).  
By Definition \ref{lzeta} and Lemma \ref{lemma:gamma_set}, 
\begin{equation*}
\Xi_{x,v}(s) = Z_{x,v}(\Phi_{v,0},s) = C_x \int_{\Gamma_v}
|\chi(\gam)|_v^s\Phi_{v,0}(\gam x) d\gam
\end{equation*}
for some constant $C_x\neq0$. Since $\Gam_v=\coprod_j\Gam_v^j$,
\begin{equation*}
\Xi_{x,v}(s)=C_x\sum_{j=0}^{\infty}q_v^{-2js} 
\int_{\Gamma_v^j}\Phi_{v,0}(\gam x)d\gam
\end{equation*}
which implies that
\begin{equation}\label{eq:equality_for_a}
a_{x,v,2j}=C_x\int_{\Gam_v^j}\Phi_{v,0}(\gam x)\,d\gam
\end{equation}
for all $j\geq0$. (Recall that $a_{x,v,n}=0$ if $n$ is odd by
Corollary \ref{611fulfilled}.)

Computing, we find that $\gam x=(M_1,M_2)$ where
\begin{equation*}
\begin{aligned}
M_1 & = \begin{pmatrix} 0 &  t_2t_3\\ 
t_1t_3 & t_1t_2t_3(u_1+u_2)\end{pmatrix}, \\
M_2 & = \begin{pmatrix}
t_4 & t_2(t_4u_2+t_3u_3)\\
t_1(t_4u_1+t_3u_3) & m(t,u)\end{pmatrix}
\end{aligned}
\end{equation*}
with
\begin{equation*}
m(t,u) =t_1t_2t_3(u_1+u_2)u_3+
t_1t_2t_4(u_1u_2-a_2).
\end{equation*}
If we make $t_1,\dots,t_4$ units and $u_1,\dots,u_3$
integers then $\gam x\in V_{\co_v}$ and the volume of the set
$\{\gam\mid t_j\in\calo_v^{\times},u_j\in\calo_v\}$ under $d\gam$
is $1$ and so 
It follows from this, Condition \ref{xicoeff} and 
(\ref{eq:equality_for_a}) that
\begin{equation*}
1=a_{x,v,0}=C_x\int_{\Gam_v^0}\Phi_{v,0}(\gam x)\,d\gam\geq C_x\,.
\end{equation*}
Therefore, from (\ref{eq:equality_for_a}) again,
\begin{equation}\label{eq:inequality_for_A}
a_{x,v,2j}\leq\int_{\Gamma_v^j}\Phi_{v,0}(\gam x)\,d\gam
\end{equation}
for all $j\geq0$.

We introduce new variables defined by
\begin{equation*}
\bar t_1= t_4,\; 
\bar t_2 = t_2t_3,\; 
\bar t_3 = t_1t_3,\; 
\bar t_4 = t_1t_2t_3t_4\; 
.\end{equation*}
Then 
\begin{equation*}
t_1 = \bar t_1^{-1}\bar t_2^{-1}\bar t_4 ,\; 
t_2 = \bar t_1^{-1}\bar t_3^{-1}\bar t_4,\; 
t_3 = \bar t_1\bar t_2\bar t_3\bar t_4^{-1},\; 
t_4 = \bar t_1 
.\end{equation*}
Note that  $\bar t_1,\dots,\bar t_4$ are monomials 
of $t_1,\dots,t_4$.  So they correspond to a lattice in 
$\Z^4$.  Since the correspondence between $(t_1\ccd t_4)$
and $(\bar t_1\ccd \bar t_4)$ is bijective, this lattice
must be unimodular.  This implies that 
\begin{equation} \label{change-variables-t1}
\md \bar t_1\md \bar t_2\md \bar t_3\md \bar t_4=
\md t_1\md t_2\md t_3\md t_4
.\end{equation}
Suppose that $\gam x\in V_{\co_v}$. Then $\bar t_1,\bar t_2,\bar t_3\in\co_v$.
Since $|P(x)|_v$ is the maximum of $|P(y)|_v$ for $y\in
G_{k_v}x\cap V_{\calo_v}$,
$|P(x)|_v\geq |P(\gam x)|_v=|\bar t_4|_v^2|P(x)|_v$, which implies
that $\bar t_4\in\co_v$. The conditions that the $(2,2)$ entry in $M_1$
and the $(2,1)$ and $(1,2)$ entries in $M_2$ are integers may be
expressed as $N\,{}^{t}(u_1,u_2,u_3)\in\co_v^3$ where
\begin{equation*}
N=\begin{pmatrix}
t_1t_2t_3 & t_1t_2t_3 & 0\\
t_1t_4 & 0 & t_1t_3\\
0 & t_2t_4 & t_2t_3 \end{pmatrix}
= \begin{pmatrix}
\bar t_1^{-1}\bar t_4  & \bar t_1^{-1}\bar t_4  & 0\\
\bar t_2^{-1}\bar t_4  & 0 & \bar t_3\\
0 & \bar t_3^{-1}\bar t_4  & \bar t_2\end{pmatrix}.
\end{equation*}
This matrix factorizes as $N=D_1^{-1}CD_2$ where we have set
$D_1=\diag(\bar t_1,\bar t_2,\bar t_3)$, 
$D_2=\diag(\bar t_4,\bar t_4,\bar t_2\bar t_3)$ and
\begin{equation*}
C=\begin{pmatrix}
1 & 1 & 0\\ 1 & 0 & 1\\ 0 & 1 & 1\end{pmatrix}.
\end{equation*}
Let 
\begin{equation*} 
\pmatrix \bar u_1\\ \bar u_2\\ \bar u_3\endpmatrix 
= CD_2\pmatrix u_1\\ u_2\\ u_3\endpmatrix 
.\end{equation*}
Then the three conditions are equivalent to
${}^t(\bar u_1,\bar u_2,\bar u_3)\in D_1\co_v^3$, 
which in turn is equivalent to the conditions
\begin{equation}\label{eq:u_conditions}
\bar u_1\in \bar t_1\co_v\,,\qquad
\bar u_2\in \bar t_2\co_v\,,\qquad
\bar u_3\in \bar t_3\co_v\,.
\end{equation}
By computation, 
\begin{equation*}
\begin{pmatrix}
u_1\\ u_2\\ u_3\end{pmatrix}= (1/2)
\begin{pmatrix} \bar t_4^{-1}(\bar u_1+\bar u_2-\bar u_3)\\
\bar t_4^{-1}(\bar u_1-\bar u_2+\bar u_3)\\
\bar t_2^{-1}\bar t_3^{-1}(-\bar u_1+\bar u_2+\bar u_3)
\end{pmatrix}
\end{equation*}
and 
\begin{equation}  \label{change-variables-u1}
du_1du_2du_3=|\bar t_2\bar t_3\bar t_4^2|_v^{-1}
d\bar u_1d\bar u_2d\bar u_3
.\end{equation}
The remaining condition for $\gam x\in V_{\co_v}$ is that
$m(t,u)\in\co_v$. Expressing $m(t,u)$ in terms of the coordinates
$(\bar t_1,\bar t_2,\bar t_3,\bar t_4,\bar u_1,\bar u_2,\bar u_3)$ 
we find that
\begin{equation*}
m(t,u)= (1/4) \bar t_1^{-1}\bar t_2^{-1}\bar t_3^{-1}
[-Q(\bar u_1,\bar u_2,\bar u_3)-4\bar t_4^2a_2]
\end{equation*}
where $Q(\bar u_1,\bar u_2,\bar u_3)=\bar u_1^2+\bar u_2^2+\bar u_3^2
-2(\bar u_1\bar u_2+\bar u_1\bar u_3+\bar u_2\bar u_3)$.
Since $v\notin \gMdy$ and $P(x)=-4a_2$, 
$m(t,u)\in\co_v$ if and only if
\begin{equation}
\label{eq:quadratic_cond1}
Q(\bar u_1,\bar u_2,\bar u_3)-\bar t_4^2P(x)\in \bar t_1\bar t_2\bar t_3\co_v.
\end{equation}

We claim that at least one
of $|\bar t_1|_v$, $|\bar t_3|_v$ and $|\bar t_2|_v$ 
must be greater than or equal to
$|\bar t_4|_v$. Suppose to the contrary that
$|\bar t_1|_v,|\bar t_2|_v,|\bar t_3|_v<|\bar t_4|_v$. 
Then $|\bar u_1|_v,|\bar u_2|_v,|\bar u_3|_v<|\bar t_4|_v$ 
also, by (\ref{eq:u_conditions}),
and so $|Q(\bar u_1,\bar u_2,\bar u_3)|_v\leq 
|\bar t_4|_v^2q_v^{-2}$.
Furthermore, since $\bar t_4\in\co_v$, 
\begin{equation*}
|\bar t_1\bar t_2\bar t_3|_v\leq |\bar t_4|_v^3q_v^{-3}
<|\bar t_4|_v^2q_v^{-2}
\end{equation*}
and it follows from (\ref{eq:quadratic_cond1}) that
$|\bar t_4|_v^2|P(x)|_v\leq|\bar t_4|_v^2q_v^{-2}$ 
and so $|P(x)|_v\leq q_v^{-2}$. 
However, by the choice of the standard orbital
representatives, $|P(x)|_v\geq q_v^{-1}$ and we have a
contradiction. This establishes our claim.

Next we claim that $|\bar t_1|_v,|\bar t_2|_v,|\bar t_3|_v
\geq|\bar t_4|_v^2q_v^{-1}$.
Suppose to the contrary that one of these quantities is less than
$|\bar t_4|_v^2q_v^{-1}$. In light of the symmetry between the
roles of the pairs $(\bar t_1,\bar u_1)$, $(\bar t_2,\bar u_2)$ 
and $(\bar t_3,\bar u_3)$ we
may suppose without loss of generality that $|\bar t_3|_v$ is the greatest
of $|\bar t_1|_v$, $|\bar t_2|_v$ and $|\bar t_3|_v$ and that
$|\bar t_1|_v<|\bar t_4|_v^2q_v^{-1}$. By the previous paragraph,
$|\bar t_3|_v\geq |\bar t_4|_v$. Dividing (\ref{eq:quadratic_cond1})
through by $\bar t_3^2$ we obtain
\begin{equation*}
Q({\bar t_3}^{-1}{\bar u_1},{\bar t_3}^{-1}{\bar u_2},
{\bar t_3}^{-1}{\bar u_3})-
({\bar t_3}^{-1}{\bar t_4})^2P(x)\in
{\bar t_3}^{-1}{\bar t_1\bar t_2}\co_v\subseteq
{\bar t_3}^{-1}{\bar t_1}\co_v.
\end{equation*}
We have $\bar u_1/\bar t_3\in (\bar t_1/\bar t_3)\co_v$ 
and so we may drop the terms
involving $\bar u_1/\bar t_3$ to obtain
\begin{equation}  \label{use-hensel1}
({\bar t_3}^{-1}(\bar u_2-\bar u_3))^2
-({\bar t_3}^{-1}{\bar t_4})^2P(x)\in\
{\bar t_3}^{-1}\bar t_1\co_v.
\end{equation}
Now
\begin{equation*}
|({\bar t_3}^{-1}{\bar t_4})^2P(x)|_v\geq 
{|\bar t_3|_v^{-2}}{|\bar t_4|_v^2q_v^{-1}}
>{|\bar t_3|^{-2}_v}{|\bar t_1|_v}
\geq {|\bar t_3|^{-1}_v}{|\bar t_1|_v}
\end{equation*}
and hence $|{\bar t_3}^{-1}({\bar u_2-\bar u_3})|_v^2
=|({\bar t_3}^{-1}{\bar t_4})^2P(x)|_v$.
This implies that $|\bar t_4^{-1}(\bar u_2-\bar u_3)|_v^2
=|P(x)|_v\geq q_v^{-1}$ and so 
$\ord_{k_v}(\bar t_4^{-1}(\bar u_2-\bar u_3))\leq 0$.  By 
(\ref{use-hensel1}), 
\begin{equation*}
({\bar t_4}^{-1}({\bar u_2-\bar u_3}))^2-P(x)\in
{\bar t_4^{-2}}{\bar t_1\bar t_3}\co_v
\subseteq {\bar t_4^{-2}}{\bar t_1}\co_v \subseteq \gp_v^2.
\end{equation*}
These last two facts allow us to apply Hensel's lemma to conclude
that $P(x)\in (k_v^{\times})^2$, which contradicts the assumption that
$k_v(x)\neq k_v$. Thus
$|\bar t_1|_v,|\bar t_2|_v,|\bar t_3|_v\geq|\bar t_4|_v^2q_v^{-1}$, 
as claimed.

Changing variables to $(\bar t_1,\bar t_2,\bar t_3,\bar t_4,
\bar u_1,\bar u_2,\bar u_3)$ in
(\ref{eq:inequality_for_A}) and using (\ref{change-variables-t1}),
(\ref{change-variables-u1}),
we obtain
\begin{align*}
a_{x,v,2j} & \leq \int |\bar t_2\bar t_3\bar t_4^2|_v^{-1}
\md \bar t_1\md \bar t_2\md \bar t_3\md \bar t_4
d\bar u_1d\bar u_2d\bar u_3 \\
&= q_v^{2j} \int|\bar t_2\bar t_3|_v^{-1}
\md \bar t_1\md \bar t_2\md \bar t_3
d\bar u_1d\bar u_2d\bar u_3
\end{align*}
where, on the domain of integration, $|\bar u_i|_v \leq|\bar t_i|_v$
and  $1\geq |\bar t_i|_v\geq |\bar t_4|_v^2q_v^{-1}=q_v^{-2j-1}$ 
for $i=1,2,3$. 
Note that $|\bar t_4|_v=q_v^{-j}$ on $\Gam_v^{j}$.  
Carrying out the integration with respect
to $\bar u_1$, $\bar u_2$ and $\bar u_3$ we get
\begin{align*}
a_{x,v,2j} & \leq q_v^{2j} \int|\bar t_1|_v
\md \bar t_1\,\md \bar t_2\, \md \bar t_3 \\
& \leq q_v^{2j}(1-q_v^{-1})^{-1}
\int_{1\geq|\bar t_3|_v,|\bar t_2|_v\geq q_v^{-2j-1}}
\md \bar t_2\md \bar t_3 \\
& \leq 2q_v^{2j}(2j+2)^2\\
& = 8q_v^{2j}(j+1)^2.
\end{align*}
Note that the volume of the set $\cup_{i=0}^{2j+1} \pi_v^i \co_v^{\times}$ 
is $2j+2$ and $(1-q_v^{-1})^{-1}\leq 2$.  
Put $B_j(v) = 8q_v^{2j}(j+1)^2$. Using the formul\ae
\begin{equation}  \label{formulae}
\begin{aligned} 
\sum_{j=1}^{\infty}q_v^{-js} & =q_v^{-s}(1-q_v^{-s})^{-1}, \\
\sum_{j=1}^{\infty}jq_v^{-js} & =q_v^{-s}(1-q_v^{-s})^{-2}, \\
\sum_{j=1}^{\infty}j^2q_v^{-js} & =
q_v^{-s}(1+q_v^{-s})(1-q_v^{-s})^{-3}
,\end{aligned}
\end{equation}
valid for $\re(s)>0$, we obtain
\begin{equation*}
\sum_{j=1}^{\infty}  B_j(v)q_v^{-2js} = L_v(s)-1
,\end{equation*}
valid for $\re(s)>1$, where $L_v(s)$ is given in the 
statement of the proposition.  
This completes the proof. 
\end{proof}

\begin{prop}\label{nonsplv}
Let {\upshape $v\in\gM_{\text{in}}$} and suppose that
$x$ is the standard orbital representative for an orbit with
$k_v(x)\neq k_v$. If
\begin{equation*}
L_v(s)= 1 + 4(1-q_v^{-2(s-1)})^{-2}q_v^{-2(s-1)}(2-q_v^{-2(s-1)})
\end{equation*}
then $\Xi_{x,v}(s)\cleq L_v(s)$. 
\end{prop}
\begin{proof}
The structure of this proof will be very similar to that of the
proof of proposition \ref{splv} and so we shall abbreviate somewhat.
We have $x=w_p$ for some irreducible quadratic polynomial
$p(z)=z^2+a_2\in k_v[z]$.
Let $\gam$, $d\gam$, $\Gam_v$ and $\Gam_v^j$ be as in (\ref{gamdefn}).  
Arguing as in the previous proposition we obtain the inequality
\begin{equation}
\label{eq:inequality_for_A2}
a_{x,v,2j}\leq \int_{\Gamma_v^j}\Phi_{v,0}(\gam x)\,d\gam
\end{equation}
for all $j\geq0$.

Calculation gives $\gam x=(M_1,M_2)$ where
\begin{equation*}
\begin{aligned}
M_1 & = \begin{pmatrix}
0 & t_1^{\sig}t_2\\ t_1t_2  & t_2\normv{t_1}\trv{u_1}
\end{pmatrix}, \\
M_2 & = \begin{pmatrix}
t_3 & t_1^{\sigma}(t_3u_1^{\sig}+t_2u_2)\\
t_1(t_3u_1+t_2u_2) & m(t,u) \end{pmatrix}
\end{aligned}
\end{equation*}
with 
\begin{equation*}
m(t,u)=
t_2\normv{t_1}\trv{u_1}u_2+
t_3\normv{t_1}[\normv{u_1}-a_2].
\end{equation*}

We introduce new variables defined by
\begin{equation*}
\bar t_1 = t_1t_2,\; \bar t_2 = t_3,\; \bar t_3 = t_2t_3\normv{t_1}
.\end{equation*}
Then 
\begin{equation*}
t_1 = \bar t_1^{-\sigma}\bar t_2^{-1}\bar t_3,\; 
t_2 = \bar t_2\bar t_3^{-1}\normv{\bar t_1},\; 
t_3 = \bar t_2
.\end{equation*}
Since we are dealing with coordinates in two different fields,
$k_v$ and $\ti k_v$, a small digression is required to calculate the
relationship between $\md \bar t_1\md \bar t_2\md \bar t_3$ 
and $\md t_1\md  t_2\md t_3$. Let us fix an element $\beta\in\ti
k_v^{\times}$ which satisfies $\beta^{\sigma}=-\beta$.
For
$u\in\ti k_v$, we define $u^+=u+u^{\sigma}$ and
$u^-=(u-u^{\sigma})/\beta$. Both $u^+$ and $u^-$ lie
in $k_v$ and since $u=(1/2)(u^++\beta u^-)$, $u^+$
and $u^-$ serve as $k_v$ coordinates for $\ti k_v$. 
We use this notation replacing $u$ by  other letters.  
The measure corresponding to $dt_1^+\,dt_1^-$ is 
invariant under addition and hence
there is a constant $C_v$, depending only on $k$, $\ti k$ and $v$,
such that
\begin{equation*}
\md t_1=C_v\frac{dt_1^ + dt_1^-}{|\normv{t_1}|_v}
.\end{equation*}
We also have
$\normv{t_1}= (1/4)[(t_1^+)^2-\beta^2(t_1^-)^2]$ and a
calculation gives
\begin{equation*}
\left|\frac{\partial(\bar t_1^+,\bar t_1^-,\bar t_2,\bar t_3)}
{\partial(t_1^+,t_1^-,t_2,t_3)}\right|_v=|t_3\normv{\bar t_1}|_v
\end{equation*}
so that $d\bar t_1^+d\bar t_1^-\md \bar t_2\md \bar t_3/
|\normv{\bar t_1}|_v= d t_1^+dt_1^-\md t_2\md t_3\, /|\normv{t_1}|_v$. 
Multiplying both sides
by $C_v$ we obtain 
\begin{equation} \label{change-variables-t2}
\md \bar t_1\md \bar t_2\md \bar t_3 =
\md t_1\md t_2\md t_3
.\end{equation}

Suppose that $\gam x\in V_{\co_v}$. Then $\bar t_1\in \ti\co_v$ and
$\bar t_2\in\co_v$. Also $\bar t_3\in\co_v$ by Lemma \ref{lessthanlem}.
If we set
\begin{align*}
\bar u_1 & =\normv{\bar t_1}u_2+\bar t_3u_1\\
\bar u_2 & =\bar t_3\trv{u_1}
\end{align*}
then the $(2,2)$ entry in $M_1$ is $\bar t_2^{-1}\bar u_2$ and the $(2,1)$
entry in $M_2$ is $\bar t_1^{-\sigma}\bar u_1$ and it follows that
\begin{equation*}
\bar u_1\in \bar t_1^{\sigma}\ti\co_v \qquad\text{and}\qquad
\bar u_2\in \bar t_2\co_v
.\end{equation*}
We have
\begin{align*}
u_1 & = (1/2)\bar t_3^{-1}(\bar u_1-\bar u_1^{\sigma}
+\bar u_2), \\
u_2 & = (1/2)\normv{\bar t_1}^{-1}(\bar u_1+\bar u_1^{\sigma}
-\bar u_2)
\end{align*}
and so
\begin{align*}
u_1^+ & = \bar t_3^{-1} \bar u_2, \\
u_1^{-} & = \bar t_3^{-1}\bar u_1^-, \\
u_2 & = (1/2)\normv{\bar t_1}^{-1}(\bar u_1^+-\bar u_2)
.\end{align*}
Hence
$du_1^+du_1^-du_2= |\bar t_3|_v^{-2}
|\normv{\bar t_1}|_v^{-1}d\bar u_1^+d\bar u_1^-d\bar u_2$ 
which implies that  
\begin{equation} \label{change-variables-u2}
du_1du_2=
|\bar t_3|_v^{-2}|\normv{\bar t_1}|_v^{-1}d\bar u_1d\bar u_2
.\end{equation}

The remaining condition for $\gam x\in V_{\co_v}$ is that
$m(t,u)\in \co_v$. In the coordinates 
$(\bar t_1,\bar t_2,\bar t_3,\bar u_1,\bar u_2)$ we have
\begin{equation*}
m(t,u)= (1/4)\bar t_2^{-1}\normv{\bar t_1}^{-1}
[-Q(\bar u_1,\bar u_2)+\bar t_3^2 P(x)]
\end{equation*}
where
\begin{equation*}
Q(\bar u_1,\bar u_2)=\bar u_1^2+\bar u_2^2
+(\bar u_1^{\sigma})^2-2(\bar u_1\bar u_2+\bar u_1^{\sigma}\bar u_2
+\bar u_1\bar u_1^{\sigma})
.\end{equation*} 
Thus $m(t,u)\in\co_v$ if and only if
\begin{equation}
\label{eq:quadratic_cond2}
Q(\bar u_1,\bar u_2)-\bar t_3^2P(x)\in \bar t_2\normv{\bar t_1}\co_v
.\end{equation}

Note that for any $a\in k_v$ we have
$|a|_{\ti k_v}=|a|_v^2$. 
We claim that either
$|\bar t_2|_v\geq|\bar t_3|_v$ or
$|\bar t_1|_{\ti k_v}\geq|\bar t_3|_v^2q_v^{-1}$. 
Suppose to the
contrary that $|\bar t_2|_v\leq|\bar t_3|_vq_v^{-1}$ and
$|\bar t_1|_{\ti k_v}\leq|\bar t_3|_v^2q_v^{-2}$, so that
$|\bar t_2|_{\ti k_v}\leq|\bar t_3|_v^2q_v^{-2}$. 
Then $|\bar u_1|_{\ti k_v},|\bar u_2|_{\ti k_v}
\leq |\bar t_3|^2_{k_v}q_v^{-2}$ and so
$|Q(\bar u_1,\bar u_2)|_{\ti k_v}\leq|\bar t_3|_v^4q_v^{-4}$.
Also
\begin{equation*}
|\bar t_2\normv{\bar t_1}|_{\ti k_v}\leq
|\bar t_3|_v^2q_v^{-2}|\bar t_3|_v^4q_v^{-4}
<|\bar t_3|_v^4q_v^{-4}
.\end{equation*}
So, from (\ref{eq:quadratic_cond2}), $|\bar t_3^2P(x)|_{\ti k_v}\leq
|\bar t_3|_v^4q_v^{-4}$. Thus $|P(x)|_v\leq q_v^{-2}$,
which is a contradiction. The claim follows.

Next we claim that $|\bar t_1|_{\ti k_v}\geq|\bar t_3|_v^4q_v^{-2}$.
Suppose to the contrary that
$|\bar t_1|_{\ti k_v}\leq|\bar t_3|_v^4q_v^{-3}$. Then, from the previous
paragraph, $|\bar t_2|_{\ti k_v}\geq|\bar t_3|_v^2$.  Dividing
(\ref{eq:quadratic_cond2}) by $\bar t_2^2$ we obtain 
\begin{equation*}
Q({\bar t_2}^{-1}{\bar u_2},{\bar t_2}^{-1}{\bar u_1})
-\big({\bar t_2}^{-1}{\bar t_3})^2P(x)
\in {\bar t_2}^{-1}{\normv{\bar t_1}}
\co_v\subseteq {\bar t_2}^{-1}{\bar t_1}\ti\co_v
.\end{equation*} 
Since
$\bar u_1/\bar t_2,\bar u_1^{\sig}/\bar t_2\in (\bar t_1/\bar t_2)\ti\co_v$ 
this containment implies that
\begin{equation}  \label{use-hensel2}
({\bar t_2}^{-1}{\bar u_2})^2
-({\bar t_2}^{-1}{\bar t_3})^2P(x)
\in {\bar t_2}^{-1}{\bar t_1}\ti\co_v
.\end{equation} 
Now 
\begin{equation*}
|({\bar t_2}^{-1}{\bar t_3})^2P(x)|_{\ti k_v}
\geq {|\bar t_2|_{\ti k_v}^{-2}}{|\bar t_3|_v^4q_v^{-2}}>
{|\bar t_2|_{\ti k_v}^{-2}}{|\bar t_1|_{\ti k_v}} 
\geq |\bar t_2|_{\ti k_v}^{-1}{|\bar t_1|_{\ti k_v}}
.\end{equation*} 
Hence
\begin{equation*}
|({\bar t_2}^{-1}{\bar u_2})^2|_{\ti k_v}
=|({\bar t_2}^{-1}{\bar t_3})^2P(x)|_{\ti k_v}
.\end{equation*} 
This implies that
$|\bar u_2/\bar t_3|^2_v=|P(x)|_v\geq q_v^{-1}$ and
so $\ord_{k_v}(\bar u_2/\bar t_3)\leq 0$.  By (\ref{use-hensel2}), 
\begin{equation*}
({\bar t_3}^{-1}{\bar u_2})^2-P(x)
\in {\bar t_3^{-2}}{\bar t_1\bar t_2}\ti\co_v
\subseteq {\bar t_3^{-2}}{\bar t_1}\ti\co_v
.\end{equation*} 
Thus
$|(\bar u_2/\bar t_3)^2-P(x)|_{\ti k_v}
\leq|\bar t_1/\bar t_3^2|_{\ti k_v}\leq q_v^{-3}$
and so $|(\bar u_2/\bar t_3)^2-P(x)|_v\leq q_v^{-2}$. 
We may now apply
Hensel's lemma to conclude that $P(x)\in (k_v^{\times})^2$, which
contradicts the assumption that $k_v(x)\neq k_v$. 
Thus $|\bar t_1|_{\ti k_v}\geq
|\bar t_3|_v^4 q_v^{-2}$.

Changing variables to $(\bar t_1,\bar t_2,\bar t_3,\bar u_1,\bar u_2)$ in
(\ref{eq:inequality_for_A2}) and using (\ref{change-variables-t2}), 
(\ref{change-variables-u2}), we obtain
\begin{align*}
a_{x,v,2j} & \leq \int |\bar t_3|_v^{-2}
|\normv{\bar t_1}|_v^{-1}
\md \bar t_1\md \bar t_2\md \bar t_3d\bar u_1d\bar u_2 \\
&  = q_v^{2j}\int|\normv{\bar t_1}|_v^{-1}
\md \bar t_1\md \bar t_2 d\bar u_1d\bar u_2
\end{align*}
where, on the domain of integration, 
$|\bar u_2|_v\leq|\bar t_2|_v$,
$|\bar u_1|_{\ti k_v}\leq|\bar t_1|_{\ti k_v}=|\normv{\bar t_1}|_v$, 
$|\bar t_2|\leq1$ and
$|\bar t_3|^4_{k_v}q_v^{-2}\leq|\bar t_1|_{\ti k_v}\leq1$.
Carrying out the integration with respect to $\bar u_1$ and $\bar u_2$ we get
\begin{align*}
a_{x,v,2j} & \leq q_v^{2j}\int|\bar t_2|
\md \bar t_1\md \bar t_2 \\
& \leq q_v^{2j}(1-q_v^{-1})^{-1}
\int_{1\geq|\bar t_1|_{\ti k_v}\geq q_v^{-2}|\bar t_3|_v^4}
\md \bar t_1 \\
& \leq 2q_v^{2j}(2j+2) 
\end{align*}
since $\ti k_v/k_v$ is unramified.
Set $B_j(v)=4q_v^{2j}(j+1)$. 
Using (\ref{formulae}),  we obtain
\begin{equation*}
\sum_{j=1}^{\infty}  B_j(v)q_v^{-2js} = L_v(s)-1
,\end{equation*}
valid for $\re(s)>1$, where $L_v(s)$ is given in the 
statement of the proposition.  
This completes the proof. 
\end{proof}

We define 
\begin{equation} \label{lvdefn}
L_v(s) = \begin{cases}
\displaystyle\frac {1+29q_v^{-2(s-1)}-21q_v^{-4(s-1)}+7q_v^{-6(s-1)}} 
{(1-q_v^{-(2s-1)})(1-q_v^{-2(s-1)})^4} & v\in \gM_{\text{sp}},
\vphantom{\intl_{\intl_\int}}\\
\displaystyle\frac {1+6q_v^{-2(s-1)}-3q_v^{-4(s-1)}}
{(1-q_v^{-(2s-1)})(1-q_v^{-2(s-1)})^3} & v\in\gM_{\text{in}}
.\end{cases}
\end{equation}
\begin{prop}
Let $L_v(s)$ be defined by  (\ref{lvdefn}). Then
$\Xi_{x,v}(s)\cleq L_v(s)$ for all $v\in\gM\setminus S_0$ and all
$x\in V_{k_v}^{\text{\upshape ss}}$. The product
$\prod_{v\in\gM\setminus S_0} L_v(s)$ converges absolutely and
locally uniformly in the region $\re(s)>3/2$. Moreover, if
$L_v(s)=\sum_{n=0}^{\infty} \ell_{v,n}q_v^{-ns}$ then
$\ell_{v,0}=1$, $\ell_{v,n}\geq0$ for all $n$ and the series is
convergent in the region $\re(s)>1$. Thus Condition 
\ref{aproperty} is satisfied.
\end{prop}
\begin{proof}
Suppose we have two series
\begin{equation*}
L_{i,v}(s)=1+\sum_{j=1}^{\infty}B_{i,j}(v)q_v^{-js}\,,\qquad 
i=1,2\,,
\end{equation*}
with $B_{i,j}(v)\geq0$ for all $i$ and $j$. Then
\begin{equation*}
L_{1,v}(s)L_{2,v}(s)=1+\sum_{j=1}^{\infty}C_j(v)q_v^{-js}
\end{equation*}
with
\begin{equation*}
C_j(v)=B_{1,j}(v)+B_{2,j}(v)+\sum_{m=1}^{j-1}B_{1,m}(v)B_{2,j-m}(v)
\end{equation*}
and so if we set $L_v(s)=L_{1,v}(s)L_{2,v}(s)$ then
$L_{1,v}(s)\cleq L_v(s)$, $L_{2,v}(s)\cleq L_v(s)$ and
$C_j(v)\geq0$ for all $j$.

We have shown that if $v\in\gMsp$ then
\begin{equation*}
\Xi_{x,v}(\Phi_{v,0},s)=(1-q_v^{-(2s-1)})^{-1}
(1-q_v^{-(2s-2)})^{-1}
\end{equation*}
if $k_v(x)=k_v$ and
\begin{equation*}
\Xi_{x,v}(\Phi_{v,0},s)\cleq(1-q_v^{-2(s-1)})^{-3}
[1+29q_v^{-2(s-1)}-21q_v^{-4(s-1)}+7q_v^{-6(s-1)}]
\end{equation*}
if $k_v(x)\neq k_v$ (the right hand side comes from writing
$L_v(s)$ in Proposition \ref{splv} over a common denominator).
Multiplying these two gives the value of $L_v(s)$ recorded in
(\ref{lvdefn}). The case $v\in\gM_{\text{in}}$ is similar.

From their construction, the series for $L_v(s)$ in (\ref{lvdefn})
have non-negative coefficients and constant term $1$. It follows
by inspection that these series converge when $\re(s)>1$. The
discussion in the first paragraph shows that $\Xi_{x,v}(s)\cleq
L_v(s)$ for all $v\in\gM\setminus S_0$ and $x\in V_{k_v}^{\sst}$.
Finally, it is well-known that the series $\sum_{v\in\gM\setminus
S_0}q_v^{-s}$ is absolutely and locally uniformly convergent in the
region $\re(s)>1$. The usual convergence test for products now
shows that $\prod_{v\in\gM\setminus S_0}L_v(s)$ has the stated
convergence properties.
\end{proof}

\section{The volume of the Integral Points of the Stabilizer}%
\label{sec-iv}

Suppose that $v\in\gM_{\text{f}}$ and that $x\in V_{k_v}^{\sst}$ is
a standard orbital representative. In this section we shall compute
the volume of the intersection $K_v\cap G_{x\,k_v}^{\circ}$ under
the canonical measure $dg_{x,v}''$ introduced in Definition 
\ref{measure''}, unless $v\in\gM_{\text{dy}}$ and $x$ has type
(rm~rm~ur) or (rm~rm~rm). These cases, which involve additional
technicalities, will be dealt with in \cite{kable-yukie-pbh-II}.

As in section \ref{final}, we let $\Del_{\ti
k_v/k_v}=\gp_v^{\ti\del_v}$. We shall use the notation for elements
of the stabilizer introduced in section \ref{space}. Note that
since $w,w_p\in V_{\calo_v}$ in all cases, $G_w$ and $G_{w_p}$ are
defined over $\calo_v$. The canonical measure on
$G_{x\,k_v}^{\circ}$ is induced via the map $\theta:
G_{x\,k_v}^{\circ}\to H_{x\,k_v}$ from the measure $dt_{x,v}$ (we
are suppressing the element $g_x$ since it plays no essential role
here). Now $\theta(K_v\cap G_{x\,k_v}^{\circ})\subseteq
H_{x\,\calo_v}$ and the measure of $H_{x\,\calo_v}$ under
$dt_{x,v}$ is $1$. Below we shall consider the subgroup
$\theta(K_v\cap G_{x\,k_v}^{\circ})$ of $H_{x\,\calo_v}$ in order to
compute $\vol(K_v\cap G_{x\,k_v}^{\circ})$.
\begin{lem}\label{inKcond}
Let $x=w_p$ with $p(z)=z^2+a_1z+a_2$ be a standard orbital
representative. If {\rm $v\in\gM_{\text{sp}}$} then
$g=(A_p(c_1,d_1),A_p(c_2,d_2),A_p(c_3,d_3))\in G_{x\,k_v}^{\circ}$
lies in $K_v$ if and only if $c_1,c_2,d_1,d_2\in\calo_v$ and $\det
A_p(c_j,d_j)\in\calo_v^{\times}$ for $j=1,2$. If
{\rm $v\notin\gM_{\text{sp}}$} then $g=(A_p(c_1,d_1),A_p(c_2,d_2))\in
G_{x\,k_v}^{\circ}$ lies in $K_v$ if and only if
$c_1,d_1\in\ti\calo_v$ and $\det A_p(c_1,d_1)\in\ti
\calo_v^{\times}$.
\end{lem}
\begin{proof}
The conditions on $c_j$, $d_j$ and $\det A_p(c_j,d_j)$ which are
proposed in the statement simply say that $A_p(c_j,d_j)$ lies in
the standard maximal compact subgroup in that factor for $j=1,2$
(respectively $j=1$) and so they are certainly necessary. To see
that they are also sufficient we must show that they imply that the
last entry in $g$ also lies in the relevant maximal compact
subgroup. But this follows immediately from Lemmas \ref{stabexpl(1)}
and \ref{stabexpl(2)}.
\end{proof}

In the following, when we consider $w_p$ we let
$\al=\{\al_1,\al_2\}$ be the set of roots of $p(z)$ as usual. Also,
when $k_v(w_p)/k_v$ is quadratic, we denote by $\nu$ the
non-trivial element of $\gal(k_v(w_p)/k_v)$.
\begin{prop}\label{easyvolume}
Suppose that the type of the orbit containing the standard
representative $x$ is on the list {\upshape(sp), (in),
(rm), (sp ur), (sp rm)}.
Then 
$\vol(K_v\cap G_{x\,k_v}^{\circ})=1$.
\end{prop}
\begin{proof}
We have to show that in every case listed above $\theta(K_v\cap
G_{x\,k_v}^{\circ})=H_{x\,\co_v}$. Of the triple (sp), (in), (rm)
we deal with the first, since the other two are very similar. For
this orbit $p(z)=z^2-z$, we have
$H_{x\,\calo_v}=(\calo_v^{\times})^4$ and the map $\theta$ is
$\theta(A_p(c_1,d_1),A_p(c_2,d_2),A_p(c_3,d_3))=
(c_1,c_1+d_1,c_2,c_2+d_2)$. Note that in this case we may choose
$h_{\al} = \left(\begin{smallmatrix} 1 & -1\\ 
0 & 1\end{smallmatrix}\right)$ in (\ref{stabelm1}).  
Then computing $g_ps_x(t_x)g_p^{-1}$ in (\ref{stabelm1})
explicitly for $t_x=(c_1,c_1+d_1,c_2,c_2+d_2)$, we get the formula
for $\theta$.  If $h=(t_1,t_2,t_3,t_4)\in
H_{x\,\calo_v}$ then $\theta(g)=h$ where $g=(A_p(t_1,t_2-t_1), 
A_p(t_3,t_4-t_3),\ast)$ and, since $\det A_p(t_i,t_j-t_i)=t_it_j$,
it follows from Lemma \ref{inKcond} that $g\in K_v$.

This leaves the cases (sp~ur) and (sp~rm). 
Recall that $p(z)$ was chosen
so that
$\co_{k_v(x)}=\co_v[\al_1]$. We have $H_{x\,
\co_v}=(\co_{k_v(x)}^{\times})^2$ and the map $\theta$ is
\begin{equation*}
\theta(A_p(c_1,d_1),A_p(c_2,d_2),*)=(c_1+d_1\al_1,c_2+d_2\al_1)\,.
\end{equation*}
If $h=(t_1,t_2)\in H_{x\,\co_v}$ then 
$t_1,t_2\in\co_{k_v(x)}^{\times}$
and so we may find $c_1,d_1,c_2,d_2\in\co_v$ so that 
$t_1=c_1+d_1\al_1$ and $t_2=c_2+d_2\al_1$. If we set
$g=(A_p(c_1,d_1),A_p(c_2,d_2),*)\in G_{x\,k_v}^{\circ}$ then
$\det A_p(c_j,d_j)=(c_j+d_j\al_1)(c_j+d_j\al_2)=t_jt_j^{\nu}\in
\co_{v}^{\times}$.   Thus $g\in K_v$ and
$\theta(g)=h$.
\end{proof}
\begin{prop}\label{iv-volume-(in_ur)}
If the type of the orbit containing the standard representative
$x$ is {\upshape(in ur)} then $\vol(K_v\cap G_{x\,k_v}^{\circ})=1$.
If it is {\upshape(rm rm)*} then $\vol(K_v\cap
G_{x\,k_v}^{\circ})=(1-q_v^{-1})^{-1}q_v^{-\tl{\delta}_v}$.
\end{prop}
\begin{proof}
For these cases $H_{x\;\co_v}=(\ti\co_v^{\times})^2$ and $\theta$ is
$\theta(A_p(c,d),*)=(c+d\al_1,c+d\al_2)$. 
Note that we get this formula by explicitly computing 
$g_ps_x(t_x)g_p^{-1}$ for $t_x=(c+d\al_1,c+d\al_2)$ 
in (\ref{stabelm3}).  
Given $t_1,t_2\in\ti\co_v^{\times}$ 
we wish to determine whether there are
integers $c,d\in\ti\co_v$ such that
$\theta(A_p(c,d),*)=(t_1,t_2)$, because
$(t_1,t_2)\in\theta(K_v\cap G_{x\,k_v}^{\circ})$ if and only if
this is possible. Elementary linear algebra shows that
\begin{equation*}
c=(\al_2t_1-\al_1t_2)/(\al_2-\al_1)\,,\qquad
d=(t_2-t_1)/(\al_2-\al_1)\,.
\end{equation*}
Thus $(t_1,t_2)\in\theta(K_v\cap G_{x\,k_v}^{\circ})$ if and only
if $\al_2-\al_1$ divides $t_2-t_1$ and $\al_2t_1-\al_1t_2$. If
$\kt_v=k_v(x)$ is unramified over $k_v$ then $\al_2-\al_1$ is a
unit and this is no condition on $t_1$ and $t_2$. Thus
$\theta(K_v\cap G_{x\,k_v}^{\circ})=H_{x\,\co_v}$ in this case
and the first claim follows.

If $\kt_v=k_v(x)$ is ramified over $k_v$ then
$(\al_2-\al_1)=\ti\gp_v^{\tl{\delta}_v}$ and so we must have $t_1\equiv
t_2 \; (\ti\gp_v^{\tl{\delta}_v})$. With this condition
\begin{equation*}
\al_2t_1-\al_1t_2=(\al_2-\al_1)t_1+\al_1(t_1-t_2)\in
\ti\gp_v^{\tl{\delta}_v}
\end{equation*}
and so the second condition also holds. Thus
\begin{align*}
\theta(K_v\cap G_{x\,k_v}^{\circ})&=
\{(t_1,t_2)\in H_{x\,\co_v}\mid
t_1\equiv t_2\; (\ti\gp_v^{\tl{\delta}_v})\} \\
&\cong\{(\bar t_1,\bar t_2)\in H_{x\,\co_v}\mid
\bar t_1\equiv 1\; (\ti\gp_v^{\tl{\delta}_v})\}\,,
\end{align*}
where the isomorphism is by the measure preserving map
$(t_1,t_2)\mapsto(t_1t_2^{-1},t_2)$. Thus $\vol(K_v\cap
G_{x\,k_v}^{\circ})=\vol(1+\ti\gp_v^{\tl{\delta}_v})$ under the
normalized multiplicative Haar measure on $\kt_v^{\times}$. The
second claim follows. \end{proof}

In the remaining cases, the fields $\kt_v$ and $k_v(x)$ are
distinct quadratic extensions of $k_v$ and we shall assume this to
be so for the rest of this section. 
\begin{lem}\label{imageofK}
We have $\theta(K_v\cap G_{x\,k_v}^{\circ})=
\ti\co_v[\al_1]^{\times}$.
\end{lem}
\begin{proof}
The map $\theta$ is
$\theta(A_p(c,d),*)=c+d\al_1\in\kt_v(x)^{\times}$
as pointed out in the proof of Lemma \ref{stabexpl(2)}.
If $\ka\in K_v\cap
G_{x\,k_v}^{\circ}$ and $\ka=(A_p(c,d),*)$ then $c,d\in\ti\co_v$ and
so $\theta(\ka)\in\ti\co_v[\al_1]$. Moreover, $\det
A_p(c,d)=\theta(\ka)\theta(\ka)^{\nu}$ and
since $\det
A_p(c,d)\in\ti\co_v^{\times}$ and $\ti\calo_v[\al_1]$ is stable
under $\nu$, it follows that $\theta(\ka)$ is a unit
in the ring $\ti\co_v[\al_1]$. Thus $\theta(K_v\cap
G_{x\,k_v}^{\circ})\subseteq\ti\co_v[\al_1]^{\times}$. Suppose that
$c+d\al\in\ti\co_v[\al_1]^{\times}$. Then $c,d\in\ti\co_v$ and $\det
A_p(c,d)=(c+d\al_1)(c+d\al_2)\in\ti\co_v^{\times}$ and it follows
from Lemma \ref{inKcond} that
$(A_p(c,d),*)\in K_v\cap G_{x\,k_v}^{\circ}$. This establishes
the reverse inclusion. \end{proof}
\begin{prop}\label{iv-volume-(rm_ur)}
If the type of the orbit containing the standard representative
$x$ is {\upshape (rm ur)} or {\upshape (in rm)} then $\vol(K_v\cap
G_{x\,k_v}^{\circ})=1$.
\end{prop}
\begin{proof}
If the type is (rm~ur) then $p(z)\in k_v[z]$ either has the form
$p(z)=z^2-r$ with $r$ a non-square unit, if
$v\notin\gM_{\text{dy}}$, or is an Artin-Schreier polynomial, if
$v\in\gM_{\text{dy}}$. By hypothesis, $p(z)$ is irreducible when
regarded as an element of $\ti k_v[z]$. Thinking of $\ti k_v$ as
the ground field, these facts imply that $\calo_{\ti
k_v(x)}=\ti\calo_v[\al_1]$ and so $\theta(K_v\cap
G_{x\,k_v}^{\circ})=\calo_{\ti k_v(x)}^{\times}=H_{x\,\calo_v}$, by
Lemma \ref{imageofK}.

If the type is {\upshape (in rm)} then $p(z)\in k_v[z]$ is an
Eisenstein polynomial and, since $\kt_v/k_v$ is unramified, $p(z)$
is still an Eisenstein polynomial when regarded as an element of
$\kt_v[z]$. Thus $\co_{\kt_v(x)}=\ti\co_v[\al_1]$ and, again by the
Lemma, $\theta(K_v\cap G_{x\,k_v}^{\circ})=H_{x\,\co_v}$.
\end{proof}

We are left with the type (rm~rm~ur) with $v\notin\gM_{\text{dy}}$. In
this case $\ti k_v(x)=k_v(\sqrt{\pi_v},\sqrt{\eta})$ where $\pi_v$
is a uniformizer and $\eta$ is a non-square unit. We may assume
without loss of generality that $\ti k_v=k_v(\sqrt{\pi_v})$ and
$k_v(x)=k_v(\sqrt{\pi_v\eta})$.
\begin{prop}\label{iv-volume-(rm_rm_ur)}
If {\rm $v\in\gM_{\text{f}}\setminus\gM_{\text{dy}}$} and $x$ 
is the standard orbital representative for an orbit with
type {\upshape (rm rm ur)} then $\vol(K_v\cap G_{x\,k_v}^{\circ})=
(q_v+1)^{-1}$.   
\end{prop}
\begin{proof}
We may choose $p(z)=z^2-\pi_v\eta$ as the Eisenstein polynomial
associated to the extension $k_v(x)/k_v$. Then
$\al_1=\sqrt{\pi_v\eta}$ and so $\theta(K_v\cap
G_{x\,k_v}^{\circ})= \ti\calo_v[\sqrt{\pi_v\eta}]^{\times}$ by
Lemma \ref{imageofK}. Since $\ti\calo_v=\calo_v[\sqrt{\pi_v}]$, we
have
\begin{equation*}
\ti\calo_v[\sqrt{\pi_v\eta}]^{\times}=
\{e_1+e_2\pi_v\sqrt{\eta}+e_3\sqrt{\pi_v}+e_4\sqrt{\pi_v\eta}
\mid
e_1\in\calo_v^{\times},e_2,e_3,e_4\in\calo_v\}\,.
\end{equation*}
On the other hand, $\calo_{\ti k_v(x)}=\ti\calo_v[\sqrt{\eta}]$ and
so $H_{x\,\calo_v}=\calo_{\ti k_v(x)}^{\times}$ is
\begin{equation*}
\{f_1+f_2\sqrt{\eta}+f_3\sqrt{\pi_v}+f_4\sqrt{\pi_v\eta}
\mid
f_1,\dots,f_4\in\calo_v,\text{ $f_1$ or $f_2$ a unit}\}\,.
\end{equation*}
From this it is clear that
\begin{equation*}
[H_{x\,\calo_v}:\theta(K_v\cap G_{x\,k_v}^{\circ})]=q_v+1
\end{equation*}
and so the volume has the indicated value.
\end{proof}

\section{Orbital Volumes at the Finite Places}%
\label{bxsif} 

In this section we compute $\vol(K_vx)$ for all standard orbital
representatives, $x$, when $v$ is a finite, non-dyadic place. When
$v$ is dyadic we compute the sum of $\vol(K_vx)$ over the
equivalence class of $x$ under the relation $\asymp$ introduced in
section \ref{final}, unless $x$ has type (rm~rm)*, (rm~rm~ur) or
(rm~rm~rm). These cases will be treated in
\cite{kable-yukie-pbh-II}. Throughout this section, $v$ will be a
finite place of $k$.

Our strategy is to find a subset $\cD$ of $K_vx$, defined by
congruence conditions, whose translates cover the orbit and then to
find its stabilizer modulo a certain high power of the prime ideal.
Since the entries of elements of $V$ may lie in different fields we
need a notation for congruences which takes this into account.
Suppose that $v\notin\gM_{\text{sp}}$. We use the coordinate system
(\ref{case(2)coords}) on $V_{k_v}$. If $x=(x_{ij})$ and
$y=(y_{ij})$ are written with respect to this coordinate system
then $x\equiv y\;(\gp_v^{n_1},\ti\gp_v^{n_2})$ means that
$x_{i0}\equiv y_{i0}\;(\gp_v^{n_1})$, $x_{i2}\equiv
y_{i2}\;(\gp_v^{n_1})$ and $x_{i1}\equiv y_{i1}\;(\ti\gp_v^{n_2})$
for $i=1,2$. If $v\in\gM_{\text{sp}}$ then $x\equiv y\;(\gp_v^n)$
will signify congruences to the modulus $\gp_v^n$ on all the
corresponding entries of the two pairs of matrices.

In addition to this, we shall require some notation from earlier
sections, which we now briefly recall. We put
$2\calo_v=\gp_v^{m_v}$ (so that $m_v=0$ unless $v$ is dyadic) and
if $k_v(x)/k_v$ is quadratic then we put
$\Del_{k_v(x)/k_v}=\gp_v^{\del_{x,v}}$. If $p(z)=z^2+a_1z+a_2\in
k_v[z]$ then $\al=\{\al_1,\al_2\}$ will be its set of roots. Our
standard orbital representative for orbits with $k_v(x)=k_v$ is
$x=w_p$ with $p(z)=z^2-z$. However, we shall find it convenient to
make use of $x=w$ instead. Note that if $g_p$ is as in 
(\ref{g_pdefn}) then $g_p\in K_v$ and so $K_v w_p=K_v w$ and this
substitution is permissible. So we shall use
$w$ as the standard orbital representative of the orbit of $w$.  

Since $w,w_p\in G_{\co_v}$, $G_w$ and $G_{w_p}$ are defined 
over $\co_v$.  So for
any $j\geq 0$ we may consider $G_{x\,\co_v/\gp_v^{j+1}}$ and
$G^{\circ}_{x\,\co_v/\gp_v^{j+1}}$.
We may also  consider the reduction $\bar x$ of $x$ 
modulo $\gp_v^{j+1}$ and its
stabilizer,  $G_{\bar x\, \co_v/\gp_v^{j+1}}$.   In general, 
$G_{x\, \co_v/\gp_v^{j+1}}$ 
and $G_{\bar x\, \co_v/\gp_v^{j+1}}$ may not coincide
and we have to proceed carefully. However,
$G_{x\,\calo_v/\gp_v^{j+1}}\subseteq G_{\bar
x\,\calo_v/\gp_v^{j+1}}$ in all cases. 

The following lemma describe $G_{x\,\calo_v/\gp_v^{j+1}}$. 
when $v\in\gM_{\text{sp}}$.  
\begin{lem} \label{stabreduction(1)}
Suppose $v\in\gM_{\text{{\upshape sp}}}$ and 
$x=w_p$ is a standard orbital
representative.  Then $G^{\circ}_{x\, \co_v/\gp_v^{j+1}}$ 
consists of elements of the form 
{\rm (\ref{stabform(1)})} such that 
\begin{equation*}
\det A_p(c_1,d_1),\det A_p(c_2,d_2) 
\in (\co_v/\gp_v^{j+1})^{\times}
\end{equation*}
and $c_3,d_3\in \co_v/\gp_v^{j+1}$ are related to 
$c_1,c_2,d_1,d_2\in\co_v/\gp_v^{j+1}$  by {\rm (\ref{explicit(1)})}.  
Moreover,
\begin{equation*}
[G_{x\, \co_v/\gp_v^{j+1}}:
G^{\circ}_{x\,\co_v/\gp_v^{j+1}}]=2
\end{equation*}
and the non-trivial class in
$G_{x\,\co_v/\gp_v^{j+1}}/
G^{\circ}_{x\, \co_v/\gp_v^{j+1}}$ is represented by 
$(\tau_p,\tau_p,\tau_p)$.     
\end{lem}
\begin{proof} We briefly sketch the proof here.  
The conditions in Lemma \ref{stabexpl(1)} determine
the structure of $G_x$ as a scheme over $\co_v$
regarding $c_1,d_1,c_2,d_2$ as variables.    
Note that the inequalities $A_p(c_1,d_1),A_p(c_2,d_2)\not=0$
can be regarded as equations $x_1A_p(c_1,d_1)=x_2A_p(c_2,d_2)=1$
after adding variables $x_1,x_2$.  Then $G_{x\, \co_v/\gp_v^{j+1}}$
is by definition the set of $(\co_v/\gp_v^{j+1})$-valued points of this 
scheme and so the equation \ref{explicit(1)} 
will be regarded as an equation over 
$\co_v/\gp_v^{j+1}$ and the above inequalities, after reduction modulo
$\gp_v^{j+1}$, can be regarded  as the condition that 
$A_p(c_1,d_1),A_p(c_2,d_2)$ are units.  Since $(\tau_p,\tau_p,\tau_p)$
is defined over $\co_v$, this proves the lemma.   
\end{proof} 

Similarly, the following lemma follows from Lemma \ref{stabexpl(2)}
and we will not give a proof.  
\begin{lem} \label{stabreduction(2)}
Suppose $v\notin\gM_{\text{{\upshape sp}}}$ 
and $x=w_p$ is a standard orbital
representative.  Then $G_{x\,\co_v/\gp_v^{j+1}}$ consists of
elements of the form  {\rm (\ref{stabform(2)})}
such that $\det A_p(c_1,d_1)\in (\ti\co_v/\ti\gp_v^{j+1})^{\times}$ or 
$(\ti\co_v/\ti\gp_v^{2(j+1)})^{\times}$, 
according as $\ti k_v/k_v$ is ramified
or unramified,  and $c_2,d_2\in \co_v/\gp_v^{j+1}$ 
are related to $c_1,d_1\in \ti\co_v/\ti\gp_v^{j+1}$ or 
$\ti\co_v/\ti\gp_v^{2(j+1)}$
by {\rm (\ref{explicit(2)})}.
Moreover, 
\begin{equation*}
[G_{x\, \co_v/\gp_v^{j+1}}
:G^{\circ}_{x\, \co_v/\gp_v^{j+1}}]=2
\end{equation*}
and the non-trivial class in
$G_{x\, \co_v/\gp_v^{j+1}}/G^{\circ}_{x\, \co_v/\gp_v^{j+1}}$ 
is represented by $(\tau_p,\tau_p)$.    
\end{lem}
\begin{defn}
For $x=(x_1,x_2)\in V_{\co_v/\gp_v^{j+1}}$ let 
$\xspan(x)$ be the $\co_v/\gp_v^{j+1}$-module generated by $x_1$
and $x_2$. 
\end{defn}
The following simple observation will be useful below. We will
usually apply it with $x$ being the reduction modulo $\gp_v^{j+1}$
of a standard orbital representative.
\begin{lem}\label{observation}
Let $v\in\gM_{\text{{\upshape sp}}}$. Given
$g_1,g_2\in\gl(2)_{\calo_v/\gp_v^{j+1}}$, there exists an element
$g_3\in\gl(2)_{\calo_v/\gp_v^{j+1}}$ such that $(g_1,g_2,g_3)\in
G_{x\,\calo_v/\gp_v^{j+1}}$ if and only if
$\xspan((g_1,g_2,1)x)=\xspan(x)$. Similarly, let
$v\notin\gM_{\text{{\upshape sp}}}$. Given
$g_1\in\gl(2)_{\ti\calo_v/\ti\gp_v^{j+1}}$ or
$\gl(2)_{\ti\calo_v/\ti\gp_v^{2(j+1)}}$, according as
$v\in\gM_{\text{{\upshape in}}}$ or $v\in\gM_{\text{{\upshape rm}}}$,
there exists
$g_2\in\gl(2)_{\calo_v/\gp_v^{j+1}}$ such that $(g_1,g_2)\in
G_{x\,\calo_v/\gp_v^{j+1}}$ if and only if
$\xspan((g_1,1)x)=\xspan(x)$.
\end{lem}

Before continuing, we would like to mention a simple fact which we
shall have to use frequently below. If $v$ is a finite, non-dyadic
place of $k$ and $v\notin\gM_{\text{sp}}$ then any class in
$\ti\calo_v/\ti\gp_v^j$ which is fixed by $\sigma$ has a
representative which lies in $\calo_v$. To see this, simply observe
that if $u\in\ti\calo_v$ represents such a class then
$u^{\sigma}\equiv u\;(\ti\gp_v^j)$ and so $u'=\tfrac12(u+u^{\sigma})$
satisfies $u'\equiv u\;(\ti\gp_v^j)$ and $u'\in\calo_v$. Of course, the
corresponding claim when $v$ is a dyadic place is false.
\begin{prop}\label{coinstab}
If $k_v(x)/k_v$ is not ramified then
$G_{\bar x\,\co_v/\gp_v^{j+1}}=G_{x\,\co_v/\gp_v^{j+1}}$.  
\end{prop}
\begin{proof}
We first consider the case $x=w$. 
Suppose $v\in \gM_{\text{sp}}$.  Let $g= (g_1,g_2,g_3)\in G_{\bar w\,
\calo_v/\gp_v^{j+1}}$
be as in (\ref{groupcoords})
where the entries are elements of $\co_v/\gp_v^{j+1}$.  
By computation,  
\begin{equation*} 
(g_1,g_2,1)w  = \left(
\pmatrix g_{111}g_{211} & g_{111}g_{221}\\ 
g_{121}g_{211} & g_{121}g_{221}\endpmatrix, 
\pmatrix g_{112}g_{212} & g_{112}g_{222}\\ 
g_{122}g_{212} & g_{122}g_{222}\endpmatrix\right)
.\end{equation*}
Since $\xspan((g_1,g_2,1)\bar w )=\xspan(\bar w)$, 
\begin{equation*} 
g_{111}g_{221}=g_{121}g_{211}=g_{112}g_{222}=g_{122}g_{212}=0
.\end{equation*}
If $g_{111}$ is a unit then $g_{221}=0$.
Thus $g_{211},g_{222}$ are units.
So $g_{112}=g_{121}=0$.  This implies $g_{122}$ is a unit and so
$g_{212}=0$.  If $g_{111}$ is not a unit, $g_{121}$ must be a unit.
By a similar argument, we can conclude that 
$g_{111}=g_{122}=g_{211}=g_{222}=0$ and $g_{112},g_{121},g_{212},g_{221}$
are units.  Multiplying by an element of 
$G_{x\,\co_v/\gp_v^{j+1}}$ if necessary,
we may assume that $g_1=g_2=1$.  Then it is easy to see that 
$g_3=1$. It follows that $G_{\bar
w\,\calo_v/\gp_v^{j+1}}=G_{w\,\calo_v/\gp_v^{j+1}}$.
The case $v\notin\gM_{\text{sp}}$ is similar.    

We now assume $k_v(x)/k_v$ is quadratic and unramified.  
Suppose $v\in  \gM_{\text{sp}}$.  
Let $g=(g_1,g_2,g_3)\in G_{\bar x\, \co_v/\gp_v^{j+1}}$. 
We choose representatives of $g_1,g_2,g_3$ in $\gl(2)_{\co_v}$
and use the same notation.  
Since $k_v(x)/k_v$ is unramified, if $c,d\in \co_v$,
$c+d\al_1$ is a unit if and only if 
either $c$ or $d$ is a unit.   
The $(1,1)$-entry of $A_p(c_i,d_i)g_i$ is $c_ig_{i11}-d_ig_{i21}$
and the $(1,2)$-entry is $c_ig_{i12}-d_ig_{i22}$. Since 
either $g_{i12}$ or $g_{i22}$ is a unit, we may choose
$c_i=g_{i22}$ and $d_i=g_{i12}$ to make the $(1,2)$-entry of 
$A_p(c_i,d_i)g_i$ zero.  Replacing $g_i$ by an element of the form
$(A_p(c_1,d_1), A_p(c_2,d_2), A_p(c_3,d_3))g$ 
(with $c_3,d_3$ determined by $c_1,c_2,d_1,d_2$), we may assume 
that
\begin{equation}\label{cosetrep}
g_1= \pmatrix 1 & 0\\ u_1 & t_1\endpmatrix,\;  
g_2= \pmatrix 1 & 0\\ u_2 & t_2\endpmatrix
.\end{equation}
If $x=w_p$ and 
$y = \left(\begin{smallmatrix} y_0 & y_1\\ 
y_1^{\sig} & y_2 \end{smallmatrix}\right)$, 
$\bar y\in\xspan(\bar x)$ if and only if $y_1=y_1^{\sig}$ and 
$y_2-a_1y_1+a_2y_0=0$.  By computation,
$(g_1,g_2,1)x = (M_1,M_2)$ where 
\begin{equation} \label{module}
\begin{aligned}
M_1 & = \pmatrix 0 & t_2\\ t_1 & a_1t_1t_2+t_2u_1+t_1u_2 \endpmatrix, \\
M_2 & = \pmatrix 1 & a_1t_2+u_2\\ a_1t_1+u_1 & 
u_1u_2+a_1(t_2u_1+t_1u_2)+(a_1^2-a_2)t_1t_2\endpmatrix 
.\end{aligned}
\end{equation} 
Therefore $t_1=t_2,\; u_1=u_2$, and  
\begin{equation} \label{urstab}
a_1(t_1-1)+2u_1=0,\; 
u_1^2+a_1(2t_1-1)u_1+a_1^2(t_1^2-t_1)+a_2(1-t_1^2)=0
.\end{equation} 
If $v$ is not dyadic, $a_1=0$.  So $u_1=t_1^2-1=0$.  
Then $g_1,g_2$ are both the identity matrix or both $-\tau_p$. 
If $v$ is dyadic, $p$ is an Artin-Schreier 
polynomial.  This means $a_1=-1$ and $a_2$ is a unit.  
By the first condition of (\ref{urstab}),  $t_1=2u_1+1$.  
Substituting in the second condition and simplifying, 
$(4a_2-1)(u_1^2+u_1)=0$.  Since $4a_2-1$ is a unit, 
$u_1^2+u_1=0$.  If $u_1$ is a unit, $u_1=-1$, 
and if $u_1$ is not a unit, 
$u_1=0$.  Note that this argument is valid in 
$\co_v/\gp_v^{j+1}$ (even though it is not a field). 
If $u_1=-1$ then $t_1=-1$ and $g_1,g_2$ are both $-\tau_p$.
If $u_1=0$ then $t_1=1$ and $g_1,g_2=1$.  Then since
$(-1,-1,1),(\tau_p,\tau_p,\tau_p)\in G_{x\, \co_v/\gp_v^{j+1}}$, 
$(-\tau_p,-\tau_p,\tau_p)\in G_{x\, \co_v/\gp_v^{j+1}}$.  
This proves the proposition when $v\in\gM_{\text{sp}}$.

Suppose now that $v\notin\gM_{\text{sp}}$. 
Let $g= (g_1,g_2)\in G_{\bar x\,\calo_v/\gp_v^{j+1}}$ be as 
in (\ref{groupcoords})
where entries of $g_2$ are in $\co_v/\gp_v^{j+1}$
and entries of $g_1$ are in $\ti\co_v/\ti\gp_v^{j+1}$ or
$\ti\co_v/\ti\gp_v^{2(j+1)}$, according as $v\in \gM_{\text{in}}$
or $\gM_{\text{rm}}$.
We first show that the right $G_{x\, \co_v/\gp_v^{j+1}}$-coset 
of $g$ contains an element of the form 
\begin{equation}
\left(\pmatrix 1 & 0\\ u_1 & t_1\endpmatrix,* \right)
.\end{equation}
If $v\in \gM_{\text{rm}}$, roots of $p(z)$ still
generates the unique unramified quadratic extension of $\ti k_v$  
and so an element $c_1+d_1\al_1$ with 
$c_1,d_1\in \ti\co_v$ is a unit if and only if 
$c_1$ or $d_1$ is a unit.  
We  choose 
$c_1=g_{122},d_1=g_{112}$.  Then either $c_1$ or $d_1$ is a unit.  
So there exist $c_2,d_2\in \co_v$ such that 
$(A_p(c_1,d_1),A_p(c_2,d_2))\in G_{x\, \co_v/\gp_v^{j+1}}$.  
Multiplying this element by $g$, we may assume that $g_{112}=0$.
  
If $v\in\gM_{\text{in}}$ then $k_v(x)=\ti k_v$ 
is unramified over $k_v$.  
Note that $p(z)$ is irreducible modulo $\gp_v$ by assumption.  
Therefore, $\bar x\in V^{\sst}_{\co_v/\gp_v}$.  Since 
$\co_v/\gp_v$ is a perfect field, one can use the same argument
as in Proposition (2.10)(2) \cite{kayu}, p. 323 to determine the stabilizer 
$G_{\bar x\,\co_v/\gp_v}$.  So the  proposition is true 
if $j=0$.  

By the previous step, we can multiply $g$ by an element of 
$G_{x\, \co_v/\gp_v}$ to assume that
$g_{112},g_{121}\in\ti\gp_v$.  Let $c_1=g_{122},d_1=g_{112}$.  
Then $c_1\in\ti\co_v^{\times}$ and $d_1\in\ti\gp_v$.  
So $c_1-d_1\al_1\equiv c_1-d_1\al_2\equiv c_1 \; (\ti\gp_v)$. 
Therefore, there exist $c_2,d_2\in \co_v$ such that 
$(A_p(c_1,d_1),A_p(c_2,d_2))\in G_{x\, \co_v/\gp_v^{j+1}}$.  
Multiplying $g$ by this element, we may assume $g_{112}=0$. 

In both cases,  $g_{112}=0$ and so 
$g_{111},g_{122}\in\ti\co_v^{\times}$.  Multiplying by
an element of the form 
$(A_p(c_1,0),A_p(c_2,d_2))\in G_{x\, \co_v/\gp_v^{j+1}}$,
we may further assume that $g_{111}=1$.  

This done, $(g_1,1)w_p$ is given by (\ref{module}) with $t_2,u_2$
replaced by $t_1^{\sig},u_1^{\sig}$.  Therefore, by the same argument,
$t_1=t_1^{\sig},\; u_1^{\sig}=u_1$, and the condition (\ref{urstab})
holds also.  After this, the argument is the same as in the case
$v\in\gM_{\text{sp}}$.
\end{proof}
\begin{defn}\label{easyddefn} 
Let $x$ denote any of the standard orbital representatives. In each
of the cases enumerated below, we define $\cD$ to be the set of
$y\in V_{\calo_v}$ such that $y\equiv
x\;(\gp_v^{n_1},\ti\gp_v^{n_2})$ if $v\notin\gM_{\text{sp}}$ or
$y\equiv x\;(\gp_v^n)$ if $v\in\gM_{\text{sp}}$, where $n_1$, $n_2$
and $n$ are as indicated.
\begin{enumerate}
\item[(1)] $n=2m_v+1$ if $x$ has type (sp) or (sp~ur),
\item[(2)] $n=2$ if $x$ has type (sp~rm) and
$v\notin\gM_{\text{dy}}$,
\item[(3)] $n_1=n_2=2m_v+1$ if $x$ has type (in) or (in~ur),
\item[(4)] $n_1=n_2=2$ if $x$ has type (in~rm) and
$v\notin\gM_{\text{dy}}$,
\item[(5)] $n_1=2m_v+1$, $n_2=4m_v+2$ if $x$ has type (rm) or
(rm~ur),
\item[(6)] $n_1=2$, $n_2=4$ if $x$ has type (rm~rm)* or (rm~rm~ur)
and $v\notin\gM_{\text{dy}}$.
\end{enumerate}
\end{defn}
\begin{prop} Let $x$ have one of the types enumerated in Definition
\ref{easyddefn} and $\cD$ be the corresponding set.
If $y\in {\mathcal D}$ then $k_v(y)=k_v(x)$.
\end{prop} 
\begin{proof}
By the choice of orbital representatives, $P(x)$ is a unit when
$k_v(x)/k_v$ is not ramified (including the case $k_v(x)=k_v$) 
and $(P(x))= \gp_v^{\del_{x,v}}$  when $k_v(x)/k_v$ 
is ramified. In cases (1), (3) and (5), 
we have $P(y)\equiv P(x) \; (\gp_v^{2m_v+1})$ and $P(x)$ is a unit. 
This congruence may be rewritten as 
$P(y)/P(x)\equiv 1\; (\gp_v^{2m_v+1})$.
By Hensel's lemma, a unit is a square if and only if it is a square
modulo $\gp_v^{2m_v+1}$.  Note that this is true whether 
or not $v$ is dyadic. So 
$P(y)/P(x)$ is a square. This implies that $k_v(y)=k_v(x)$ 
in all these cases.

In case (2), the condition implies that 
$P(y)\equiv P(x) \;(\gp_v^2)$
and hence that we have
$P(y)/P(x)\equiv 1 \; (\gp_v)$, which again implies
that $k_v(y)=k_v(x)$, since $v$ is not dyadic.  
The same argument works for case (4) because
$\ti\gp_v^{2}\cap\co_v=\gp^2_v$ in this case. 

In case (6), the assumption implies that 
$P(y)\equiv P(x) \; (\gp_v^2)$ and $\ord_{k_v}(P(x))=1$. 
So $P(y)P(x)^{-1}\equiv 1\; (\gp_v)$.  Since $v$ is not dyadic, 
$P(y)P(x)^{-1}\in (k_v^{\times})^2$.  Therefore, 
$k_v(y)=k_v(x)$.
\end{proof}
\begin{prop}
Let $x$ have one of the types enumerated in Definition
\ref{easyddefn} and $\cD$ be the corresponding set.
Then ${\mathcal D}\sub K_vx$.
\end{prop}
\begin{proof}
We first deal with (1)--(5) in Definition \ref{easyddefn}. For every
one of these orbit types, we showed in section \ref{sec-omega} that
there exists an omega set for $x$. Suppose that $y\in\cD$. Then
$k_v(y)=k_v(x)$ and so $y\in G_{k_v}x$. Since $y$ is also integral,
it follows from the first property of omega sets that
$y\in\Om_{x,v}x$. That is, $y=gx$ for some $g\in\Om_{x,v}$. We have
seen above that the congruence conditions on $y$ also imply that
$|P(y)|_v=|P(x)|_v$ and so $g\in\Om_{x,v}^1=K_v$. This proves the
claim for (1)--(5).

Consider case (6).  Let 
\begin{equation*}
G(\pi_v^2) = \{g\in K_v\mid g_1\equiv 1 \; (\ti\gp_v^4),\; 
g_2\equiv 1 \; (\gp_v^2)\}
.\end{equation*}
Then $G(\pi_v^2)$ fixes the set ${\mathcal D}$.  
Consider the usual coordinates
$y=(y_{ij})$ as in (\ref{case(2)coords}).
Since $y_{20}\equiv 1\;(\gp_v^2)$, 
$(1,a(y_{20}^{-1},1))\in G(\pi_v^2)$.  Applying this element, 
we may assume that $y_{20}=1$.  
Since $y_{10}\equiv 0\;(\gp_v^2)$, 
$(1,{}^tn(-y_{10}))\in G(\pi_v^2)$.  Applying this element, 
we may assume that $y_{10}=0$.  
Since $y_{11}\equiv 1\;(\ti\gp_v^4)$, 
$(a(y_{11}^{-1},1),1)\in G(\pi_v^2)$. Applying this element, 
we may assume that $y_{11}= 1$.  
Since we are assuming $v\notin \gMdy$, 
we chose $p$ so that $a_1=0$.  So $y_{12}\equiv 0\; (\gp_v^2)$. 
So by a similar argument, we may assume $y_{12}=0$.  
Also $y_{21}\equiv 0\; (\ti\gp_v^4)$ and 
$y_{22}\equiv -a_2\;(\gp_v^2)$. 
We may assume that $\ti k_v=k_v(\sqrt{\pi_v})$.  
Let $y_{21}= \pi_v^2(c+d\sqrt{\pi_v})$ where 
$c,d\in\co_v$.  Then 
$(n(-d\pi_v^2\sqrt{\pi_v}),n(-\pi_v^2c))\in G(\pi_v^2)$.  
Note that $\tr_{k_v/k_v}(-d\pi_v^2\sqrt{\pi_v})=0$.  
So applying this element, we may assume that $y_{21}=0$. 
Then $-y_{22} a_2^{-1}\equiv 1\;(\gp_v)$.  
So  $-y_{22} a_2^{-1} = t^2$ with $t\in\co_v^{\times}$
by Hensel's lemma.  Applying 
$(a(1,t^{-1}),a(t,1))$, we get $y_{22}=-a_2$.  
\end{proof}
\begin{prop}\label{volume1}
\begin{enumerate}
\item[(1)] If $x$ has type {\upshape (sp)} then 
$\vol(K_v x)=(1/2)(1+q_v^{-1})(1-q_v^{-2})^2$.
\item[(2)] If $x$ has type {\upshape (sp~ur)} then
$\vol(K_v x)=(1/2)(1-q_v^{-1})^3(1-q_v^{-2})$.
\item[(3)] If $x$ has type {\upshape (in)} or {\upshape (in~ur)} then 
$\vol(K_v x)=(1/2)(1-q_v^{-1})(1-q_v^{-4})$.
\item[(4)] If $x$ has type {\upshape (rm)} then 
$\vol(K_v x)=(1/2)(1-q_v^{-2})^2$.
\item[(5)] If $x$ has type {\upshape (rm~ur)} then
$\vol(K_v x)=(1/2)(1-q_v^{-1})^2(1-q_v^{-2})$.
\end{enumerate}
\end{prop}
\begin{proof}  Let $j=2m_v$ in all  cases.  Then 
$\vol(K_vx) = \vol({\mathcal D})\#\bigl( G_{\co_v/\gp_v^{j+1}}/
G_{\bar x\, \co_v/\gp_v^{j+1}}\bigr)$. By Proposition \ref{coinstab}, 
$G_{\bar x\, \co_v/\gp_v^{j+1}}=G_{x\, \co_v/\gp_v^{j+1}}$.  

Consider (1).   In this case, $G^{\circ}_{x\, \co_v/\gp_v^{j+1}}
\cong ((\co_v/\gp_v^{j+1})^{\times})^4$.  So its order is 
$((q_v-1)q_v^j)^4$.  Therefore, 
\begin{equation*}
\vol(K_vx)=q_v^{-8(j+1)} \cdot 
\frac{(q_v^2-q_v)^3(q_v^2-1)^3(q_v^j)^{12}}
{2(q_v-1)^4(q_v^j)^4}
=(1/2)(1+q_v^{-1})(1-q_v^{-2})^2\,.
\end{equation*}

Consider (2) and (3). Let $\gp_x\sub \co_{k_v(x)}$ 
be the  prime ideal.  
In both cases, $G^{\circ}_{x\, \co_v/\gp_v^{j+1}}
\cong ((\co_x/\gp_x^{j+1})^{\times})^2$, and 
$\#(\co_x/\gp_x)= q_v^2$.  So its order is 
$((q_v^2-1)q_v^{2j})^2$.    Therefore, in case (2),   
\begin{equation*}
\vol(K_vx)=q_v^{-8(j+1)} \cdot 
\frac{(q_v^2-q_v)^3(q_v^2-1)^3(q_v^j)^{12}}
{2(q_v^2-1)^2(q_v^{2j})^2}
=(1/2)(1-q_v^{-1})^3(1-q_v^{-2})\,.
\end{equation*}
In case (3), 
\begin{equation*}
\begin{aligned}
\vol(K_vx) & = q_v^{-8(j+1)} \cdot 
\frac {(q_v^4-q_v^2)(q_v^4-1)(q_v^{2j})^4(q_v^2-q_v)(q_v^2-1)(q_v^j)^4}
{2(q_v^2-1)^2(q_v^{2j})^2}  \\
& = (1/2)(1-q_v^{-1})(1-q_v^{-4})
.\end{aligned}
\end{equation*}

Consider (4).   In this case, $G^{\circ}_{x\, \co_v/\gp_v^{j+1}}
\cong ((\ti\co_v/\ti\gp_v^{2(j+1)})^{\times})^2$ and 
$\#(\ti\co_v/\ti\gp_v)=q_v$.    So its order is 
$((q_v-1)q_v^{2j+1})^2$.  Therefore, 
\begin{equation*}
\vol(K_vx)=q_v^{-8(j+1)} \cdot 
\frac{(q_v^2-q_v)^2(q_v^2-1)^2(q_v^{2j+1})^4 (q_v^j)^4} 
{2(q_v-1)^2(q_v^{2j+1})^2}
=(1/2)(1-q_v^{-2})^2\,.
\end{equation*}

Consider (5).   Let  $\gp_x\sub \co_{k_v(x)}$ 
be the  prime ideal. 
In this case, $G^{\circ}_{x\, \co_v/\gp_v^{j+1}}
\cong (\co_x/\gp_x^{2(j+1)})^{\times}$ and 
$\#(\co_x/\gp_x)=q_v^2$.    So its order is 
$(q_v^2-1)q_v^{2(2j+1)}$.  Therefore, 
\begin{equation*}
\vol(K_vx)=q_v^{-8(j+1)} \cdot 
\frac{(q_v^2-q_v)^2(q_v^2-1)^2(q_v^{2j+1})^4 (q_v^j)^4} 
{2(q_v^2-1)(q_v^{2(2j+1)})}
=(1/2)(1-q_v^{-1})^2(1-q_v^{-2})\,.
\end{equation*}
\end{proof}

\begin{prop}\label{volume2}
\begin{enumerate}
\item[(1)] Suppose that $x$ has type {\upshape (sp~rm)} and that
$v\notin\gM_{\text{{\upshape dy}}}$. Then 
\begin{equation*}
\vol(K_vx) = (1/2)q_v^{-1}(1-q_v^{-1})(1-q_v^{-2})^3\,.
\end{equation*}
\item[(2)] Suppose that $x$ has type {\upshape (in~rm)} and that
$v\notin\gM_{\text{{\upshape dy}}}$. Then 
\begin{equation*}
\vol(K_vx) = (1/2)q_v^{-1}(1-q_v^{-1})(1-q_v^{-2})(1-q_v^{-4})\,.
\end{equation*}
\end{enumerate}
\end{prop}
\begin{proof}  We prove that 
$[G_{\bar x\, \co_v/\gp_v^2}:G^{\circ}_{x\, \co_v/\gp_v^2}]=2q_v$. 
Let $p(z)=z^2+a_1z+a_2$ be the Eisenstein polynomial corresponding
to $x$.  Since $v$ is not dyadic, we assume $a_1=0$.  
 
Consider (1).  Elements of the form 
$(A_p(c_1,d_1),A_p(c_2,d_2),A_p(c_3,d_3))$ are in
$G^{\circ}_{x\, \co_v/\gp_v^2}$ if and only if  
$c_1^2+a_2d_1^2,c_2^2+a_2d_2^2\in (\co_v/\gp_v^2)^{\times}$
(see Lemma \ref{stabreduction(1)}).   
So the order of $G^{\circ}_{x\, \co_v/\gp_v^2}$ is $(q_v-1)q_v^3$.  
Since $a_2\in\gp_v$, this is equivalent to 
$c_1,c_2\in (\co_v/\gp_v^2)^{\times}$.
Suppose $g=(g_1,g_2,g_3)\in G_{\bar x\, \co_v/\gp_v^2}$.  
Since $F_x(v)$ reduces to $v_1^2$ modulo $\gp_v$, 
$g_{321}\in \gp_v/\gp_v^2$.  Let $x=(x_1,x_2)$.  
Then  $(g_1,g_2,1)x_2$ is a unit scalar multiple of $x_2$ modulo
$\gp_v$.  
By computation, 
\begin{equation*}
(g_1,g_2,1)x_2 \equiv \pmatrix 
g_{111}g_{211} & g_{111}g_{221}\\ 
g_{121}g_{212} & g_{121}g_{222}\endpmatrix\;(\gp_v)\,.
\end{equation*}
So $g_{111},g_{211}$ are units.  
By a similar argument as in the proof of Proposition \ref{coinstab},
$g_{121},g_{221}\in \gp_v/\gp_v^2$.  This implies that 
$g_{122},g_{222}$ are units.  Since 
\begin{equation*}
A_p(c_i,d_i) g_i = \pmatrix * & c_ig_{i12}-d_ig_{i22}\\
* & * \endpmatrix
,\end{equation*}
for $i=1,2$, the right coset $G^{\circ}_{x\, \co_v/\gp_v^2} g$ contains 
an element $g=(g_1,g_2,g_3)$ where $g_1,g_2$ 
are in the form (\ref{cosetrep}).  Moreover, 
it is easy to see that 
$t_1,t_2,u_1,u_2$ are determined by the coset 
$G^{\circ}_{x\, \co_v/\gp_v^2} g$.  We use Lemma \ref{observation}
to determine the possibilities for $g$. 
By computation, 
\begin{equation*}
(g_1,g_2,1)x = 
\left( \pmatrix 0 & t_2\\ t_1 & t_1u_2+t_2u_1\endpmatrix,
\pmatrix 1 & u_2\\ u_1 & u_1u_2-a_2t_1t_2\endpmatrix \right)
.\end{equation*}
The condition  $\xspan((g_1,g_2,1)x)=\xspan(x)$ is 
equivalent to the condition that 
$t_1=t_2,u_1=u_2,t_1u_2+t_2u_1=0$, and $u_1u_2-a_2t_1t_2 = -a_2$. 
Since $v$ is not dyadic, $u_1=u_2=0$ and $a_2(t_1^2-1)=0$.  
Therefore, $t_1\equiv\pm 1 \;(\gp_v)$.  So there are $2q_v$ 
possibilities for $t_1$ modulo $\gp_v^2$. This proves that  
$[G_{\bar x\, \co_v/\gp_v^2}:G^{\circ}_{x\, \co_v/\gp_v^2}]=2q_v$.

In case (2), by a similar argument as in case (1), 
we can assume that $g_1$ is in the form  (\ref{cosetrep})
with $t_1$ unit and $u_1\in\ti\co_v/\ti\gp_v^2$.  
Note that the order of $G^{\circ}_{x\, \co_v/\gp_v^2}$ 
is $(q_v^2-1)q_v^6$ in this case.  
By the same 
consideration, we get $t_1=t_1^{\sig},u_1=u_1^{\sig}$, and
the rest of the conditions turn out to be the same.  
Since $\ti\gp_v^2\cap\co_v = \gp_v^2$, there are $2q_v$ possibilities
for $t_1$.  So 
$[G_{\bar x\, \co_v/\gp_v^2}:G^{\circ}_{x\, \co_v/\gp_v^2}]=2q_v$
in this case also.  

Since the volume of ${\mathcal D}$ in Definition \ref{easyddefn}(2)
is $q_v^{-16}$, if $v\in \gM_{\text{sp}}$, 
\begin{equation*}
\vol(K_v x)=q_v^{-16}\cdot  
\frac {(q_v^2-q_v)^3(q_v^2-1)^3q_v^{12}}
{2q_v((q_v-1)q_v^3)^2}
=(1/2)q_v^{-1}(1-q_v^{-1})(1-q_v^{-2})^3
\end{equation*}
and if $v\in \gM_{\text{in}}$, 
\begin{equation*}
\begin{aligned}
\vol(K_vx)&=q_v^{-16}\cdot  
\frac{(q_v^4-q_v^2)(q_v^4-1)q_v^8 (q_v^2-q_v)(q_v^2-1)q_v^4}
{2q_v(q_v^2-1)q_v^6} \\
&=(1/2)q_v^{-1}(1-q_v^{-1})(1-q_v^{-2})(1-q_v^{-4})\,.  
\end{aligned}
\end{equation*}
\end{proof}  

Next we assume that $v$ is dyadic and that $x$ has type (sp~rm) or
(in~rm). As explained in section \ref{final}, we shall compute the
sum of $\vol(K_vx)$ over the equivalence class of $x$ under the
relation $\asymp$. For these orbits, this amounts to summing
$\vol(K_vx)$ over all $x$ having a given value of $\del_{x,v}$. Let
$x$ be coordinatized as in (\ref{case(1)coords}) or
(\ref{case(2)coords}).
If $v\in \gM_{\text{sp}}$ then $F_x(v) = a_0(x)v_1^2+
a_1(x)v_1v_2+a_2(x)v_2^2$ where
\begin{equation}\label{spaform}
\begin{aligned}
a_0(x) & = x_{112}x_{121}-x_{111}x_{122}, \\
a_1(x) & = x_{112}x_{221}+x_{121}x_{212}-x_{111}x_{222}-x_{122}x_{211}, \\
a_2(x) & = x_{212}x_{221}-x_{211}x_{222}\,
.\end{aligned}
\end{equation}
If $v\notin \gM_{\text{sp}}$ then $F_x(v) = a_0(x)v_1^2+
a_1(x)v_1v_2+a_2(x)v_2^2$ where 
\begin{equation}  \label{aform}
\begin{aligned}
a_0(x) & = \n_{\ti k_v/k_v}(x_{11})-x_{10}x_{12}, \\
a_1(x) & = \tr_{\ti k_v/k_v}(x_{11}x_{21}^{\sig})
-x_{10}x_{22}-x_{12}x_{20}, \\
a_2(x) & = \n_{\ti k_v/k_v}(x_{21})-x_{20}x_{22}\,
.\end{aligned}
\end{equation}
\begin{defn}\label{ddefn1}
\begin{enumerate}
\item[(1)] If $v\in\gM_{\text{{\upshape sp}}}$ then
we define ${\mathcal D}_{\ell}$ for $\ell=1,\dots,m_v+1$
to be the set of  $x$ such that 
\begin{equation*}
\begin{aligned}
{} & x_{112},x_{121},x_{211}\in \co_v^{\times},\; 
x_{122},x_{212},x_{221}\in\gp_v, \\
& \ord_{k_v}(x_{222})=1,\; \ord_{k_v}(a_1(x))
\begin{cases} = \ell &\text{ if } \ell\leq m_v, \\ 
\geq m_v+1 &\text{ if } \ell=m_v+1.\end{cases}
\end{aligned}
\end{equation*}
\item[(2)] If $v\in\gM_{\text{{\upshape in}}}$  then
we define ${\mathcal D}_{\ell}$ for $\ell=1,\dots,m_v+1$
to be the set of $x$ such that 
\begin{equation*}
\begin{aligned}
{} & x_{11}\in \ti\co_v^{\times},\; 
x_{20}\in \co_v^{\times},\;  
x_{12}\in\gp_v, x_{21}\in\ti\gp_v, \\ 
& \ord_{k_v}(x_{22})=1,\;  \ord_{k_v}(a_1(x))
\begin{cases} = \ell &\text{ if } \ell\leq m_v, \\ 
\geq m_v+1 &\text{ if } \ell=m_v+1.\end{cases}
\end{aligned}
\end{equation*}
\end{enumerate}
\end{defn}
\begin{prop}\label{dyadDvol}
\begin{enumerate}
\item[(1)] If $v\in\gM_{\text{{\upshape sp}}}$ then
\begin{equation*}
\vol({\mathcal D}_{\ell}) = \begin{cases}
q_v^{-3}(1-q_v^{-1})^5q_v^{-\ell} &\text{ if } \ell\leq m_v, \\
q_v^{-3}(1-q_v^{-1})^4q_v^{-(m_v+1)} &\text{ if } \ell=m_v+1\,
.\end{cases}
\end{equation*}
\item[(2)] If $v\in\gM_{\text{{\upshape in}}}$ then
\begin{equation*}
\vol({\mathcal D}_{\ell}) = \begin{cases}
q_v^{-3}(1-q_v^{-1})^3(1-q_v^{-2})q_v^{-\ell} &\text{ if } 
\ell\leq m_v, \\
q_v^{-3}(1-q_v^{-1})^2(1-q_v^{-2})q_v^{-(m_v+1)} &\text{ if } 
\ell=m_v+1\,
.\end{cases}
\end{equation*}
\end{enumerate}
\end{prop}
\begin{proof}
Suppose $v\in \gM_{\text{sp}}$.  It is easy to see that
if $x\in V_{\co_v}$ satisfies all the conditions
in Definition \ref{ddefn1}(1) except for the last condition, 
then  $a_0(x)\in\co_v^{\times}$ and $\ord_{k_v}(a_2(x))=1$.  
Suppose 
\begin{equation*}
x_{112}x_{221}+x_{121}x_{212}-x_{122}x_{211}
= e_1\pi_v + \dots + e_{\ell}\pi_v^{\ell}+\cdots 
,\end{equation*}
where $e_j$ is a unit or zero for all $j$ 
(see I.4, Corollary 2 in \cite{weilc}, p. 14).  
Let $x_{222}=\pi_v x_{222}'$ where $x_{222}'\in\co_v^{\times}$.  
Then by (\ref{spaform}), 
$a_1(x)$ satisfies the condition in Definition \ref{ddefn1}(1)
if and only if $x_{111}$ has an expansion of the form 
\begin{equation*}
x_{111} = (x_{222}')^{-1}(e_1 +\dots + e_{\ell-1}\pi_v^{\ell-2}
+ f \pi_v^{\ell-1}+\cdots)
\end{equation*}
such that $f\not\equiv e_{\ell}\;(\gp_v)$ if $\ell\leq m_v$ 
(and no condition on $f$ if $\ell=m_v+1$).   
If $\ell\leq m_v$, the volume of the set of such $x_{111}$ is 
$(q_v-1)q_v^{-\ell}=(1-q_v^{-1})q_v^{-\ell+1}$ if $\ell\leq m_v$ and 
$q_v^{-m_v}$ if $\ell=m_v+1$.  
The volumes of the sets of $x_{112},x_{121},x_{211}$
are $1-q_v^{-1}$, the volumes of the sets of 
$x_{122},x_{212},x_{221}$ are 
$q_v^{-1}$, and the volume of the set of $x_{222}$ is 
$q_v^{-1}(1-q_v^{-1})$.  
Therefore, 
\begin{equation*}
\vol({\mathcal D}_{\ell}) = 
(1-q_v^{-1})^3q_v^{-3}q_v^{-1}(1-q_v^{-1}) \times 
\begin{cases} (1-q_v^{-1})q_v^{-\ell+1} &\text{ if } \ell\leq m_v, \\
q_v^{-m_v} &\text{ if } \ell=m_v+1.\end{cases}
\end{equation*}
Simplifying, we get (1).  

If $v\in \gM_{\text{in}}$, similar
considerations apply to $x_{10}$ as to
$x_{111}$ in the case $v\in\gM_{\text{sp}}$.  
The volumes of the sets of $x_{11},x_{20}$ are $1-q_v^{-2}$ 
and $1-q_v^{-1}$, respectively,  
the volumes of the sets of $x_{12},x_{21}$ are 
$q_v^{-2}$ and $q_v^{-1}$, respectively,
and the volume of the set of $x_{22}$ is $q_v^{-1}(1-q_v^{-1})$.
Therefore, 
\begin{equation*}
\vol({\mathcal D}_{\ell}) = 
(1-q_v^{-2})(1-q_v^{-1})q_v^{-2}q_v^{-1}q_v^{-1}(1-q_v^{-1}) \times 
\begin{cases} (1-q_v^{-1})q_v^{-\ell+1} &\text{ if } \ell\leq m_v, \\
q_v^{-m_v} &\text{ if } \ell=m_v+1.\end{cases}
\end{equation*}
Simplifying, we get (2).  
\end{proof} 

We employ coordinates on $G$ as in (\ref{groupcoords}). If
$v\in\gM_{\text{sp}}$ then we define
\begin{equation}\label{hstab1}
\begin{aligned}
\ & H_1=\{g\in K_v\mid g_{121},g_{221},g_{321}\in\gp_v\}, \\
\ &\overline{H}_1=
\{g\in(\gl(2)_{\calo_v/\gp_v})^3\mid
g_{121}=g_{221}=g_{321}=0\}
\end{aligned}
\end{equation}
and if $v\in\gM_{\text{in}}$ then we define
\begin{equation}\label{hstab2}
\begin{aligned}
\ &H_2=\{g\in K_v\mid g_{121}\in\ti\gp_v, g_{221}\in\gp_v\}, \\
\ &\overline{H}_2=\{g\in\gl(2)_{\ti\calo_v/\ti\gp_v}\times
\gl(2)_{\calo_v/\gp_v}\mid g_{121}=0, g_{221}=0\}\,.
\end{aligned}
\end{equation}
Note that
\begin{equation}\label{hstab1ord}
\#(K_v/H_1)=\#(G_{\calo_v/\gp_v}/\overline{H}_1)=
\frac{(q_v^2-1)^3(q_v^2-q_v)^3}
{(q_v-1)^6q_v^3}=(q_v+1)^3
\end{equation}
and
\begin{equation}\label{hstab2ord}
\begin{aligned}
\#(K_v/H_2)=\#(G_{\calo_v/\gp_v}/\overline{H}_2)&=
\frac{(q_v^4-1)(q_v^4-q_v^2)(q_v^2-1)(q_v^2-q_v)}
{(q_v^2-1)^2q_v^2(q_v-1)^2q_v} \\
&=(q_v^2+1)(q_v+1)\,.
\end{aligned}
\end{equation}

\begin{lem}\label{spdvol}
\begin{enumerate}
\item[(1)]  If $v\in\gM_{\text{{\upshape sp}}}$ and $g\in H_1$ then
$g{\mathcal D}_{\ell}={\mathcal D}_{\ell}$.
\item[(2)] If $v\notin\gM_{\text{{\upshape sp}}}$ and $g\in H_2$ then
$g{\mathcal D}_{\ell}={\mathcal D}_{\ell}$.
\end{enumerate}
\end{lem}
\begin{proof} A simple, direct calculation shows that elements of
$H_1$ and $H_2$ preserve all the conditions for membership in
$\cD_{\ell}$, with the possible exception of the last. Consider (1)
and suppose that $y=gx$ with $g\in H_1$ and $x\in\cD_{\ell}$. Then
$F_y(v)=\det(g_1)\det(g_2)F_x(vg_3)$ and so $a_1(y)$ is equal to
\begin{equation*}
\det(g_1)\det(g_2)\big(a_1(x)[g_{311}g_{322}+g_{321}g_{312}]
+2a_0(x)g_{311}g_{321}+2a_2(x)g_{312}g_{322}\big)\,.
\end{equation*}
We must investigate the order of this number. Since $\det(g_1)$ and
$\det(g_2)$ are units, they may be ignored. Both $g_{311}$ and
$g_{322}$ are units and, since $g_{321}\in\gp_v$,
$g_{311}g_{322}+g_{321}g_{312}$ is also a unit. Thus the order of
$a_1(x)[g_{311}g_{322}+g_{321}g_{312}]$ equals $\ell$ if $\ell\leq
m_v$ and is greater than or equal to $m_v+1$ if $\ell=m_v+1$. Also,
$\ord_{k_v}(2a_0(x)g_{311}g_{321})\geq m_v+1$ and
$\ord_{k_v}(2a_2(x)g_{312}g_{322})\geq m_v+1$, in the first case
because $g_{321}\in\gp_v$ and in the second because
$a_2(x)\in\gp_v$. It follows that the order of $a_1(y)$ is $\ell$
if $\ell\leq m_v$ and is greater than or equal to $m_v+1$ if
$\ell=m_v+1$. Thus $y\in\cD_{\ell}$, as required. Similar arguments
apply in case (2).
\end{proof}
\begin{prop}\label{volume3}
\begin{enumerate}
\item[(1)]  Suppose $v\in\gM_{\text{{\upshape sp}}}$ is dyadic. Then 
\begin{equation*}
\begin{aligned}
\sum_{2\leq \del_{x,v}= 2\ell\leq 2m_v}\vol(K_vx) & = 
(1-q_v^{-1})^2(1-q_v^{-2})^3q_v^{-\ell}, \\
\sum_{\del_{x,v}= 2m_v+1}\vol(K_vx) & = 
(1-q_v^{-1})(1-q_v^{-2})^3q_v^{-(m_v+1)}
.\end{aligned}
\end{equation*}
\item[(2)]  Suppose $v\in\gM_{\text{{\upshape in}}}$ is dyadic. Then 
\begin{equation*}
\begin{aligned}
\sum_{\del_{x,v}= 2\ell\leq 2m_v}\vol(K_vx) & = 
(1-q_v^{-1})^2(1-q_v^{-2})(1-q_v^{-4})q_v^{-\ell}, \\
\sum_{\del_{x,v}= 2m_v+1}\vol(K_vx) & = 
(1-q_v^{-1})(1-q_v^{-2})(1-q_v^{-4})q_v^{-(m_v+1)}
.\end{aligned}
\end{equation*}
\end{enumerate}
\end{prop}  
\begin{proof}
Consider (1). We first show that $\cD_{\ell}\subseteq\cup_xK_vx$,
where $x$ runs through the standard orbital representatives for the
orbits of type (sp~rm) whose corresponding fields have a fixed
discriminant. As in the statement, the discriminant is related to
$\ell$ by $\ell=\lfloor\tfrac12(\del_{x,v}+1)\rfloor$. By
construction, if $y\in\cD_{\ell}$ then $|P(y)|_v=q_v^{-\del_{x,v}}$
and so if $x$ is the standard representative for the orbit
corresponding to $k_v(y)$ then $|P(y)|_v=|P(x)|_v$. Since $y\in
G_{k_v}x\cap V_{\calo_v}$, the theory of omega sets implies that
there exists $g\in\Om_{x,v}$ such that $y=gx$. But then
$|\chi(g)|_v=1$ and so $g\in\Om_{x,v}^1=K_v$. That is, $y\in
K_vx$.

Each one of the standard orbital representatives itself lies in
$\cD_{\ell}$ and so the $K_v$-translates of $\cD_{\ell}$ cover
$\cup_xK_vx$. In order to find the volume of $\cup_xK_vx$ it
remains to determine the number of distinct $K_v$-translates of
$\cD_{\ell}$. Suppose that $g\in K_v$ and
$g\cD_{\ell}\cap\cD_{\ell}\neq\emptyset$; say $x,y\in\cD_{\ell}$
and $y=gx$. We shall show that $g\in H_1$.
Since $F_x(v),F_y(v)$
both reduce to unit scalar multiples of $v_1^2$ modulo $\gp_v$, 
the assumption implies that $g_{321}\in\gp_v$.  Since 
$(1,1,g_3)\in H_1$, we may assume $g_3=1$.  
Let $x=(x_{ijk})$ and $\bar x=(\bar x_{ijk})$ be 
the reduction of $x$ modulo $\gp_v$.  
We define $\bar g_{ijk}$, $\bar y$ similarly.   
Then 
\begin{equation*}
(g_1,g_2,1) \bar x  = \bar y = 
\left(*,  \bar x_{211} \pmatrix 
\bar g_{111}\bar g_{211} & \bar g_{111}\bar g_{221}\\
\bar g_{121}\bar g_{212} & \bar g_{121}\bar g_{222}\endpmatrix \right)
\,.\end{equation*}  
Since $\bar g_{212}\neq0$ or $\bar g_{222}\neq0$, $\bar g_{121}=0$.
This implies that $\bar g_{111}\neq0$ and therefore $\bar g_{221}=0$. 
Thus $g\in H_1$, as claimed.

It follows from this and Lemma \ref{spdvol} that two
$K_v$-translates of $\cD_{\ell}$ are either disjoint or equal and
that the number of them is $\#(K_v/H_1)=(q_v+1)^3$ by
(\ref{hstab1ord}). Using Proposition \ref{dyadDvol}, the result
follows. Similar arguments apply to prove the formula in case (2).
\end{proof} 

\begin{prop}\label{volume4}
If $v$ is not dyadic and $x$ has type {\upshape (rm~rm~ur)} or 
{\upshape (rm~rm)*} then 
\begin{equation*}
\vol(K_v x) = (1/2)q_v^{-1}(1-q_v^{-1})
(1-q_v^{-2})^2\,.
\end{equation*}
\end{prop} 
\begin{proof}  Consider ${\mathcal D}$ in Definition \ref{easyddefn}(6).
Obviously $\vol({\mathcal D})=q_v^{-16}$.  
In both cases, $G^{\circ}_{x\, \co_v/\gp_v^2}$ consists of elements 
of the form $(A_p(c,d),*)$ where $c\in (\ti\co_v/\ti\gp_v^4)^{\times},\;
d\in \ti\co_v/\ti\gp_v^4$.  
So $\# G^{\circ}_{x\, \co_v/\gp_v^2}= (q_v-1)q_v^7$. 

We show that 
$[G_{\bar x\, \co_v/\gp_v^2}:G^{\circ}_{x\, \co_v/\gp_v^2}]=2q_v$.
Suppose $g\in G_{\bar x\, \co_v/\gp_v^2}$.  We consider the
coordinates (\ref{groupcoords}) again.
Since $F_x(v)$ is a unit scalar multiple of 
$v_1^2$ modulo $\gp_v$, $g_{221}\equiv 0\;(\gp_v)$.  
Then since $g_1 \left(\begin{smallmatrix} 1 & 0\\ 
0 & 0\end{smallmatrix} \right){}^tg_1^{\sig}$ is a unit 
scalar multiple of 
$\left(\begin{smallmatrix} 1 & 0\\ 0 & 0\end{smallmatrix} \right)$, 
$g_{121}\equiv 0\;(\ti\gp_v)$.  
So $g_{111},g_{122},g_{211},g_{222}$ are units.  
This implies that the right $G^{\circ}_{x\, \co_v/\gp_v^2}$-coset of
$g$ contains an element of the form 
$\left(\left(\begin{smallmatrix} 1 & 0\\ 
u & t\end{smallmatrix} \right),*\right)$.  
We use Lemma \ref{observation} again. 
By computation, 
\begin{equation*}
\left(\pmatrix 1 & 0\\  u & t\endpmatrix,1\right)x
= \left(\pmatrix 0 & t^{\sig}\\ t 
& \tr_{\ti k_v/k_v}(tu^{\sig})\endpmatrix,
\pmatrix 1 & u^{\sig}\\ 
u & \n_{\ti k_v/k_v}(u)-a_2\n_{\ti k_v/k_v}(t)\endpmatrix\right)\,
.\end{equation*}
So $t=t^{\sig},u=u^{\sig},2tu=0$, and 
$a_2(t^2-1)-u^2=0$.  Since $v\notin \gM_{\text{dy}}$, 
$u=0$, and $t \equiv \pm 1\; (\gp_v)$.
Since $t\in (\co_v/\gp_v^2)^{\times}$,
there are $2q_v$ possibilities for $t$.  
Therefore, in both cases, 
\begin{equation*}
\vol(K_vx)=q_v^{-16}\cdot \frac {(q_v^2-1)^2(q_v^2-q_v)^2 q_v^{16}}
{2(q_v-1)q_v^8}
=(1/2)q_v^{-1}(1-q_v^{-1})(1-q_v^{-2})^2\,.
\end{equation*}
\end{proof}

\section{Computation of $b_{x,v}$ at the infinite places}%
\label{sec-infinite}

In this section we shall compute the constants $b_{x,v}$ when $v$ is an
infinite place of $k$. By Proposition \ref{bxindep}, we know that we are free
to make use of any orbital representative to carry out this calculation
and we shall take advantage of this freedom to simplify our task
as much as possible. In the proof of each of the four propositions
in this section we shall have to calculate the $8\times 8$ Jacobian
determinant associated with the map $g\mapsto gx$ in some
coordinate system. Each of these calculations was carried out using
the \textsc{Maple} computer algebra package \cite{maple}.

Before we begin, we recall a fact about
the Haar measure on $\gl(2)_{F}$ where $F=\mathbb{R}$ or
$\mathbb{C}$. As is usual, we shall take Lebesgue measure to be the
standard measure on the real numbers and twice Lebesgue measure to
be the standard measure on the complex numbers.
If $g=\left(\begin{smallmatrix}a_{11}&a_{12}\\
a_{21}&a_{22}\end{smallmatrix}\right)$ then
$d\mu(g)=da_{11}\,da_{12}\,da_{21}\,da_{22}/|\det(g)|_F^2$ defines a
Haar measure on $\gl(2)_{F}$ and it is well known that if $dg$
denotes our standard choice of Haar measure on this group then
$dg=p_F\,d\mu(g)$ where $p_F=1/\pi$ if $F=\mathbb{R}$ and
$p_F=1/2\pi$ if $F=\mathbb{C}$.  For example if $F=\mathbb{R}$, 
using the usual paramatrization
of each factor of the Iwasawa decomposition, 
all one has to do is to find the Jacobian of the map 
$\orth(2)\times T_{\R}\times N_{\R}\to G_{\R}$ and the degree of
this covering and this can be done easily.  
When $F=k_v$ we shall use the
notation $p_v$ in place of $p_{F}$. We shall refer to $d\mu(g)$ as
the \emph{coordinate measure} on $\gl(2)_F$. 
\begin{prop} \label{bxinf-(sp)}
Let $v\in\mathfrak{M}_{\infty}$ and
suppose that $x$ is a representative for the orbit of type 
\upshape{(sp)}.
Then
\begin{equation*}
b_{x,v}=\begin{cases} \tfrac{2}{\pi^3} &\text{\quad if
$v\in\mathfrak{M}_{\mathbb{R}}$,} \\
\tfrac{1}{4\pi^3} &\text{\quad if $v\in\mathfrak{M}_{\mathbb{C}}$.}
\end{cases}
\end{equation*}
\end{prop}
\begin{proof}
As noted above, it suffices to compute $b_{x,v}$ where $x=w$, the
element introduced in (\ref{wdefn}). Let $d\mu_v$ be the product of
the coordinate measures on each of the three factors of $G_{k_v}$
so that $dg_v=p_v^3\,d\mu_v(g)$.

Let $\mathcal{B}$ denote the set of $(g_1,g_2,g_3)\in G_{k_v}$ such
that $g_1$ and $g_2$ lie in the big Bruhat cell. Then $\mathcal{B}$
is dense in $G_{k_v}$ and it is on this set that we shall carry out
the comparison of measures. Any element $g$ of $\mathcal{B}$ may be
written as
\begin{equation}\label{eq:infty_coords1}
g=({}^tn(u_1)n(y_1)a(t_{11},t_{12}),
{}^tn(u_2)n(y_2)a(t_{21},t_{22}),g_3a(t_{11}^{-1}t_{21}^{-1},
t_{12}^{-1}t_{22}^{-1}))
\end{equation}
with $g_3=\left(\begin{smallmatrix}a_{11}&a_{12}\\a_{21}&a_{22}
\end{smallmatrix}\right)$ and when $g$ is written in this form,
$u_i$, $y_i$ and $a_{ij}$ for $i,j=1,2$ may be regarded as
coordinates on $G_{k_v}/G^{\circ}_{x\,k_v}$. 
Note that the map 
\begin{equation*}
k_v^2\times (k_v^{\times})^2\ni 
(u,y,t_1,t_2)\to {}^tn(u)n(y)a(t_1,t_2)\in \gl(2)_{k_v}
\end{equation*}
is injective.  With respect to these
coordinates, the Jacobian determinant of the map 
$g\mapsto gx$ is found to be $|\det(g_3)|_v^2$. 
Note that this map is a double cover because 
$[G_{x\, k_v}:G^{\circ}_{x\, k_v}]=2$.  
Since  $P(gx)=\chi(g)P(x)=\det(g_3)^2P(x)$ and
$P(x)=1$, it follows that the pullback
of the measure $dy/|P(y)|_v^2$ to $G_{k_v}/G^{\circ}_{x\,k_v}$ is
$dg_{x,v}'=\tfrac12\,du_1\,dy_1\,du_2\,dy_2\,d\mu_{3,v}(g_3)$,
where $d\mu_{3,v}(g_3)$ denotes the coordinate measure on
the third factor. (The measure has been divided by $2$ because the
map $g\mapsto gx$ is a double cover and, in (\ref{lint}), the
measure $dg_{x,v}'$ was defined via an integral over $V_{k_v}$.)
In (\ref{eq:infty_coords1}), $t_{ij}$ may be regarded as
coordinates on $G_{x\,k_v}^{\circ}$ and as such
$dg_{x,v}''=\prod_{i,j=1,2}dt_{ij}/|t_{ij}|_v$. It follows from the
definition and the remarks above that $d\mu_v(g)=p_v^{-3}b_{x,v}\,
dg_{x,v}'\,dg_{x,v}''$ and so all that remains is for us to
determine the relationship between the coordinate measure and the
measure $du\,dy\,dt_1\,dt_2/|t_1t_2|_v$ on the big cell inside
$\gl(2)$ when it is coordinatized as ${}^tn(u)n(y)a(t_1,t_2)$. A
simple Jacobian calculation shows that in fact these two measures
are equal. This gives $\tfrac12p_v^{-3}b_{x,v}=1$ or equivalently
$b_{x,v}=2p_v^3$. Using the remarks before the proof, this gives
the stated values. \end{proof}
\begin{prop}\label{bxinf-(rm)}
Let $v\in\mathfrak{M}_{\mathbb{R}}$ and suppose that $x$ is a
representative for the orbit of type \upshape{(rm)}. Then
\begin{equation*}
b_{x,v}=\tfrac{1}{\pi^2}\,.
\end{equation*}
\end{prop}
\begin{proof}
We shall again use $x=w$ as the orbital representative. 
Let $\ti v$ be the complex place of $\ti k$ which divides
$v$.   Let  $d\mu_v$ be the product of 
the coordinate measures on the two
factors of $G_{k_v}$ so that $dg_v=p_{\tilde{v}}p_v\,d\mu_v(g)$.

We let $\mathcal{B}$ denote the set of $(g_1,g_2)\in G_{k_v}$ such
that $g_1$ lies in the big Bruhat cell. An element of $\mathcal{B}$
may be written as
\begin{equation}\label{eq:infty_coords2}
g=({}^tn(u)n(y)a(t_1,t_2),g_2a(
\operatorname{N}_{\kt_{\tilde{v}}/k_v}(t_1)^{-1},
\operatorname{N}_{\kt_{\tilde{v}}/k_v}(t_2)^{-1}))
\end{equation}
with $g_2=\left(\begin{smallmatrix}a_{11}&a_{12}\\a_{21}&a_{22}
\end{smallmatrix}\right)$ and when $g$ is written in this form,
$u$, $y$ and $a_{ij}$, $i,j=1,2$ may be regarded as coordinates on
$G_{k_v}/G^{\circ}_{x\,k_v}$. With respect to the real coordinates
$\re(u)$, $\im(u)$, $\re(y)$, $\im(y)$, $a_{ij}$, the Jacobian
determinant of the map $g\mapsto gx$ is found to be
$4|\det(g_2)|_v^2$. Since the map is a double cover, $P(g\cdot
x)=\det(g_2)^2$ and the canonical measures $du$ and $dy$ are
$du=2\,d\re(u)\,d\im(u)$ and $dy=2\,d\re(y)\,d\im(y)$, it follows
that the pullback of the measure $dy/|P(y)|_v^2$ to
$G_{k_v}/G^{\circ}_{x\,k_v}$ is
$dg'_{x,v}=\tfrac12\,du\,dy\,d\mu_{2,v}(g_2)$, where
$d\mu_{2,v}(g_2)$ denotes the coordinate measure on the second
factor. In (\ref{eq:infty_coords2}), $t_j$ may be regarded as
coordinates on $G^{\circ}_{x\,k_v}$ and as such
$dg''_{x,v}=\prod_{j=1,2}dt_j/|t_j|_{\tilde{v}}$. It follows from
the definition and the remarks above that
$d\mu_v(g)=p_{\tilde{v}}^{-1}p_{v}^{-1}
b_{x,v}\,dg_{x,v}'\,dg_{x,v}''$ and since the coordinate measure
restricted to the big cell inside $\gl(2)_{\kt_{\tilde{v}}}$ is
$du\,dy\,dt_1\,dt_2/|t_1t_2|_{\tilde{v}}$, we have
$b_{x,v}=2p_{\tilde{v}}p_{v}$. \end{proof}
\begin{prop}\label{bxinf-(sp_rm)}
Suppose that $v\in\mathfrak{M}_{\mathbb{R}}$ and that $x$ is a
representative for the orbit of type \upshape{(sp rm)}. Then
\begin{equation*}
b_{x,v}=\tfrac{2}{\pi^3}.
\end{equation*}
\end{prop}
\begin{proof}
We shall use $x=w_p$, where $p(z)=z^2+1$, as the orbital
representative for this orbit. 
The roots  of $p$ are $\pm \sqrt{-1}$.
Let $A_p(c,d)$ be as in (\ref{apdefn}).  
As we discussed in Lemma \ref{stabexpl(1)}, 
any $g\in G^{\circ}_{x\,k_v}$ has the form
$g=(A_p(c_1,d_1),A_p(c_2,d_2),A_p(c_3,d_3))$ where
\begin{equation*}
c_3 =  \frac{c_1c_2-d_1d_2}{(c_1^2+d_1^2)(c_2^2+d_2^2)},\; 
d_3 = \frac{-(c_1d_2+ c_2d_1)}{(c_1^2+d_1^2)(c_2^2+d_2^2)}
.\end{equation*}
The isomorphism $\theta$ from $G^{\circ}_{x k_v}$ to
$H_{x k_v}=(\mathbb{C}^{\times})^2$ may be taken as
\begin{equation*}
\theta(A_p(c_1,d_1),A_p(c_2,d_2),A_p(c_3,d_3))
=(c_1+\sqrt{-1}d_1,c_2+\sqrt{-1}d_2)
.\end{equation*}
Recalling that the canonical
measure on $\mathbb{C}^{\times}$ is $dz/|z|_{\mathbb{C}}$ where
$dz$ is twice Lebesgue measure, we see that
\begin{equation*}
dg''_{x,v}=\frac{4\,dc_1\,dd_1\,dc_2\,dd_2}
{(c_1^2+d_1^2)(c_2^2+d_2^2)}\,.
\end{equation*}

For any $r$ and any $\ka\in\operatorname{SO}(2)$ we may find $c$
and $d$ such that $a(r,r)\ka=A_p(c,d)$ and it follows from
this and the Iwasawa decomposition that any
$g\in\gl(2)_{\mathbb{R}}$ may be expressed as
$g=n(u)a(1,t)A_p(c,d)$. Thus any $g\in G_{k_v}$ may be
expressed as
\begin{equation*}
g=(n(u_1)a(1,t_1)A_p(c_1,d_1),
n(u_2)a(1,t_2)A_p(c_2,d_2),g_3A_p(c_3,d_3))
\end{equation*}
where $g_3=\left(\begin{smallmatrix}a_{11}&a_{12}\\a_{21}&a_{22}
\end{smallmatrix}\right)$. Then $a_{ij}$, $t_j$ and $u_j$, $i,j=1,2$
may be regarded as coordinates on $G_{k_v}/G^{\circ}_{x k_v}$.
With respect to these coordinates, the Jacobian determinant of the
map $g\mapsto gx$ is found to be $4|t_1^2t_2^2\det(g_3)^2|_v$.
Since the map is a double cover, $P(gx)=\chi(g)P(x)=
t_1^2t_2^2\det(g_3)^2P(x)$ and $P(x)=-4$ this shows that the
pullback of the measure $dy/|P(y)|_v^2$ to $G/G^{\circ}_{x\,k_v}$
is
\begin{equation*}
dg'_{x,v}=\frac{du_1\,dt_1\,du_2\,dt_2}{8|t_1^2t_2^2|_v}\,
d\mu_{3,v}(g_3)
\end{equation*}
where $d\mu_{3,v}(g_3)$ is the coordinate measure on the third
factor. An easy Jacobian calculation shows that if
$g=n(u)a(1,t)A_p(c,d)$ then the coordinate measure is 
$du\,dt\,dc\,dd/|t^2|_v(c^2+d^2)$
and so, with $d\mu_v$ denoting the product of the coordinate measures 
on the three factors, $d\mu_v=2\,dg'_{x,v}\,dg''_{x,v}$. Since
$dg_v=p_v^3\,d\mu_v(g)=2p_v^3\,dg'_{x,v}\,dg''_{x,v}$, it follows
that $b_{x,v}=2p_v^3$. \end{proof}
\begin{prop}\label{bxinf-(rm_rm)*}
Suppose that $v\in\mathfrak{M}_{\mathbb{R}}$ and that $x$ is a
representative for the orbit of type \upshape{(rm rm)*}. Then
\begin{equation*}
b_{x,v}=\tfrac{1}{\pi^2}.
\end{equation*}
\end{prop}
\begin{proof}
We again use $x=w_p$, where $p(z)=z^2+1$, as the orbital
representative. We let $\tl{v}$ be the complex place of $\kt$
dividing $v$. With $A_p(c,d)$ as in (\ref{apdefn})
and
\begin{equation*}
c_2 = \frac{|c_1|_{\tl{v}}-|d_1|_{\tl{v}}}{|c_1^2+d_1^2|_{\tl{v}}},\; 
d_2 = \frac{-(c_1\bar{d}_1+\bar{c}_1d_1)}{|c_1^2+d_1^2|_{\tl{v}}}\,,
\end{equation*}
any element of $G_{x\,k_v}^{\circ}$ has the form
$g=(A_p(c_1,d_1),A_p(c_1,d_2))$ where $c_1^2+d_1^2\not= 0$.
The isomorphism $\theta$ from $G_{x\,k_v}^{\circ}$ to
$H_{x\,k_v}=(\C^{\times})^2$ may be taken as
\begin{equation*}
\theta(A_p(c_1,d_1),A_p(c_2,d_2))=(c_1+d_1\sqrt{-1},c_1-d_1\sqrt{-1})\,.
\end{equation*}
From this it follows that
\begin{equation*}
dg''_{x,v}=\frac{4\,dc_1\,dd_1}{|c_1^2+d_1^2|_{\tl{v}}}\,.
\end{equation*}
It was shown during the proof of
Lemma (\ref{lemma:gamma_set}) that any matrix
$\left(\begin{smallmatrix}m_{11}&m_{12}\\m_{21}&m_{22}
\end{smallmatrix}\right)$ in $\gl(2)_{\mathbb{C}}$ may be written
in the form $n(u)a(1,t)A_p(c_1,d_1)$ provided that
$m_{11}^2+m_{12}^2\neq 0$. Since the complement of the set of
matrices satisfying this condition has measure zero it suffices to
make the comparison of measures on the set of elements of $G_{k_v}$
whose first entry satisfies this condition. Any $g$ in this set may
be expressed as
\begin{equation*}
g=(n(u)a(1,t)A_p(c_1,d_1),g_2A_p(c_2,d_2))
\end{equation*}
where $g_2=\left(\begin{smallmatrix}a_{11}&a_{12}\\a_{21}&a_{22}
\end{smallmatrix}\right)$. We may use $a_{ij}$, $i,j=1,2$,
$\re(u)$, $\im(u)$, $\re(t)$ and $\im(t)$ as real coordinates on (a
set of comeasure zero in) $G_{k_v}/G_{x k_v}^{\circ}$. With respect to
these coordinates the Jacobian determinant of the map $g\mapsto
gx$ is $16|t|^2_{\tl{v}}|\det(g_2)|_v^2$. This shows, as
usual, that the pullback of the measure $dy/|P(y)|_v^2$ to
$G_{k_v}/G^{\circ}_{x k_v}$ is
\begin{equation*}
dg_{x,v}'=\frac{du\,dt}{8|t|_{\tl{v}}^2}\,d\mu_{2,v}(g_2)
\end{equation*}
where $d\mu_{2,v}(g_2)$ is the coordinate measure on the second
factor. (Note that $du = 2d\re(u) d\im(u)$ and $dt = 2d\re(t) d\im(t)$
again.)
An easy Jacobian calculation shows that if
$g=n(u)a(1,t)A_p(c_1,d_1)$ then the coordinate measure is
$du\,dt\,dc_1\,dd_1/|t|_{\tl{v}}^2|c_1^2+d_1^2|_{\tl{v}}$ and
so, with $d\mu_{v}$ denoting the product of the coordinate measures
on the two factors, $d\mu_{v}=2\,dg_{x,v}'\,dg_{x,v}''$. Since
$dg_v=p_{v}p_{\tl{v}}\,d\mu_{v}(g)=2p_{v}p_{\tl{v}}\,dg_{x,v}'
\,dg_{x,v}''$, it follows that $b_{x,v}=2p_vp_{\tl{v}}$. 
\end{proof}

Even though we made use of results obtained from the software package
MAPLE in the above proofs, everything in this section could have been
proved manually without undue difficulty.  For example, in the proof
of Proposition \ref{bxinf-(sp)}, the only place we used MAPLE was
to determine the Jacobian of the map $g\mapsto gx$.
Using the invariance
properties of the measures, one can easily prove that this is a constant
multiple of $|\det(g_3)|_v^2$. To determine this constant we only
have to compute the Jacobian of the above map at the identity element.
At this point the Jacobian matrix is a fairly sparse $8\times 8$ matrix
and its determinant is easily found by hand. Indeed, the value of this
constant was verified manually in this case.  
However, we chose not to include the details of these manual
calculations in this paper.

\bibliographystyle{plain}
\bibliography{paper1}

\end{document}